\documentclass[11pt]{article}
\usepackage{latexsym,amssymb}
\usepackage{mathtools}
\usepackage{psfrag}
\usepackage{amsmath}
\usepackage{color}
\usepackage{subfig}
\usepackage[abbrev, nobysame]{amsrefs}
\usepackage{ulem}

\usepackage{tikz}
\usetikzlibrary{intersections, calc, angles, arrows.meta}
\usetikzlibrary{patterns}

\newtheorem{Theorem}{\bf Theorem}[section]
\newtheorem{Lemma}{\bf Lemma}[section]
\newtheorem{Proposition}{\bf Proposition}[section]
\newtheorem{Corollary}{\bf Corollary}[section]
\newtheorem{Remark}{\bf Remark}[section]
\newtheorem{Example}{\bf Example}[section]
\newtheorem{Definition}{\bf Definition}[section]

\newenvironment{theorem}{\begin{Theorem}$\!\!\!$}{\end{Theorem}}
\newenvironment{lemma}{\begin{Lemma}$\!\!\!$}{\end{Lemma}}
\newenvironment{proposition}{\begin{Proposition}$\!\!\!$}{\end{Proposition}}
\newenvironment{corollary}{\begin{Corollary}$\!\!\!$}{\end{Corollary}}
\newenvironment{remark}{\begin{Remark}$\!\!\!$}{\end{Remark}}

\newenvironment{definition}{\begin{Definition}$\!\!\!$}{\end{Definition}}

\newcommand{\K}{\mathsf{K}}
\newcommand{\G}{\mathsf{G}}
\newcommand{\dee}{{\rm{d}}}

\def\XXint#1#2#3{{\setbox0=\hbox{$#1{#2#3}{\int}$}
\vcenter{\hbox{$#2#3$}}\kern-.5\wd0}}

\numberwithin{equation}{section}

\allowdisplaybreaks[4]

\begin{document}

\title{Initial traces of solutions to a semilinear heat equation\\ 
under the Dirichlet boundary condition}
\author{Kotaro Hisa and Kazuhiro Ishige\\ \\
Graduate School of Mathematical Sciences, The University of Tokyo,\\ 
3-8-1 Komaba, Meguro-ku, Tokyo 153-8914, Japan}
\date{}
\maketitle
\begin{abstract}
We study qualitative properties of initial traces of nonnegative solutions 
to a semilinear heat equation in a smooth domain 
under the Dirichlet boundary condition. Furthermore, 
for the corresponding Cauchy--Dirichlet problem, we obtain sharp necessary conditions and sufficient conditions
on the existence of nonnegative solutions and identify optimal singularities 
of solvable nonnegative initial data. 
\end{abstract}
\vspace{25pt}
\noindent E-mail addresses:\\
K. H.: {\tt hisak@ms.u-tokyo.ac.jp}\\
K. I. : {\tt ishige@ms.u-tokyo.ac.jp}
\vspace{25pt}

\noindent
{\it MSC:} 35K58, 
35A01, 35A21, 35K20 
\vspace{3pt}
%
%
\vspace{3pt}
\newpage
\section{Introduction}
Let $u$ be a solution to the semilinear heat equation under the Dirichlet boundary condition 
\begin{equation}
\tag{E}
\label{eq:E}
\left\{
\begin{array}{ll}
	\partial_t u=\Delta u+u^p  & \mbox{in}\quad\Omega\times(0,T), \vspace{3pt}\\
	u\ge 0 & \mbox{in}\quad\Omega\times(0,T), \vspace{3pt}\\
	u = 0 \quad  & \mbox{on}\quad \partial \Omega\times(0,T),
\end{array}
\right.
\end{equation}
where $\partial_t:=\partial/\partial t$, $N\ge 1$, $T\in(0,\infty]$, $p>1$, 
and $\Omega\subset{\mathbb R}^N$ is a uniformly regular domain of class $C^2$ with non-empty boundary $\partial\Omega$.
In this paper we show the existence and the uniqueness of the initial trace of $d(\cdot)u(\cdot,t)$
and investigate its qualitative properties, in particular, the strength of singularities of the initial trace.
Here $d$ is the distance function from the boundary $\partial\Omega$, that is, 
$$
d(x):=\mbox{dist}\,(x,\partial\Omega)\equiv \inf_{z\in\partial\Omega}|x-z|\quad\mbox{for $x\in{\mathbb R}^N$}.
$$ 
Throughout this paper 
we denote by ${\mathcal M}(\overline{\Omega})$ (resp.~${\mathcal M}$) 
the set of nonnegative Radon measures on $\overline{\Omega}$ (resp.~${\mathbb R}^N$). 
We often identify $\dee \mu=\mu (x)\,\dee x$ in ${\mathcal M}(\overline{\Omega})$ for $\mu\in L^1_{\rm loc}(\overline{\Omega})$. 

In general, a nonnegative solution to the initial--boundary value problem for parabolic equations 
possesses its initial and lateral traces. 
In our case, thanks to the weight $d$, 
the lateral trace of $d(\cdot)u(\cdot,t)$ is identically zero at $t=0$, 
and we can characterize the initial trace of $d(\cdot)u(\cdot,t)$ in the framework of ${\mathcal M}(\overline{\Omega})$. 
The study of the strength of singularities of the initial trace 
enables us to obtain sharp necessary conditions for the existence of solutions to the Cauchy--Dirichlet problem
\begin{equation}
\tag{P}
  \label{eq:PP}
\left\{
\begin{array}{ll}
	\partial_t u=\Delta u+u^p  & \mbox{in}\quad\Omega\times(0,T), \vspace{3pt}\\
	u = 0 \quad  & \mbox{on}\quad \partial \Omega\times(0,T), \vspace{3pt}\\
	d u=\mu & \mbox{on}\quad\overline{\Omega}\times\{0\},
\end{array}
\right.
\end{equation}
where $\mu\in{\mathcal M}(\overline{\Omega})$. 

The study of initial and lateral traces of nonnegative solutions to parabolic equations is a classical subject 
and it has fascinated many mathematicians since the papers~\cites{HW, W1,W2} for the one dimensional heat equation. 
Qualitative properties of initial traces of nonnegative solutions strongly depend on the diffusion term,  the nonlinear term, and the boundary condition. 
Concerning relations among the initial trace and the diffusion term, 
see~e.g., \cites{A, HC} for linear parabolic equations; 
see e.g.,~\cites{AC, BCP, DK, HP} for porous medium equations; 
see e.g.,~\cites{DH, DH02} for parabolic $p$\,-Laplace equations;
see e.g.,~\cites{I, IJK, ZX} for doubly nonlinear parabolic equations; 
see e.g.,~\cite{BSV} for fractional diffusion equations; 
see e.g.,~\cite{AIS} for Finsler heat equations.
Concerning relations among the initial trace and the nonlinear term, 
see e.g.,~\cites{ADi, FHIL, HIT02, FI01, Hisa, HI02, IKO, TY} for source nonlinearity (positive nonlinearity);
see e.g.,~\cites{BCV, BD, MV, MV02, MV03} for absorption nonlinearity (negative nonlinearity). 
Notice that qualitative properties of initial traces in the case of source nonlinearity are quite different from those in the case of absorption nonlinearity.
\vspace{5pt}

We recall some results of \cite{HI01} on initial traces of solutions to the semilinear heat equation, 
which belongs to the case of source nonlinearity, 
\begin{equation}
\tag{E'}
\label{eq:E'}
\left\{
\begin{array}{ll}
\partial_t u=\Delta u+u^p\quad & \mbox{in}\quad {\mathbb R}^N\times(0,T),\vspace{3pt}\\
u\ge 0\quad & \mbox{in}\quad {\mathbb R}^N\times(0,T),
\end{array}
\right.
\end{equation}
where $p>1$ and $T\in(0,\infty]$. 
Let $p_k:=1+2/k$ for $k=1,2,\dots$.  
\begin{itemize}
  \item[(A)] 
  Assume that problem~\eqref{eq:E'} possesses a solution~$u$ in ${\mathbb R}^N\times(0,T)$ for some $T\in(0,\infty)$. 
  Then there exists a unique $\nu\in{\mathcal M}$ such that 
  $$
  \underset{\tau\to+0}{\mbox{{\rm ess lim}}}\int_{{\mathbb R}^N} \psi(y)u(y,\tau)\,\dee y=\int_{{\mathbb R}^N}\psi(y)\,\dee \nu(y),\quad \psi\in C_0({\mathbb R}^N).
  $$
  Furthermore,  
  there exists $C_1=C_1(N,p)>0$ such that
  $$
  \nu(B(z,\sigma))\le
  \left\{
  \begin{array}{ll}
  C_1 \sigma^{N-\frac{2}{p-1}} & \mbox{if}\quad p\not=p_N,\vspace{5pt}\\
  C_1\displaystyle{\left[\log\left(e+\frac{\sqrt{T}}{\sigma}\right)\right]^{-\frac{N}{2}}} & \mbox{if}\quad p=p_N,
  \end{array}
  \right.
  $$
  for $z\in{\mathbb R}^N$ and $\sigma\in(0,\sqrt{T})$. 
  Here $B(z,\sigma):=\{y\in{\mathbb R}^N\,:\,|y-z|<\sigma\}$.
   \item[(A')] 
  Let $u$ be a solution to the Cauchy problem
  \begin{equation}
  \tag{P'}
  \label{eq:P'}
  \partial_t u=\Delta u+u^p\quad\mbox{in}\quad {\mathbb R}^N\times(0,T),
  \qquad
  u=\mu\quad\mbox{in}\quad{\mathbb R}^N\times\{0\},
  \end{equation}
  where $\mu\in{\mathcal M}$ and $T\in(0,\infty]$. 
  Then $u$ is a solution to problem~\eqref{eq:E'} and the initial trace of $u$ coincides with the initial data $\mu$
  in ${\mathcal M}$. 
\end{itemize}
See \cite{HI01}*{Theorem~1.1} for assertion~(A) and \cite{HI01}*{Lemma~2.4} for assertion~(A'). 
(See e.g.,~\cites{ADi, BP, IKO, FHIL, FHIL02, FI01} for related results.) 
Thanks to assertion~(A'), assertion~(A) gives necessary conditions on the existence of solutions to problem~\eqref{eq:P'}. 
Notice that necessary conditions given in assertion~(A) are sharp (see e.g., \cite{HI01}*{Theorems 1.3, 1.4, and 1.5}). 

Next, we consider the case when $\Omega$ is a half-space of ${\mathbb R}^N$. 
Then the initial traces of solutions to problem~\eqref{eq:E} possess different properties of those in problem~\eqref{eq:E'}. 
Indeed, in~\cite{HIT02}, the authors of this paper and Takahashi 
developed the arguments in \cite{HI01} to obtain the following results.
\begin{itemize}
  \item[(B)] 
  Let $\Omega={\mathbb R}^N_+:=\{(x',x_N)\in{\mathbb R}^{N-1}\times{\mathbb R}\,:\,x_N>0\}$. 
  Assume that problem~\eqref{eq:E} possesses a solution~$u$ in $\Omega\times(0,T)$ for some $T\in(0,\infty)$. 
  Then there exists a unique $\nu\in{\mathcal M}(\overline{\Omega})$ such that 
  $$
  \underset{\tau\to+0}{\mbox{{\rm ess lim}}}\int_\Omega \psi(y) y_Nu(y,\tau)\,\dee y=\int_{\overline{\Omega}}\psi(y)\,\dee \nu(y),\quad \psi\in C_0({\mathbb R}^N),
  $$
  and $\nu(\partial\Omega)= 0$ if $p\ge 2$. 
  Furthermore, 
  there exists $C_2=C_2(N,p)>0$ such that
  \begin{equation}
  \label{eq:1.1}
  \nu(B_\Omega(z,\sigma))\leq C_2
  \sigma^{-\frac{2}{p-1}}\int_{B_\Omega(z,\sigma)}y_N\,\dee y
  \end{equation}
  for $z\in\overline{\Omega}$ and $\sigma\in(0,\sqrt{T})$. 
  Here $B_\Omega(z,\sigma):=\overline{\Omega}\cap B(z,\sigma)$.
  In addition, 
  \begin{itemize}
  \item[{\rm (1)}]
  if $p=p_N$, then
  there exists $C_3=C_3(N)>0$ such that 
  $$
  z_N^{-1}\nu(B_\Omega(z,\sigma))\le C_3 \left[\log\left(e+\frac{\sqrt{T}}{\sigma}\right)\right]^{-\frac{N}{2}}
  $$
  for $z=(z',z_N)\in\Omega$ with $z_N\ge 3\sqrt{T}$ and $\sigma\in(0,\sqrt{T})$;
  \item[{\rm (2)}]
  if $p=p_{N+1}$, then there exists $C_4=C_4(N)>0$ such that
  $$
  \nu(B_\Omega(z,\sigma))\le C_4
  \left[\log\left(e+\frac{\sqrt{T}}{\sigma}\right)\right]^{-\frac{N+1}{2}}
 $$
  for $z\in\partial\Omega$ and $\sigma\in(0,\sqrt{T})$.
  \end{itemize}
  \item[(B')] 
  Let $\Omega={\mathbb R}^N_+$. 
  Let $u$ be a solution to problem~\eqref{eq:PP} in $\Omega\times(0,T)$, where $T\in(0,\infty]$. 
  Then $u$ is a solution to problem~\eqref{eq:E} and the initial trace of $u$ coincides with the initial data of~$u$
  in ${\mathcal M}(\overline{\Omega})$. 
\end{itemize}
Similarly to problem~\eqref{eq:P'}, 
thanks to assertion~(B'), assertion~(B) gives necessary conditions on the existence of solutions to problem~\eqref{eq:PP} 
with $\Omega={\mathbb R}^N_+$. 
Notice that necessary conditions given in assertion~(B) are sharp (see \cite{HIT02}*{Section~5}). 
Recently, the first author of this paper obtained similar results for the fractional semilinear heat equation in domains (see \cite{Hisa}), 
however the arguments in \cite{Hisa} are not applicable to our problems~\eqref{eq:E} and \eqref{eq:PP}. 

In this paper we improve the results in assertions~(B) and (B') 
and generalize them to results in general uniformly regular domains of class $C^2$. 
Furthermore, we obtain a characterization of solvable initial data of problem~\eqref{eq:PP}. 
The results of this paper give new properties of the initial traces of solutions to problem~\eqref{eq:E} 
even in the case of $\Omega={\mathbb R}^N_+$. 
\subsection{Notation and definition of solutions}
We introduce some notation and formulate the definitions of solutions to problems~\eqref{eq:E} and \eqref{eq:PP}. 
For any measurable set $E$ in ${\mathbb R}^n$ and $n = 1,2,\cdots$, 
we denote by $|E|_n$ the $n$-dimensional Lebesgue measure of $E$. 
Let $\langle\cdot,\cdot\rangle$ be the inner product in ${\mathbb R}^N$. 
Throughout this paper, let $\Omega\subset{\mathbb R}^N$ be a uniformly regular domain of class $C^2$ such that $\partial\Omega\not=\emptyset$. 
(See \cite{QS}*{Section 1} for the definition of uniformly regular domains.) 
We denote by $\dee S$ the surface measure on $\partial\Omega$. 
Set
\begin{align*}
 & B(x,r):=\{y\in{\mathbb R}^N\,:\,|y-x|<r\},\quad 
 B_\Omega(x,r):=\overline{\Omega}\cap B(x,r),\\
 & \Omega(r):=\{x\in\overline{\Omega}\,:\, d(x)\le r\},\quad
 \quad
 \Omega^c(r):=\{x\in\overline{\Omega}\,:\, d(x)>r\},
\end{align*}
for $x\in{\mathbb R}^N$ and $r>0$. 
For any $T\in(0,\infty]$, we write
$Q_T:=\Omega\times(0,T)$. 

For any positive functions $f$ and $g$ in a set $X$, 
we say that $f\preceq g$ for $x\in X$ or equivalently that $f\succeq g$ for $x\in X$ if 
there exists $C>0$ such that 
$$
f(x)\le Cg(x),\quad x\in X.
$$
If $f\preceq g$ and  $g\preceq f$ for $x\in X$, we say that $f\asymp g$ for $x\in X$.
In all that follows we will use $C$ to denote generic positive constants and point out that $C$  
may take different values  within a calculation. 

We denote by $G_\Omega$ the minimal Dirichlet heat kernel in~$\Omega$. 
Then, for any fixed $y\in\Omega$, $G_\Omega(\cdot,y,\cdot)\in C^{2;1}(\Omega\times(0,\infty))\cap C^{1;0}(\overline{\Omega}\times(0,\infty))$ 
and $G_\Omega$ satisfies 
\begin{equation}
\label{eq:1.2}
\left\{
\begin{array}{ll}
\partial_t G_\Omega(x,y,t)=\Delta_x G_\Omega(x,y,t),\quad & (x,y,t)\in\Omega\times\Omega\times(0,\infty),\vspace{5pt}\\
G_\Omega(x,y,t)>0,\quad & (x,y,t)\in\Omega\times\Omega\times(0,\infty),\vspace{5pt}\\
G_\Omega(x,y,t)=0,\quad & (x,y,t)\in\partial\Omega\times\Omega\times(0,\infty),\vspace{5pt}\\
G_\Omega(x,y,t)=G_\Omega(y,x,t),\quad  & (x,y,t)\in\Omega\times\Omega\times(0,\infty),\vspace{5pt}\\
G_\Omega(x,y,0)=\delta_N(x-y), & (x,y)\in\Omega\times\Omega,
\end{array}
\right.
\end{equation}
where $\delta_N$ is the Dirac delta function in ${\mathbb R}^N$ concentrated at the origin. 
Furthermore, 
\begin{equation}
\label{eq:1.3}
\sup_{y\in\Omega}\|(\nabla_x G_\Omega)(\cdot,y,t)\|_{L^\infty(\Omega)}<\infty,\quad t>0.
\end{equation}
See e.g., \cite{LSU}*{Chapter~IV} for $C^{2;1}$ interior regularity of $G_\Omega$ and see e.g., \cite{L} for $C^{1;0}$-boundary regularity of $G_\Omega$. 
In addition,  
\begin{equation}
\label{eq:1.4}
G_\Omega(x,y,t+s)=\int_\Omega G_\Omega(x,z,t)G_\Omega(z,y,s)\,\dee z
\end{equation}
for $x\in\overline{\Omega}$, $y\in \Omega$, and $t$, $s>0$.

We define the integral kernel $K_\Omega$ on $\Omega\times\overline{\Omega}\times(0,\infty)$ by 
\begin{equation}
\label{eq:1.5}
K_\Omega(x,y,t):=
\left\{
\begin{array}{ll}
\displaystyle{\frac{G_\Omega(x,y,t)}{d(y)}}\quad & \mbox{for}\quad (x,y,t)\in \Omega\times\Omega\times(0,\infty),\vspace{7pt}\\ 
\langle\nabla_yG_\Omega(x,y,t),n_y\rangle \quad & \mbox{for}\quad (x,y,t)\in \Omega\times\partial\Omega\times(0,\infty),
\end{array}
\right.
\end{equation}
where $n_y$ is the inner normal unit vector at $y\in \partial \Omega$. 
Since $G_\Omega(x,y,t)=G_\Omega(y,x,t)=0$ for $x\in\Omega$ and $y\in\partial\Omega$, 
we have
\begin{equation}
\label{eq:1.6}
\lim_{\Omega\ni y\to y_*}K_\Omega(x,y,t)=
\lim_{\Omega\ni y\to y_*}\frac{G_\Omega(x,y,t)}{d(y)}=\langle\nabla_yG_\Omega(x,y_*,t),n_{y_*}\rangle=K_\Omega(x,y_*,t)
\end{equation}
for $x\in\Omega$ and $y_*\in\partial\Omega$. 
Then $K_\Omega\in C(\Omega\times\overline{\Omega}\times(0,\infty))$ 
and Hopf's lemma implies that
\begin{equation}
\label{eq:1.7}
K_\Omega(x,y,t)>0,\quad (x,y,t)\in \Omega\times\overline{\Omega}\times(0,\infty).
\end{equation} 
In addition, it follows from \eqref{eq:1.3}, \eqref{eq:1.4}, \eqref{eq:1.5}, and \eqref{eq:1.6} that 
\begin{equation}
\label{eq:1.8}
K_\Omega(x,y,t+s)=\int_\Omega G_\Omega(x,z,t)K_\Omega(z,y,s)\,\dee z
\end{equation}
for $x\in\Omega$, $y\in\overline{\Omega}$, and $t$, $s>0$.
Then, thanks to \eqref{eq:1.3}, the domain of $K_\Omega$ can be extended into $\overline{\Omega}\times\overline{\Omega}\times(0,\infty)$ such that 
\begin{equation}
\label{eq:1.9}
K_\Omega\in C(\overline{\Omega}\times\overline{\Omega}\times(0,\infty))\quad\mbox{and}\quad
K_\Omega(x,y,t)=0,\quad (x,y,t)\in\partial\Omega\times\overline{\Omega}\times(0,\infty). 
\end{equation}
Furthermore, we observe from \eqref{eq:1.8} that 
\begin{equation}
\label{eq:1.10}
\partial_tK_\Omega(x,y,t)=\Delta_x K_\Omega(x,y,t),\quad (x,y,t)\in\Omega\times\overline{\Omega}\times(0,\infty). 
\end{equation}
We write $G=G_\Omega$ and $K=K_\Omega$ simply 
when no confusion can arise.  
\vspace{5pt}

We formulate the definitions of solutions to problems~\eqref{eq:E} and \eqref{eq:PP}.
\begin{definition}
\label{Definition:1.1}
Let $u$ be a nonnegative, measurable, and finite almost everywhere function in~$Q_T$, where $T\in(0,\infty]$. 
\begin{itemize}
  \item[{\rm (1)}] 
  We say that  $u$ is a solution to problem~\eqref{eq:E} in~$Q_T$
  if, for a.a.~$\tau\in(0,T)$, $u$ satisfies 
  $$
  u(x,t)= \int_{\Omega} G(x,y,t-\tau)u(y,\tau) \, \dee y
  +\int_\tau^t\int_{\Omega} G(x,y,t-s)u(y,s)^p  \,\dee y\,\dee s
  $$
  for a.a.~$(x,t)\in \Omega\times(\tau,T)$. 
  \item[{\rm (2)}]
  Let $\mu\in{\mathcal M}(\overline{\Omega})$. 
  We say that $u$ is a solution to problem~\eqref{eq:PP} in~$Q_T$
  if $u$ satisfies 
  \begin{equation}
  \label{eq:1.11}
  u(x,t)= \int_{\overline{\Omega}} K(x,y,t)\, \dee \mu(y)
  +\int_0^t\int_{\Omega} G(x,y,t-s)u(y,s)^p  \,\dee y\,\dee s
  \end{equation}
  for a.a.~$(x,t)\in Q_T$.   
  If $u$ satisfies the above relation with $``="$ replaced by $``\ge"$, 
  then we say that $u$ is a supersolution to problem~\eqref{eq:PP} in~$Q_T$. 
\end{itemize}
\end{definition}
\subsection{Main results}
We state our main results on initial traces of solutions to problem~\eqref{eq:E}. 
The first theorem concerns with the existence and the uniqueness of the initial trace of $d(\cdot)u(\cdot,t)$. 
\begin{theorem}
\label{Theorem:1.1}
Let $\Omega$ be a uniformly regular domain in ${\mathbb R}^N$ of class $C^2$ such that $\partial\Omega\not=\emptyset$.
\begin{itemize}
  \item[{\rm (1)}] 
  Assume that 
  problem~\eqref{eq:E} possesses a solution~$u$ in~$Q_T$ for some $T\in(0,\infty]$. 
  Then there exists a unique initial trace $\nu\in{\mathcal M}(\overline{\Omega})$ of $d(\cdot,)u(\cdot,t)$, that is,   
  $$
  \underset{t\to +0}{\mbox{{\rm ess lim}}}
  \int_\Omega \psi(y)d(y)u(y,t)\,\dee y=\int_{\overline{\Omega}} \psi(y)\,\dee \nu(y), 
  \quad \psi\in C_0({\mathbb R}^N).
  $$
  \item[{\rm (2)}]  
  Let $u$ be a solution to problem~\eqref{eq:PP} in~$Q_T$, where $T\in(0,\infty]$. 
  Then $u$ is a solution to problem~{\rm (E)} in~$Q_T$ and it satisfies assertion~{\rm (1)} with $\nu=\mu$.
\end{itemize}
\end{theorem}
Next, we state a result on estimates of the initial trace given in Theorem~\ref{Theorem:1.1}. 
\begin{theorem}
\label{Theorem:1.2}
Let $\Omega$ be a uniformly regular domain in ${\mathbb R}^N$ of class $C^2$  such that $\partial\Omega\not=\emptyset$.
Let $R_*$ be a positive constant satisfying the following condition:
\begin{itemize}
  \item[{\rm ({\bf D})}] 
  \begin{itemize}
  \item[{\rm (1)}] 
  for any $z\in \partial\Omega$ and $R\in(0,R_*)$, 
  $$
  z_R:=z+Rn_z\in\Omega,\quad 
  d(z_R)=R,\quad \partial B(z_R,R)\cap\partial\Omega=\{z\};
  $$
  \item[{\rm (2)}]
  there exists $C_G>0$ such that 
  \begin{equation*}
  G(x,y,t)\ge C_G^{-1}t^{-\frac{N}{2}}\frac{d(x)}{d(x)+\sqrt{t}}\frac{d(y)}{d(y)+\sqrt{t}}\exp\left(-C_G\frac{|x-y|^2}{t}\right)
  \end{equation*}
  for $(x,y,t)\in \Omega\times \Omega\times(0,R_*^2)$.
  \end{itemize}
\end{itemize}
Assume that problem~\eqref{eq:E} possesses a solution $u$ in~$Q_T$ for some $T\in(0,R_*^2)$. 
Then the initial trace $\nu$ of $d(\cdot)u(\cdot,t)$ satisfies the following property:
\begin{itemize}
  \item[{\rm (1)}] 
  there exists $C_1=C_1(\Omega,p)>0$ such that
  $$
  \nu(B_\Omega(z,\sigma))\le
  C_1\inf_{s\in[\sigma,\sqrt{T})}\left\{(d(z)+s) s^{N-\frac{2}{p-1}}\right\}
 $$
  for $z\in\overline{\Omega}$ and $\sigma\in(0,\sqrt{T})$. 
  \end{itemize}
Furthermore, 
\begin{itemize}
  \item[{\rm (2)}] 
  if $p=p_N$, then there exists $C_2=C_2(\Omega)>0$ such that 
  $$
  \nu(B_\Omega(z,\sigma))\le
  C_2(d(z)+\sigma)\left[\log\left(e+\frac{\min\{d(z),\sqrt{T}\}}{\sigma}\right)\right]^{-\frac{N}{2}}
  $$
  for $z\in\overline{\Omega}$ and $\sigma\in(0,\sqrt{T})$;
  \item[{\rm (3)}] 
  if $p=p_{N+1}$, then 
  there exists $C_3=C_3(\Omega)>0$ such that 
  $$
  \nu(B_\Omega(z,\sigma))\le C_3\left[\log\left(e+\frac{\sqrt{T}}{\sigma}\right)\right]^{-\frac{N+1}{2}}
  $$
  for $z\in\partial\Omega$ and $\sigma\in(0,\sqrt{T})$;
  \item[{\rm (4)}]  
 if $p\ge 2$, then $\nu(\partial\Omega)=0$.
\end{itemize}
\end{theorem}
\begin{remark}
\label{Remark:1.1}
We give some remarks on the constant $R_*>0$.
\begin{itemize}
  \item[{\rm (1)}] 
  Let $\Omega$ be a uniformly regular domain in ${\mathbb R}^N$ of class $C^2$ such that $\partial\Omega\not=\emptyset$. 
  Then condition~{\rm ({\bf D})} holds for some $R_*>0$. Indeed, condition~{\rm ({\bf D})-(1)} easily holds for some $R_*>0$.
  Furthermore, it follows from \cite{Cho}*{Theorem 1.1} that, 
  for any $T>0$, 
  there exist positive constants $c_1$, $c_2>0$ such that 
  \begin{equation}
  \label{eq:1.12}
  \begin{split}
   &c_1^{-1}t^{-\frac{N}{2}}\frac{d(x)}{d(x)+\sqrt{t}}\frac{d(y)}{d(y)+\sqrt{t}}
  \exp\left(-\frac{c_2|x-y|^2}{t}\right)\\
   & \qquad\qquad
   \le G(x,y,t)\le c_1t^{-\frac{N}{2}}\frac{d(x)}{d(x)+\sqrt{t}}\frac{d(y)}{d(y)+\sqrt{t}}
  \exp\left(-\frac{|x-y|^2}{c_2t}\right)
  \end{split}
  \end{equation}
  for $x$, $y\in\Omega$ and $t\in(0,T)$. 
  {\rm (}See e.g., \cites{CKP, R, Z} for related results on inequality~\eqref{eq:1.12}.{\rm )}
  This implies that condition~{\rm ({\bf D})-(2)} holds for some $R_*>0$.  
  \item[{\rm (2)}] 
  Let $\Omega={\mathbb R}^N_+$. Then condition~{\rm ({\bf D})} holds with $R_*=\infty$. 
  Indeed, condition~{\rm ({\bf D})-(1)} easily holds for all $R>0$. 
  Furthermore,  since $1-e^{-ab}\ge 1-e^{-\min\{a,1\}\min\{b,1\}}\succeq\min\{a,1\}\min\{b,1\}$ for $a$, $b\ge 0$, 
  we have 
  \begin{align*}
  G_{{\mathbb R}^N_+}(x,y,t) & =(4\pi t)^{-\frac{N}{2}}\exp\left(-\frac{|x-y|^2}{4t}\right)\left(1-\exp\biggr(-\frac{x_Ny_N}{t}\biggr)\right)\\
   & \ge C t^{-\frac{N}{2}}\frac{x_N}{x_N+\sqrt{t}}\frac{y_N}{y_N+\sqrt{t}}\exp\left(-\frac{|x-y|^2}{4t}\right)
  \end{align*}
  for $(x,y,t)\in {\mathbb R}^N_+\times {\mathbb R}^N_+\times(0,\infty)$, 
  that is, condition~{\rm ({\bf D})-(2)} holds with $R_*=\infty$.
 Then, since Theorem~{\rm\ref{Theorem:1.2}} holds with $R_*=\infty$,
  Theorem~{\rm\ref{Theorem:1.2}-(2)} gives a refined estimate in assertion~{\rm (B)-(1)}. 
\end{itemize}
\end{remark}
We state another estimate of the initial trace given in Theorem~\ref{Theorem:1.1}. 
\begin{theorem}
\label{Theorem:1.3}
Let $\Omega$ be a domain in ${\mathbb R}^N$ with non-empty, $C^2$-smooth, and compact boundary~$\partial\Omega$.
Let $\phi\in C^2(\overline{\Omega})$ be such that 
\begin{equation}
\label{eq:1.13}
\mbox{$\phi>0$ in $\Omega$},\quad
\mbox{$\phi=0$ on $\partial\Omega$},\quad
\sup_{z\in\Omega}\frac{|(\Delta \phi)(z)|+d(z)|(\nabla \phi)(z)|}{\phi(z)}<\infty.
\end{equation}
Let $R\in(0,\infty]$ be such that $d\in C^2(\Omega(R))$. 
If $R=\infty$, assume that 
\begin{equation}
\label{eq:1.14}
\sup_{z\in \Omega(R)} d(z)|\Delta d(z)|<\infty
\quad\mbox{and}\quad\mbox{$\Delta\phi=0$ in $\Omega$}.
\end{equation}
Then there exists $C_4>0$ such that, 
if problem~\eqref{eq:E} possesses a solution $u$ in~$Q_T$ for some $T\in(0,R^2]$, 
then the initial trace $\nu$ of $d(\cdot) u(\cdot,t)$ satisfies
\begin{equation}
\label{eq:1.15}
\int_{\Omega(\sigma)}\frac{\phi(y)}{d(y)}\,\dee \nu(y)\le  C_4
\left(\int_{2\sigma^2}^{T}\left(\int_{\Omega(\sqrt{r})}\phi(y)\,\dee y\right)^{-(p-1)}\,\dee r\right)^{-\frac{1}{p-1}}
\end{equation}
for $\sigma\in(0,\sqrt{T/2})$. 
\end{theorem}
The typical examples of~$\phi$ we mind  in Theorem~\ref{Theorem:1.3} are the first Dirichlet eigenfunction in $\Omega$ for $-\Delta$ or 
the positive harmonic function in $\Omega$. 
Theorem~\ref{Theorem:1.3} plays a crucial role in the proof of  Theorem~\ref{Theorem:1.2}-(4). 
Furthermore, as a corollary of Theorem~\ref{Theorem:1.3}, 
we obtain more precise behavior of the initial trace near the boundary~$\partial\Omega$ than that of Theorem~\ref{Theorem:1.2}-(4).
\begin{corollary}
\label{Corollary:1.1}
Let $\Omega$ be a uniformly regular domain in ${\mathbb R}^N$ of class  $C^{2,\theta}$ for some $\theta\in(0,1)$
such that $\partial\Omega\not=\emptyset$ 
and let $p\ge 2$. 
Then there exist $C_5=C_5(\Omega,p)>0$ and $R_0>0$ such that, 
if problem~\eqref{eq:E} possesses a solution~$u$ in $Q_T$ for some $T>0$, 
then the initial trace $\nu$ of $d(\cdot)u(\cdot,t)$ satisfies 
  $$
  \nu(B(z,R_0)\cap \Omega(\sigma))\le 
  \left\{
  \begin{array}{ll}
  C_5\sigma^{2\frac{p-2}{p-1}}  & \mbox{if}\quad p>2,\vspace{7pt}\\
C_5\displaystyle{\left[\log\left(e+\frac{\sqrt{T}}{\sigma}\right)\right]^{-1}} & \mbox{if}\quad p=2,\\
  \end{array}
  \right.
  $$
for $z\in\partial\Omega$ and $\sigma\in(0,\sqrt{T/2})$. 
In particular, $\nu(\partial\Omega)=0$ if $p\ge 2$.
\end{corollary}
Theorem~\ref{Theorem:1.3} is also useful for the proof of the nonexistence of non-trivial global-in-time solutions to problem~\eqref{eq:PP}.
Indeed, the following result easily follows from Theorem~\ref{Theorem:1.3}.
\begin{corollary}
\label{Corollary:1.2}
Let $\Omega$ be the exterior of $C^2$-smooth compact set in ${\mathbb R}^N$, where $N\ge 2$. 
Then problem~\eqref{eq:PP} with $p\ge p_N$ possesses no non-trivial global-in-time solutions. 
\end{corollary}
Corollary~\ref{Corollary:1.2} has already been proved in \cite{IS02}, while 
the proof in this paper supports the sharpness of Theorem~\ref{Theorem:1.3}. 
\vspace{5pt}

Thanks to Theorem~\ref{Theorem:1.1}-(2), 
by Theorem~\ref{Theorem:1.2},  Theorem~\ref{Theorem:1.3}, and Corollary~\ref{Corollary:1.1}
we obtain necessary conditions on solvable initial data of problem~\eqref{eq:PP}. 
In Section~5 we also study sufficient conditions on the existence of solutions to problem~\eqref{eq:PP}. 
Combining our necessary conditions and sufficient conditions, 
we have:
\begin{theorem}
\label{Theorem:1.4}
Let $\Omega$ be a uniformly regular domain in ${\mathbb R}^N$ of class $C^2$ such that $\partial\Omega\not=\emptyset$.
\begin{itemize}
  \item[{\rm (1)}] 
  Let $1<p<p_{N+1}$. 
  Then problem~\eqref{eq:PP} possesses a local-in-time solution for initial data $\mu\in{\mathcal M}(\overline{\Omega})$ 
  if and only if 
  \begin{equation}
  \label{eq:1.16}
  \sup_{z\in\overline{\Omega}} \frac{\mu(B_\Omega(z,1))}{1+d(z)} <\infty.
  \end{equation}
 \item[{\rm (2)}]
 Let $\mu\in{\mathcal M}(\overline{\Omega})$ be such that $\rm{dist}\,(\rm{supp}\,\mu,\partial\Omega)>0$. 
 Then assertion~{\rm (1)} holds for $1<p<p_N$.  
 \item[{\rm (3)}]
 If problem~\eqref{eq:PP} possesses a local-in-time solution for initial data $\mu\in{\mathcal M}(\overline{\Omega})$ 
 and $p\ge 2$, then $\mu(\partial\Omega)=0$. 
\end{itemize}
\end{theorem}
\begin{theorem}
\label{Theorem:1.5}
Let  $\Omega$ be a uniformly regular domain in ${\mathbb R}^N$ of class $C^2$ such that $\partial\Omega\not=\emptyset$.
\begin{itemize}
  \item[{\rm (1)}] 
  Let $z\in\Omega$ and $p\ge p_N$. Let $\mu_1\in{\mathcal M}(\overline{\Omega})$ be such that 
  $$
  \dee \mu_1=
  \left\{
  \begin{array}{ll}
  d(x)|x-z|^{-\frac{2}{p-1}}\chi_{B_{\Omega}(z,1)}\,\dee x & \mbox{if}\quad p>p_N,\vspace{5pt}\\
  \displaystyle{d(x)|x-z|^{-N}\left[\log\left(e+\frac{1}{|x-z|}\right)\right]^{-\frac{N}{2}-1}\chi_{B_\Omega(z,1)}}\,\dee x & \mbox{if}\quad p=p_N.
  \end{array}
  \right.
  $$
  Then there exists $\kappa_1(z)>0$ such that 
  problem~\eqref{eq:PP} with $\mu=\kappa\mu_1$ possesses no local-in-time solutions if $\kappa>\kappa_1(z)$ 
  and possesses a local-in-time solution if $0<\kappa<\kappa_1(z)$. Here $\sup_{z\in\Omega}\kappa_1(z)<\infty$.
  \item[{\rm (2)}] 
  Let $z\in\partial\Omega$ and $p\ge p_{N+1}$. Let $\mu_2\in{\mathcal M}(\overline{\Omega})$ be such that 
  $$
  \dee \mu_2=
  \left\{
  \begin{array}{ll}
  d(x)|x-z|^{-\frac{2}{p-1}}\chi_{B_{\Omega}(z,1)}\,\dee x & \mbox{if}\quad p>p_{N+1},\vspace{5pt}\\
  d(x)\displaystyle{|x-z|^{-N-1}\left[\log\left(e+\frac{1}{|x-z|}\right)\right]^{-\frac{N+1}{2}-1}\chi_{B_\Omega(z,1)}}\,\dee x & \mbox{if}\quad p=p_{N+1}.
  \end{array}
  \right.
  $$
  Then there exists $\kappa_2(z)>0$ such that 
  problem~\eqref{eq:PP} with $\mu=\kappa\mu_2$ possesses no local-in-time solutions if $\kappa>\kappa_2(z)$ 
  and possesses a local-in-time solution if $0<\kappa<\kappa_2(z)$. 
  Here $\sup_{z\in\partial\Omega}\kappa_2(z)<\infty$.
  \item[{\rm (3)}] 
  Let $z\in\partial\Omega$  and $p_{N+1}\le p<2$. Let $\mu_3\in{\mathcal M}(\overline{\Omega})$ be such that 
  $$
  \dee \mu_3=
  \left\{
  \begin{array}{ll}
  |x-z|^{-\frac{2(2-p)}{p-1}}\chi_{B_{\Omega}(z,1)}\,\dee S & \mbox{if}\quad p>p_{N+1},\vspace{5pt}\\
  \displaystyle{|x-z|^{-N+1}\left[\log\left(e+\frac{1}{|x-z|}\right)\right]^{-\frac{N+1}{2}-1}\chi_{B_\Omega(z,1)}}\,\dee S & \mbox{if}\quad p=p_{N+1}.
  \end{array}
  \right.
  $$
  Then there exists $\kappa_3(z)>0$ such that 
  problem~\eqref{eq:PP} with $\mu=\kappa\mu_3$ possesses no local-in-time solutions if $\kappa>\kappa_3(z)$ 
  and possesses a local-in-time solution if $0<\kappa<\kappa_3(z)$. 
  Here $\sup_{z\in\partial\Omega}\kappa_3(z)<\infty$.
\end{itemize}
\end{theorem}
We observe from Theorems~\ref{Theorem:1.4} and \ref{Theorem:1.5} that
the solvability of problem~\eqref{eq:PP} varies depending on 
whether $1<p<p_{N+1}$, $p=p_{N+1}$, $p_{N+1}<p<p_N$, $p=p_N$, $p>p_N$, or $p\ge 2$. 
Furthermore, Theorem~\ref{Theorem:1.5} enables us to identify the optimal singularity of the solvable initial data 
of problem~\eqref{eq:PP}. 
\vspace{5pt}

This paper is a generalization and a refinement of \cite{HIT02}. 
We modify arguments in \cites{HIT02} to prove Theorem~\ref{Theorem:1.1}.
On the other hand, it seems difficult to apply the arguments in \cites{HIT02} 
to obtain estimates of initial traces of solutions to problem~\eqref{eq:E}. 
Indeed, in \cites{HIT02} the authors employed the explicit representation of the Dirichlet heat kernel in ${\mathbb R}^N_+$, 
and found $C>0$ such that 
\begin{equation}
\label{eq:1.17}
G_{{\mathbb R}^N_+}(x,y,2t-s)\ge C\left(\frac{s}{2t}\right)^{\frac{N}{2}}G_{{\mathbb R}^N_+}(x,y,s)
\end{equation}
holds for $x\in\Omega(\sqrt{T})$, $y\in\Omega$, and $s,t\in(0,T/C)$ with $s<t$ (see \cite{HIT02}*{Lemma~3.2}). 
This inequality plays a crucial role in obtaining interior and boundary estimates of initial traces in their arguments. 
However, due to the exponential decay of $G=G(x,y,t)$ as $t\to+0$ for $x\not=y$, 
one cannot expect inequality~\eqref{eq:1.17} to hold for general domains.
Notice that, in the fractional case, 
the Dirichlet fractional heat kernel has polynomial decay as $t\to +0$ and 
we can obtain similar estimates to that of \eqref{eq:1.17} for the Dirichlet fractional heat kernel in domains 
(see \cite{Hisa}*{Lemma~3.2}). Then arguments in \cite{HIT02} are applicable to the fractional semilinear heat equation in domains (see \cite{Hisa}).

In this paper, applying arguments in \cite{IKO} (see also \cites{IS, MP}), 
we find suitable cut-off functions to obtain interior estimates of solvable initial data of problem~\eqref{eq:PP} (see Lemma~\ref{Lemma:3.2}). 
Then, thanks to condition~({\bf D}), we obtain interior and boundary estimates of initial traces of solutions to problem~\eqref{eq:E} 
(see Theorem~\ref{Theorem:1.2}-(1), (2)). 
We also modify these arguments to prove Theorem~\ref{Theorem:1.3}, which leads to Theorem~\ref{Theorem:1.2}-(4). 
On the other hand, 
the proof of boundary estimates of initial traces in the case of $p=p_{N+1}$ (see Theorem~\ref{Theorem:1.2}-(3)) is delicate. 
In the proof, we modify the arguments in the proof of Lemma~\ref{Lemma:3.2} while we slide the supports of cut-off functions in~$Q_T$. 
Furthermore, we apply lower estimates of solutions to establish boundary estimates of initial traces in the case of $p=p_{N+1}$ (see Lemma~\ref{Lemma:4.1}). 
This enables us to obtain Theorem~\ref{Theorem:1.2}-(3), and Theorem~\ref{Theorem:1.2} follows. 

We also develop arguments in \cite{HIT02} to obtain sufficient conditions for the existence of solutions to problem~\eqref{eq:PP}. 
Then, thanks to Theorems~\ref{Theorem:1.1} and \ref{Theorem:1.2}, 
we obtain Theorems~\ref{Theorem:1.4} and \ref{Theorem:1.5}, which identify the optimal singularity 
of solvable initial data of problem~\eqref{eq:PP}. 
The proofs of Theorems~\ref{Theorem:1.4} and \ref{Theorem:1.5} show that 
our estimates of initial traces of solutions to problem~\eqref{eq:E}
and sufficient conditions in Section~5 are sharp. 

The rest of this paper is organized as follows. 
In Section~2 we modify arguments in \cites{HI01, HIT02} and prove some preliminary lemmas. 
Then we prove Theorem~\ref{Theorem:1.1}. 
In Section~3 we obtain interior and boundary estimates of solvable initial data of problem~\eqref{eq:PP}. 
Furthermore, we modify the arguments to prove Theorem~\ref{Theorem:1.3}. 
Then we obtain  Theorem~\ref{Theorem:1.2}-(1), (2), (4), Corollary~\ref{Corollary:1.1}, and Corollary~\ref{Corollary:1.2}. 
In Section~4 we prove Theorem~\ref{Theorem:1.2}-(3), and complete the proof of Theorem~\ref{Theorem:1.2}.
In Section~5 we obtain sufficient conditions 
on the existence of solutions to problem~\eqref{eq:PP}, and prove Theorem~\ref{Theorem:1.5}. 
\section{Proof of Theorem~\ref{Theorem:1.1}}
In this section we prove some preliminary lemmas on solutions to problems~\eqref{eq:E} and \eqref{eq:PP}. 
Furthermore, we prove Theorem~\ref{Theorem:1.1}.  
We first show two preliminary lemmas on solutions to problem~\eqref{eq:PP}. 
\begin{lemma}
\label{Lemma:2.1}
Assume that there exists a supersolution $v$ to problem \eqref{eq:PP} in~$Q_T$ for some $T\in(0,\infty]$.
Then problem~\eqref{eq:PP} possesses a solution $u$ in~$Q_T$ such that $u\le v$ almost everywhere in~$Q_T$.
\end{lemma}
{\bf Proof.}
This lemma can be proved by the same argument as in \cite{HI01}*{Lemma 2.2}.
Define 
\begin{align*}
& u_1(x,t) := \int_{\overline{\Omega}} K(x,y,t) \, \dee \mu (y),\\
& u_{j+1}(x,t):= u_1(x,t) + \int_0^t\int_{\Omega} G(x,y,t-s)u_j(y,s)^p\,\dee y\,\dee s, \quad j=1,2,\cdots,
\end{align*}
for a.a.~$(x,t)\in Q_T$.
Thanks to Definition~\ref{Definition:1.1}-(2) and the nonnegativity of integral kernels $K$ and $G$, by induction we obtain
\[
0\le u_1(x,t)\le u_2(x,t)\le\cdots\le u_j(x,t)\le \cdots \le v(x,t)<\infty
\]
for a.a.~$(x,t)\in Q_T$. Then the limit function
$u(x,t) := \lim_{j\to\infty}u_j(x,t)$
is well-defined for a.a.~$(x,t)\in Q_T$ and it satisfies integral equation~\eqref{eq:1.11} for a.a.~$(x,t)\in Q_T$ 
such that $u(x,t)\le v(x,t)<\infty$ for a.a.~$(x,t)\in Q_T$. 
Then $u$ is the desired solution to problem~\eqref{eq:PP} in~$Q_T$, 
and the proof is complete. 
$\Box$
\begin{lemma}
\label{Lemma:2.2}
{\rm (1)}
  Let $u$ be a solution to problem~\eqref{eq:E} in~$Q_T$, where $T\in(0,\infty]$. 
  Then there exists a measurable set $I\subset(0,T)$ with $|(0,T)\setminus I|_1=0$ such that, 
  for any $t_1$, $t_2\in I$ with $t_1<t_2$, 
  \begin{equation}
  \label{eq:2.1}
  u(x,t_2)
  =\int_\Omega G(x,y,t_2-t_1)u(y,t_1)\,\dee y+\int_{t_1}^{t_2}\int_\Omega G(x,y,t_2-s)u(y,s)^p\,\dee y\,\dee s
  \end{equation}
 holds for a.a.~$x\in\Omega$. 
  Furthermore, 
  \begin{equation}
  \label{eq:2.2}
  \sup_{t\in I\cap (0,T')}\int_{B_\Omega(z,R)}d(y)u(y,t)\,\dee y<\infty,
  \quad
  \int_0^{T'} \int_{B_\Omega(z,R)} d(y) u(y,s)^p\,\dee y\,\dee s<\infty,
  \end{equation}
  for $z\in\overline{\Omega}$, $R>0$, and $T'\in(0,T)$. 
\vspace{5pt}
\newline
{\rm (2)}
  Let $u$ be a solution to problem~\eqref{eq:PP} in~$Q_T$, where $T\in(0,\infty]$. 
  Then there exists a measurable set $I'\subset[0,T)$ with $|[0,T)\setminus I'|_1=0$ and $0\in I'$ with the following properties:
  \begin{itemize}
  \item[{\rm (a)}] 
  for any $t_1$, $t_2\in I'$ with $t_1<t_2$, \eqref{eq:2.1} holds for a.a.~$x\in\Omega$;
  \item[{\rm (b)}] 
  $u$ is a solution to problem~\eqref{eq:E} in~$Q_T$ 
  and for any $\tau\in I'\setminus\{0\}$, the function $u_\tau$ defined by $u_\tau(x,t):=u(x,t+\tau)$
  is a solution to problem~\eqref{eq:PP} in $Q_{T-\tau}$ with the initial data $\mu=d(\cdot)u(\cdot,\tau)$.
\end{itemize}
\end{lemma}
{\bf Proof.}
We prove assertion~(1). 
Let $u$ be a solution to problem~\eqref{eq:E} in $Q_T$, where $T\in(0,\infty]$. 
Thanks to Definition~\ref{Definition:1.1}-(1) and Fubini's theorem, 
there exists a measurable set $I_0\subset(0,T)$ with $|(0,T)\setminus I_0|_1=0$ with the following property: 
for any $\tau\in I_0$, there exists a measurable set $I_\tau\subset(\tau,T)$ with $|(\tau,T)\setminus I_\tau|_1=0$ such that, 
for any $t\in I_\tau$, 
$$
\infty>u(x,t)=\int_\Omega G(x,y,t-\tau)u(y,\tau)\,\dee y
+\int_\tau^t\int_\Omega G(x,y,t-s)u(y,s)^p\,\dee y\,\dee s
$$
for a.a.~$x\in\Omega$. 
Then, for any $t_1$, $t_2\in I_\tau$ with $t_1<t_2$, 
it follows from Fubini's theorem and \eqref{eq:1.4} that
\begin{align*}
 & \int_\Omega G(x,y,t_2-t_1)u(y,t_1)\,\dee y\\
 & =\int_\Omega G(x,y,t_2-t_1)\left(\int_\Omega G(y,\xi,t_1-\tau)u(\xi,\tau)\,\dee \xi+\int_\tau^{t_1}\int_\Omega G(y,\xi,t_1-s)u(\xi,s)^p\,\dee \xi\,\dee s\right)\,\dee y\\
 & =\int_\Omega G(x,\xi,t_2-\tau)u(\xi,\tau)\,\dee \xi+\int_\tau^{t_1}\int_\Omega G(x,\xi,t_2-s)u(\xi,s)^p\,\dee \xi\,\dee s
\end{align*}
for a.a~$x\in\Omega$. 
This implies that
\begin{equation}
\label{eq:2.3}
\begin{split}
 & \int_\Omega G(x,y,t_2-t_1)u(y,t_1)\,\dee y+\int_{t_1}^{t_2}\int_\Omega G(x,y,t_2-s)u(y,s)^p\,\dee y\,\dee s\\
 & =\int_\Omega G(x,\xi,t_2-\tau)u(\xi,\tau)\,\dee \xi+\int_\tau^{t_2}\int_\Omega G(x,\xi,t_2-s)u(\xi,s)^p\,\dee \xi\,\dee s\\
 & =u(x,t_2)<\infty
 \end{split}
\end{equation}
for a.a.~$x\in\Omega$. 
Let $\{\tau_j\}_{j=1}^\infty$ be a non-increasing sequence in $I_0$ such that $\lim_{j\to\infty}\tau_j=0$. 
Set 
$$
I:=\bigcup_{k=1}^\infty\biggr(\bigcap_{j=k}^\infty I_{\tau_j}\biggr).
$$
Then $|(0,T)\setminus I|_1=0$. 
Furthermore, by \eqref{eq:2.3}, for any $t_1$, $t_2\in I$ with $t_1<t_2$, 
we have
\begin{equation*}
\infty>u(x,t_2)
=\int_\Omega G(x,y,t_2-t_1)u(y,t_1)\,\dee y+\int_{t_1}^{t_2}\int_\Omega G(x,y,t_2-t_1)u(y,s)^p\,\dee y\,\dee s
\end{equation*}
for a.a.~$x\in\Omega$. Thus \eqref{eq:2.1} holds. 
Furthermore, for any $T'\in(0,T)$, taking $t_2\in I\cap (T',T)$, 
we see that 
$$
\infty>u(x,t_2)=\int_\Omega G(x,y,t_2-t)u(y,t)\,\dee y+\int_t^{t_2}\int_\Omega G(x,y,t_2-t)u(y,s)^p\,\dee y\,\dee s
$$
for $t\in I\cap(0,T')$ and a.a.~$x\in\Omega$. 
This together with \eqref{eq:1.7} implies that 
$$
\sup_{t\in I\cap (0,T')}\int_{B_\Omega(z,R)} d(y)u(y,t)\,\dee y+\int_0^{T'}\int_{B_\Omega(z,R)} d(y)u(y,s)^p\,\dee y\,\dee s<\infty
$$
for $z\in\overline{\Omega}$ and $R>0$. 
Thus \eqref{eq:2.2} holds, and assertion~(1) follows. 

We prove assertion~(2). 
Let $u$ be a solution to problem~\eqref{eq:PP} in~$Q_T$, where $T\in(0,\infty]$. 
Thanks to Definition~\ref{Definition:1.1}-(2) and Fubini's theorem, 
there exists a measurable set $\tilde{I}\subset(0,T)$ with $|(0,T)\setminus \tilde{I}|_1=0$ such that, 
for any $t\in \tilde{I}$, 
$$
\infty>u(x,t)=\int_{\overline{\Omega}} K(x,y,t)\,\dee \mu(y)
+\int_0^t\int_\Omega G(x,y,t-s)u(y,s)^p\,\dee y\,\dee s
$$
holds for a.a.~$x\in\Omega$. 
Then, by the same argument as in the proof of \eqref{eq:2.3}, 
for any $t_1$, $t_2\in I':=\tilde{I}\cup\{0\}$ with $t_1<t_2$, 
we see that \eqref{eq:2.3} holds for a.a.~$x\in\Omega$. 
This implies that $u$ is a solution to problem~\eqref{eq:E} in~$Q_T$. 
Furthermore, for any $\tau\in I'$, we see that $u_\tau$ is a solution to problem~\eqref{eq:PP} in $Q_{T-\tau}$ 
with $\mu=d(\cdot)u(\cdot,\tau)$. 
Thus assertion~(2) follows, and the proof is complete.
$\Box$\vspace{5pt}

Next, we prove that a solution to problem~\eqref{eq:PP}  
is a distributional solution to problem~\eqref{eq:PP}.
\begin{lemma}
\label{Lemma:2.3}
Let $u$ be a solution to problem~\eqref{eq:PP} in~$Q_T$, where $T\in(0,\infty]$. 
Then 
\begin{equation}
\label{eq:2.4}
\int_0^T\int_\Omega u(-\partial_t \psi-\Delta \psi )\,\dee x\,\dee t = 
\int_{\overline{\Omega}}\psi_*(x,0)\,\dee \mu(x)+\int_0^T\int_\Omega  u^p\psi\,\dee x\,\dee t
\end{equation}
for $\psi\in C^{2;1}(Q_T)\cap C^{1;0}(\overline{\Omega}\times[0,T))\cap C_0(\overline{\Omega}\times[0,T))$ 
with $|\partial_t\psi|+|\Delta\psi|\le Cd$ in $\Omega(\delta)$ for some $\delta>0$ and $C>0$ and with 
$\psi=0$ on $\partial\Omega\times(0,T)$. Here 
$$
\psi_*(x,0)=d(x)^{-1}\psi(x,0)\quad\mbox{for}\quad x\in\Omega,\quad
\psi_*(x,0)=\langle\nabla\psi(x,0),n_x\rangle\quad\mbox{for}\quad x\in\partial\Omega.
$$
\end{lemma}
{\bf Proof.}
Let $u$ be a solution to problem~\eqref{eq:PP} in~$Q_T$, where $T\in(0,\infty]$. 
Let $\psi$ be as in the above.
It follows from Lemma~\ref{Lemma:2.2} that all of the integrals in \eqref{eq:2.4} are well-defined. 
Furthermore, it follows from Definition~\ref{Definition:1.1}-(2) that 
\begin{equation}
\label{eq:2.5}
\begin{split}
 & \int_0^T\int_\Omega u(-\partial_t \psi-\Delta \psi )\,\dee x\,\dee t\\
 & =\int_0^T\int_\Omega\left(\int_{\overline{\Omega}}
 K(x,y,t)\,\dee\mu(y)+\int_0^t\int_\Omega G(x,y,t-s)u(y,s)^p\,\dee y\,\dee s\right)\\
 & \qquad\times(-\partial_t \psi-\Delta \psi)
 \,\dee x\,\dee t.
\end{split}
\end{equation}
On the other hand, by \eqref{eq:1.9} and \eqref{eq:1.10} we have 
\begin{equation}
\label{eq:2.6}
\begin{split}
 & \int_0^T\int_\Omega \int_{\overline{\Omega}} K(x,y,t)(-\partial_t \psi-\Delta \psi)\,\dee\mu(y)\,\dee x\,\dee t\\
 & =\lim_{\tau\to +0}\int_\tau^T\int_{\overline{\Omega}}\int_\Omega  K(x,y,t)(-\partial_t \psi-\Delta \psi)\,\dee x\,\dee\mu(y)\,\dee t\\
 & =\lim_{\tau\to +0}\int_{\overline{\Omega}}\int_\Omega K(x,y,\tau)\psi(x,\tau)\,\dee x\,\dee\mu(y).
\end{split}
\end{equation} 
Since 
\begin{align*}
 & \lim_{\tau\to +0}\int_\Omega G(x,y,\tau)\psi(x,\tau)\,\dee x=\psi(y,0)\quad\mbox{for}\quad y\in\Omega,\\
 & \lim_{\tau\to +0}\left\langle\nabla_y\int_\Omega G(x,y,\tau)\psi(x,\tau)\,\dee x,n_y\right\rangle=\langle\nabla_y\psi(y,0),n_y\rangle\quad\mbox{for}\quad y\in\partial\Omega,
\end{align*} 
we observe from \eqref{eq:2.6} and the definition of $K$ (see \eqref{eq:1.5})  that 
\begin{equation}
\label{eq:2.7}
\int_0^T\int_\Omega \int_{\overline{\Omega}} K(x,y,t)(-\partial_t \psi-\Delta \psi)\,\dee\mu(y)\,\dee x\,\dee t
=\int_{\overline{\Omega}}\psi_*(y,0)\,\dee \mu(y).
\end{equation}
Furthermore, setting $G(x,y,t)=0$ if $t<0$, by \eqref{eq:1.2} we obtain 
\begin{align*}
 & \int_0^T\int_\Omega\int_0^t\int_\Omega G(x,y,t-s)u(y,s)^p(-\partial_t \psi-\Delta \psi)\,\dee y\,\dee s \,\dee x\,\dee t\\
 & =\int_0^T\int_\Omega\int_0^T\int_\Omega G(x,y,t-s)u(y,s)^p(-\partial_t \psi-\Delta \psi)\,\dee y\,\dee s \,\dee x\,\dee t\\
 & =\int_0^T\int_\Omega\left(\int_0^T\int_\Omega G(x,y,t-s)(-\partial_t \psi-\Delta \psi)\,\dee x\,\dee t\right)u(y,s)^p\,\dee y\,\dee s\\
 & =\int_0^T\int_\Omega
 \lim_{\tau\to +0}\left(\int_{s+\tau}^T\int_\Omega G(x,y,t-s)(-\partial_t \psi-\Delta \psi)\,\dee x\,\dee t\right)u(y,s)^p\,\dee y\,\dee s\\
 & =\int_0^T\int_\Omega  \lim_{\tau\to +0}\left(\int_\Omega G(x,y,\tau)\psi(x,s+\tau)\,\dee x\right)u(y,s)^p\,\dee y\,\dee s\\
 & =\int_0^T\int_\Omega\psi(y,s)u(y,s)^p\,\dee y\,\dee s.
\end{align*}
This together with \eqref{eq:2.5} and \eqref{eq:2.7} implies \eqref{eq:2.4}, and Lemma~\ref{Lemma:2.3} follows. 
$\Box$\vspace{5pt}

Next, we prove Theorem~\ref{Theorem:1.1}. 
For the proof, we prepare a lemma on the short time behavior of solutions to the heat equation under the Dirichlet boundary condition. 
\begin{lemma}
\label{Lemma:2.4} 
Let $\psi\in C_0^\infty({\mathbb R}^N)$, and set 
$$
\psi_d(x,t):=\int_\Omega G(x,y,t)d(y)\psi(y)\,\dee y,\quad (x,t)\in Q_\infty.
$$
Set
$$
\psi_*(x,t):=
\left\{
\begin{array}{ll}
d(x)^{-1}\psi_d(x,t) & \mbox{for}\quad (x,t)\in\Omega\times(0,\infty),\vspace{7pt}\\
\langle\nabla_x\psi_d(x,t),n_x\rangle & \mbox{for}\quad (x,t)\in\partial\Omega\times(0,\infty).
\end{array}
\right.
$$
Then
\begin{equation}
\label{eq:2.8}
\lim_{t\to +0}\sup_{x\in B_\Omega(z,R)}|\psi_*(x,t)-\psi(x)|=0,
\quad z\in\overline{\Omega},\,\,R>0. 
\end{equation}
\end{lemma}
{\bf Proof.}
Since $d$ is Lipschitz continuous in ${\mathbb R}^N$ and $d\psi=0$ on $\partial\Omega$, 
by  parabolic regularity theorems (see e.g., \cite{LSU}*{Chapter III, Theorem 10.1})
we see that $\psi_d\in C(\overline{\Omega}\times[0,\infty))$ and 
\begin{equation}
\label{eq:2.9}
\lim_{t\to +0}\sup_{x\in B_\Omega(z,R)}|\psi_d(x,t)-d(x)\psi(x)|=0,
\quad z\in\overline{\Omega},\,\,R>0. 
\end{equation}
Let $z_*\in\partial\Omega$. 
Since $d\in C^2(B_\Omega(z_*,\delta))$ for some $\delta>0$, 
using parabolic regularity theorems again (see e.g. \cite{L}), 
we see that $\psi_d\in C^{1;0}(B_\Omega(z,\delta)\times[0,\infty))$, 
which implies that 
\begin{equation}
\label{eq:2.10}
\lim_{t\to +0}\sup_{x\in B_\Omega(z_*,\delta)}|(\nabla_x\psi_d)(x,t)-(\nabla_x(d\psi))(x)|=0. 
\end{equation}
Taking $\delta'\in(0,\delta)$ small enough, 
for any $x\in B_\Omega(z_*,\delta')$, we find $z_x\in\partial\Omega\cap B_\Omega(z_*,\delta)$ such that 
$x=z_x+d(x)n_{z_x}$. 
Then 
\begin{align*}
\psi_*(x,t) & =d(x)^{-1}\int_0^1 \frac{\dee}{\dee s}\psi_d(z_x+sd(x)n_{z_x},t)\,\dee s
=\int_0^1 \langle \nabla_x \psi_d(z_x+sd(x)n_{z_x},t),n_{z_x}\rangle\,\dee s,\\
\psi(x) & =d(x)^{-1}\int_0^1 \frac{\dee}{\dee s}(d\psi)(z_x+sd(x)n_{z_x})\,\dee s=\int_0^1 \langle\nabla_x(d\psi)(z_x+sd(x)n_{z_x}),n_{z_x}\rangle\,\dee s,
\end{align*}
for $x\in \Omega\cap B(z_*,\delta')$ and $t>0$.
These together with \eqref{eq:2.10} imply that 
\begin{equation}
\label{eq:2.11}
\lim_{t\to +0}\sup_{x\in B_\Omega(z_*,\delta')}\left|\psi_*(x,t)-\psi(x)\right|=0.
\end{equation}
Since $z_*\in\partial\Omega$ is arbitrary and $\overline{B_\Omega(z,R)}$ is compact, 
by \eqref{eq:2.9} and \eqref{eq:2.11} we obtain \eqref{eq:2.8}. 
Thus Lemma~\ref{Lemma:2.4}  follows.
$\Box$
\vspace{5pt}

Now we are ready to prove Theorem~\ref{Theorem:1.1}.
\vspace{5pt}\newline
{\bf Proof of Theorem~\ref{Theorem:1.1}.}
We prove Theorem~\ref{Theorem:1.1}-(1). 
Assume that problem~\eqref{eq:E} possesses a solution~$u$ in~$Q_T$ for some $T\in(0,\infty]$. 
Let $I$ be as in Lemma~\ref{Lemma:2.2}-(1). 
Let $\{t_j\}_{j=1}^\infty\subset I$ be such that $t_j\to +0$ as $j\to\infty$. 
By \eqref{eq:2.2} we have 
$$
\sup_{j=1,2,\dots}\,\int_{B_\Omega(z,R)} d(y)u(y,t_j)\,\dee y <\infty
$$
for $z\in\overline{\Omega}$ and $R>0$.
Applying the weak compactness of Radon measures (see e.g., \cite{EG}*{Section~1.9}) 
and taking a subsequence if necessary, 
we find $\mu\in{\mathcal M}(\overline{\Omega})$ such that 
\begin{equation}
\label{eq:2.12}
\lim_{j\to\infty}\int_\Omega \psi(y)d(y)u(y,t_j)\,\dee y=\int_{\overline{\Omega}}\psi(y)\,\dee \mu(y),
\quad \psi\in C_0({\mathbb R}^N).
\end{equation}

Let $\{s_j\}$ be a sequence in $I$ such that $\lim_{j\to\infty}s_j=0$. 
Similarly to \eqref{eq:2.12}, 
taking a subsequence if necessary, 
we find $\mu'\in{\mathcal M}(\overline{\Omega})$ such that 
\begin{equation}
\label{eq:2.13}
\lim_{j\to\infty}\int_{\Omega} \psi(y)d(y)u(y,s_j)\,\dee y=\int_{\overline{\Omega}}\psi(y)\,\dee \mu'(y),
\quad \psi\in C_0({\mathbb R}^N).
\end{equation}
Taking a subsequence again if necessary, we can assume that $t_j>s_j$ for $j=1,2,\dots$.

Let $\psi\in C_0^\infty({\mathbb R}^N)$ be such that $\psi\ge 0$ in ${\mathbb R}^N$ 
and $\mbox{supp}\,\psi\subset B(0,R)$ for some $R>0$. 
We use the same notation as in Lemma~\ref{Lemma:2.4}. 
Since $t_j$, $s_j\in I$, 
by Lemma~\ref{Lemma:2.2}-(1) we have 
$$
u(x,t_j)\ge\int_\Omega G(x,y,t_j-s_j)u(y,s_j)\,\dee y,\quad \mbox{a.a.~$x\in\Omega$},
$$
which implies that 
\begin{equation}
\label{eq:2.14}
  \begin{split}
    & \int_\Omega \psi(x)d(x)u(x,t_j)\,\dee x
   \ge \int_\Omega\int_\Omega G(x,y,t_j-s_j)u(y,s_j)d(x)\psi(x)\,\dee y\,\dee x\\
    & =\int_\Omega \psi_d(y,t_j-s_j)u(y,s_j)\,\dee y\ge\int_{B_\Omega(0,R)} \psi_d(y,t_j-s_j)u(y,s_j)\,\dee y\\
    & \ge\int_\Omega d(y)\psi(y)u(y,s_j)\,\dee y-\sup_{x\in B_\Omega(0,R)}\left|\psi_*(x,t_j-s_j)-\psi(x)\right|
    \int_{B_\Omega(0,R)}d(y)u(y,s_j)\,\dee y.
   \end{split}
\end{equation}
Combing \eqref{eq:2.2}, \eqref{eq:2.8}, \eqref{eq:2.12}, \eqref{eq:2.13}, and \eqref{eq:2.14}, 
we obtain 
$$
\int_{\overline{\Omega}}\psi(y)\,\dee \mu(y)\ge\int_{\overline{\Omega}}\psi(y)\,\dee \mu'(y).
$$
Since $\psi$ is arbitrary, we deduce that $\mu\ge\mu'$ in ${\mathcal M}(\overline{\Omega})$. 
Similarly, we have  $\mu'\ge\mu$ in ${\mathcal M}(\overline{\Omega})$. Thus we see that $\mu=\mu'$ in ${\mathcal M}(\overline{\Omega})$. 
Since $\{s_j\}\subset I$ is arbitrary, we obtain Theorem~\ref{Theorem:1.1}-(1). 

Next, we prove Theorem~\ref{Theorem:1.1}-(2). 
Let $\mu\in{\mathcal M}(\overline{\Omega})$ and let $u$ be a solution to problem~\eqref{eq:PP} in~$Q_T$, where $T\in(0,\infty]$. 
By Theorem~\ref{Theorem:1.1}-(1) and Lemma~\ref{Lemma:2.2}-(2)
we find a unique $\nu\in{\mathcal M}(\overline{\Omega})$ such that 
\begin{equation}
\label{eq:2.15}
\underset{t\to +0}{\mbox{{\rm ess lim}}}
\int_\Omega \psi(y)d(y)u(y,t)\,\dee y=\int_{\overline{\Omega}}\psi(y)\,\dee \nu(y),
\quad \psi\in C_0({\mathbb R}^N).
\end{equation}

We prove that $\mu=\nu$ in ${\mathcal M}(\overline{\Omega})$. 
Let $\psi\in C_0^\infty({\mathbb R}^N)$ be such that $\psi\ge 0$ in ${\mathbb R}^N$ and $\mbox{supp}\,\psi\subset B(0,R)$ for some $R>0$. 
We use the same notation as in Lemma~\ref{Lemma:2.4} again.  
By Lemma~\ref{Lemma:2.2} we have
\begin{equation}
\label{eq:2.16}
\begin{split}
\int_\Omega \psi(x)d(x)u(x,t)\,\dee x
 & =\int_\Omega\int_\Omega G(x,y,t-\tau)u(y,\tau)d(x)\psi(x)\,\dee x\,\dee y\\
 & +\int_\tau^t\int_\Omega\int_\Omega G(x,y,t-s)u(y,s)^pd(x)\psi(x)\,\dee x\,\dee y\,\dee s
\end{split}
\end{equation}
for a.a.~$t$, $\tau\in(0,T)$ with $t>\tau$. 
For any $T'\in(0,T)$, it follows from \eqref{eq:2.2} that
\begin{equation*}
\begin{split}
 & \int_\Omega\int_\Omega G(x,y,t-\tau)u(y,\tau)d(x)\psi(x)\,\dee x\,\dee y
\ge \int_{B_\Omega(0,R)} \psi_d(y,t-\tau)u(y,\tau)\,\dee y\\
 & \ge \int_{B_\Omega(0,R)} \psi(y)d(y)u(y,\tau)\,\dee y\\
 & \qquad\quad
 -\sup_{y\in B_\Omega(0,R)}\left|\psi_*(y,t-\tau)-\psi(y)\right|
\underset{\tau\in(0,T')}{\mbox{{\rm ess sup}}}
\int_{B_\Omega(0,R)} d(y)u(y,\tau)\,\dee y\\
 & \ge \int_\Omega \psi(y)d(y)u(y,\tau)\,\dee y-C
\sup_{x\in B_\Omega(0,R)}\left|\psi_*(x,t-\tau)-\psi(x)\right|
\end{split}
\end{equation*}
for a.a.~$t$, $\tau\in(0,T')$ with $t>\tau$. 
This together with \eqref{eq:2.15} and \eqref{eq:2.16} implies that 
\begin{align*}
\int_\Omega \psi(x)d(x)u(x,t)\,\dee x
 & \ge\int_{\overline{\Omega}} \psi(y)\,\dee \nu-C\sup_{y\in  B_\Omega(0,R)}\left|\psi_*(y,t)-\psi(y)\right|\\
 & +\int_0^t\int_\Omega \int_\Omega G(x,y,t-s)u(y,s)^pd(x)\psi(x)\,\dee x\,\dee y\,\dee s
\end{align*}
for a.a.~$t\in(0,T')$.
Then, by \eqref{eq:2.8} and \eqref{eq:2.15}  we see that 
\begin{equation}
\label{eq:2.17}
\underset{t\to +0}{\mbox{{\rm ess lim}}}
\int_0^t\int_\Omega \int_\Omega G(x,y,t-s)u(y,s)^pd(x)\psi(x)\,\dee x\,\dee y\,\dee s=0.
\end{equation}
Since 
\begin{equation*}
\begin{split}
\int_\Omega \psi(x)d(x)u(x,t)\,\dee x & =\int_\Omega \int_{\overline{\Omega}} K(x,y,t)d(x)\psi(x)\,\dee \mu(y)\,\dee x\\
 & +\int_0^t\int_\Omega\int_\Omega G(x,y,t-s)u(y,s)^p d(x)\psi(x)\,\dee x\,\dee y\,\dee s
\end{split}
\end{equation*}
for a.a.~$t\in(0,T)$ (see Definition~\ref{Definition:1.1}-(2)), by \eqref{eq:2.17} we have
\begin{equation}
\label{eq:2.18}
\int_{\overline{\Omega}} \psi(x)\,\dee \nu(x)
=\underset{t\to +0}{\mbox{{\rm ess lim}}}\,\int_\Omega \int_{\overline{\Omega}} K(x,y,t)d(x)\psi(x)\,\dee \mu(y)\,\dee x.
\end{equation}
We observe that
\begin{align*}
 & \int_\Omega \int_{\overline{\Omega}} K(x,y,t)d(x)\psi(x)\,\dee \mu(y)\,\dee x\\
 & =\int_\Omega d(y)^{-1}\left(\int_\Omega G(x,y,t)d(x)\psi(x)\,\dee x\right)\,\dee \mu(y)\\
 & \qquad\quad
 +\int_{\partial\Omega}\left(\int_\Omega \langle\nabla_y G(x,y,t),n_y\rangle d(x)\psi(x)\,\dee x\right)\,\dee \mu(y)\\
 & =\int_\Omega \psi_*(y,t)\,\dee \mu(y)+\int_{\partial\Omega}\langle\nabla\psi_d(y,t),n_y\rangle\,\dee \mu(y)
 =\int_{\overline{\Omega}} \psi_*(y,t)\,\dee \mu(y).
\end{align*}
This together with \eqref{eq:2.8} and \eqref{eq:2.18} implies that 
$$
\int_{\overline{\Omega}}\psi(x)\,\dee \nu(x)
=\lim_{t\to +0}\int_\Omega \int_{\overline{\Omega}} K(x,y,t)d(x)\psi(x)\,\dee \mu(y)\,\dee x
=\int_{\overline{\Omega}}\psi(y)\,\dee \mu(y).
$$
Since $\psi$ is arbitrary, we see that $\mu=\nu$ in ${\mathcal M}(\overline{\Omega})$. Thus Theorem~\ref{Theorem:1.1}-(2) follows. 
The proof of Theorem~\ref{Theorem:1.1} is complete.
$\Box$
\section{Interior and boundary estimates of solvable initial data}
In this section we study interior and boundary estimates of solvable initial data of problem~\eqref{eq:PP}, 
and prove Theorem~\ref{Theorem:1.2}-(1), (2). 
Furthermore, we prove Theorem~\ref{Theorem:1.3} 
to obtain Theorem~\ref{Theorem:1.2}-(4), Corollary~\ref{Corollary:1.1}, and Corollary~\ref{Corollary:1.2}.
\subsection{Preliminaries}
Let $p>1$ and 
$$
f(s):=e^{-\frac{1}{s}}\quad\mbox{if}\quad s>0,
\qquad 
f(s):=0\quad\mbox{if}\quad s\le 0.
$$
Then $f\in C^\infty({\mathbb R})$. Set 
$$
\eta(s):=\frac{f(2-s)}{f(2-s)+f(s-1)}\quad\mbox{for}\quad s\in{\mathbb R}.
$$
Then $\eta\in C^\infty([0,\infty))$, $\eta=1$ in $[0,1]$, $\eta=0$ in $[2,\infty)$, and 
$$
\eta'(s)=\frac{-f'(2-s)f(s-1)-f(2-s)f'(s-1)}{[f(2-s)+f(s-1)]^2}\le 0\quad\text{in}\quad [0, \infty).
$$
Setting 
$\eta^*(s)=0$ for $s\in[0,1)$ and $\eta^*(s)=\eta(s)$ for $s\in[1, \infty)$,  
we have
\begin{equation}
\label{eq:3.1}
|\eta'(s)|+|\eta''(s)|\le C\eta^*(s)^{\frac{1}{p}}\quad\mbox{for $s\ge 1$}. 
\end{equation}
This follows from the fact that 
$$
|f'(s)|+|f''(s)|\le f(s)^{\frac{1}{p}}\quad\mbox{for $s>0$ small enough}. 
$$ 
For any $r>0$, set 
\begin{equation}
\label{eq:3.2}
\psi_r(x,t):=\eta\left(\frac{2|x|^2+2t}{r}\right),
\quad	
\psi_r^*(x,t):=\eta^*\left(\frac{2|x|^2+2t}{r}\right), 
\end{equation}
for $(x,t)\in{\mathbb R}^N\times(0,\infty)$.
Taking into the account that $\eta^*=0$ outside $[1,2]$, 
by \eqref{eq:3.1} we see that
\begin{equation}
\label{eq:3.3}
\begin{split}
|(\partial_t\psi_r)(x,t)| & \le\frac{2}{r}\left|\eta'\left(\frac{2|x|^{2}+2t}{r}\right)\right|\le\frac{C}{r}\psi_r^*(x,t)^{\frac{1}{p}},\\
|(\nabla\psi_r)(x,t)| & \le\frac{4|x|}{r}\left|\eta'\left(\frac{2|x|^{2}+2t}{r}\right)\right|\le\frac{C|x|}{r}\psi_r^*(x,t)^{\frac{1}{p}},\\
|(\Delta\psi_r)(x,t)| & \le C\frac{|x|^2}{r^2}\left|\eta''\left(\frac{2|x|^2+2t}{r}\right)\right|
+\frac{C}{r}\left|\eta'\left(\frac{2|x|^2+2t}{r}\right)\right|
\le\frac{C}{r}\psi_r^*(x,t)^{\frac{1}{p}},
\end{split}
\end{equation}
for $(x,t)\in{\mathbb R}^N\times(0,\infty)$. 
Furthermore, since 
\begin{equation*}
\begin{array}{ll}
\displaystyle{\frac{2|x|^2+2t}{r}\ge\frac{2|x|^2}{r}\ge 2} & \quad\mbox{ for $x\in {\mathbb R}^N\setminus B(0,\sqrt{r})$ and $t\ge 0$},\vspace{5pt}\\
\displaystyle{\frac{2|x|^2+2t}{r}\ge \frac{2t}{r}\ge 2} & \quad\mbox{ for $x\in{\mathbb R}^N$ and $t\in[r,\infty)$},
\end{array}
\end{equation*}
we see that 
\begin{equation}
\label{eq:3.4}
\psi_r=0\quad\mbox{on $[{\mathbb R}^N\setminus B(0,\sqrt{r})]\times[0,\infty)$ and on ${\mathbb R}^N\times[r,\infty)$}.
\end{equation}
We use the cut-off function $\psi_r$ to obtain interior estimates of solvable initial data of problem~\eqref{eq:PP}. 

Next, we prove a lemma on a differential inequality.
\begin{lemma}
\label{Lemma:3.1}
Let $0<a<b<\infty$. 
Let $\xi\in C^1([a,b])$ be such that $\xi\ge 0$ and $\xi'\ge 0$ in $[a,b]$. 
Let $\eta\in C([a,b])$ be such that $\eta>0$ in $[a,b]$.
Assume that $\xi$ satisfies a differential inequality 
\begin{equation}
\label{eq:3.5}
m+\xi(r)\le c_*\eta(r)(\xi'(r))^{\frac{1}{\alpha}},\quad r\in (a,b),
\end{equation}
for some $c_*>0$, $m>0$, and $\alpha>1$.  
Then 
\begin{equation}
\label{eq:3.6}
m\le c_*^{\frac{\alpha}{\alpha-1}}\left(\frac{1}{\alpha-1}\right)^{\frac{1}{\alpha-1}}\left(\int_a^b \eta(r)^{-\alpha}\,\dee r\right)^{-\frac{1}{\alpha-1}}.
\end{equation}
\end{lemma}
{\bf Proof.} 
It follows from \eqref{eq:3.5} that 
$$
c_*^{-\alpha}\eta(r)^{-\alpha}\le\frac{\xi'(r)}{(m+\xi(r))^\alpha},\quad r\in(a,b).
$$
This together with $\alpha>1$ implies that
$$
c_*^{-\alpha}\int_a^b \eta(r)^{-\alpha}\,\dee r\le 
\int_a^b\frac{\xi'(r)}{(m+\xi(r))^\alpha}\,\dee r 
\le \frac{1}{\alpha-1}(m+\xi(a))^{-(\alpha-1)}\le \frac{1}{\alpha-1}m^{-(\alpha-1)}.
$$
This implies \eqref{eq:3.6}, and Lemma~\ref{Lemma:3.1} follows. 
$\Box$
\subsection{Proofs of Theorem~\ref{Theorem:1.2}-(1), (2)}
The following lemma gives interior estimates of solvable initial data of problem~\eqref{eq:PP}, 
and it is crucial for the proofs of Theorem~\ref{Theorem:1.2}-(1), (2). 
\begin{lemma}
\label{Lemma:3.2}
There exists $C=C(N,p)>0$  such that, 
if problem~\eqref{eq:PP} possesses a solution in~$Q_T$ for some $T>0$,
then 
\begin{equation}
\label{eq:3.7}
\int_{B_\Omega(z,\sigma)} d(x)^{-1}\,\dee \mu(x)
\le C\left(\int_{2\sigma^2}^{T} r^{-\frac{N}{2}(p-1)}\,\dee r\right)^{-\frac{1}{p-1}}
\end{equation}
for $z\in\Omega^c(\sqrt{T})$ and $\sigma\in(0,\sqrt{T/2})$. 
\end{lemma}
{\bf Proof.}
The proof is a modification of the arguments in \cites{IS, MP}, and it is similar to the proof of \cite{IKO}*{Theorem~1.2}.
Assume that problem~\eqref{eq:PP} possesses a solution~$u$ in~$Q_T$ for some $T>0$.
Let $z\in\Omega^c(\sqrt{T})$ and $\sigma\in(0,\sqrt{T/2})$, and fix them.
For any $r\in (0,T)$, let $\psi_r$ and $\psi_r^*$ be as in~\eqref{eq:3.2}, 
and set 
$$
\psi_{r,z}(x,t):=\psi_r(x-z,t),\quad \psi_{r,z}^*(x,t):=\psi_r^*(x-z,t),
\quad (x,t)\in Q_T.
$$
Since $d(z)>\sqrt{T}>\sqrt{r}$, it follows from \eqref{eq:3.4} that $\psi_{r,z}\in C_0^\infty(\Omega\times[0,T))$. 
Then, by Lemma~\ref{Lemma:2.3}, \eqref{eq:3.3}, and \eqref{eq:3.4} we have 
\begin{equation}
\label{eq:3.8}
\begin{split}
 & \int_{\overline{\Omega}} d(x)^{-1}\psi_{r,z}(x,0)\,\dee \mu(x)+\int_0^r\int_\Omega u(x,t)^p\psi_{r,z}(x,t)\,\dee x\,\dee t\\
 & =\int_0^r\int_\Omega u(x,t)(-\partial_t-\Delta)\psi_{r,z}(x,t)\,\dee x\,\dee t
 \le \frac{C}{r}\int_0^r\int_\Omega u(x,t)\psi_{r,z}^*(x,t)^{\frac{1}{p}}\,\dee x\,\dee t\\
 & \le \frac{C}{r}\left(\int_0^r\int_\Omega \chi_{\{\psi_{r,z}^*(x,t)>0\}}\,\dee x\,\dee t\right)^{1-\frac{1}{p}}
 \left(\int_0^r\int_\Omega u(x,t)^p\psi_{r,z}^*(x,t)\,\dee x\,\dee t\right)^{\frac{1}{p}}. 
\end{split}
\end{equation}

Let $r\in(2\sigma^2,T)$. 
Since $\sigma<\sqrt{r/2}$ and 
$\psi_{r,z}(x,0)=1$ for $x\in B_\Omega(z,\sqrt{r/2})$, 
we see that 
\begin{equation}
\label{eq:3.9}
\int_{\overline{\Omega}} d(x)^{-1}\psi_{r,z}(x,0) \,\dee x
\ge\int_{B_\Omega(z,\sqrt{r/2})}d(x)^{-1}\,\dee \mu(x)
\ge\int_{B_\Omega(z,\sigma)}d(x)^{-1}\,\dee \mu(x).
\end{equation}
On the other hand, 
since $x\in B(z,\sqrt{r})$ if $\psi_{r,z}^*(x,t)>0$, 
we have
\begin{equation}
\label{eq:3.10}
\int_0^r \int_\Omega \chi_{\{\psi_{r,z}^*(x,t)>0\}}\,\dee x\,\dee t
\le r |B(z,\sqrt{r})|_N\le Cr^{1+\frac{N}{2}}. 
\end{equation}
Combining \eqref{eq:3.8}, \eqref{eq:3.9}, and \eqref{eq:3.10}, 
we obtain 
\begin{equation}
\label{eq:3.11}
\begin{split}
 & \int_{B_\Omega(z,\sigma)}d(x)^{-1}\,\dee \mu(x)+\int_0^r\int_\Omega u(x,t)^p\psi_{r,z}(x,t)\,\dee x\,\dee t\\
 & \le Cr^{-\frac{1}{p}+\frac{N}{2}\left(1-\frac{1}{p}\right)}
 \left(\int_0^r\int_\Omega u(x,t)^p\psi_{r,z}^*(x,t)\,\dee x\,\dee t\right)^{\frac{1}{p}},\quad r\in (2\sigma^2,T).
\end{split}
\end{equation}

Let $\epsilon\in(0,2\sigma^2)$, and set
\begin{equation}
\label{eq:3.12}
W(\xi):=\int_0^\xi\int_0^{s}\int_\Omega \min\{s^{-1},\epsilon^{-1}\} u(x,t)^p\psi_{s,z}^*(x,t)\,\dee x\,\dee t\,\dee s,
\quad\xi\in(0,T).
\end{equation} 
Since $\eta$ is decreasing on $[0,\infty)$, $\eta=\eta^*$ on $[1,2]$, and $\mbox{{\rm supp}}\,\eta^*\subset [1,2]$, 
we have 
\begin{equation}
\label{eq:3.13}
\begin{split}
 & \int_0^r \psi_{s,z}^*(x,t)\min\{s^{-1},\epsilon^{-1}\}\,\dee s\\
 & \le\int_0^r\eta^*\left(\frac{2|x-z|^2+2t}{s}\right)s^{-1}\,\dee s
 =\int_{(2|x-z|^2+2t)/r}^\infty \eta^*(\tau)\tau^{-1}\,\dee\tau\\
 & \le \eta\left(\frac{2|x-z|^2+2t}{r}\right)
\int_1^2 \tau^{-1}\,\dee\tau
\le C\psi_{r,z}(x,t),\quad (x,t)\in Q_T.
\end{split}
\end{equation}
This together with \eqref{eq:3.12} implies that
\begin{equation}
\label{eq:3.14}
\begin{split}
 & \int_0^r\int_\Omega u(x,t)^p\psi_{r,z}(x,t)\,\dee x\,\dee t\\
 & \ge C\int_0^r\int_\Omega u(x,t)^p
\left(\int_0^r \psi_{s,z}^*(x,t)\min\{s^{-1},\epsilon^{-1}\}\,\dee s\right)\,\dee x\,\dee t=CW(r).
\end{split}
\end{equation}
Since
\begin{align*}
W'(r) & =\int_0^r\int_\Omega \min\{r^{-1},\epsilon^{-1}\}u(x,t)^p\psi^*_{r,z}(x,t)\,\dee x\,\dee t\\
 & =(\max\{r,\epsilon\})^{-1}\int_0^r\int_\Omega u(x,t)^p\psi^*_{r,z}(x,t)\,\dee x\,\dee t,
\end{align*}
by \eqref{eq:3.11}, \eqref{eq:3.12}, and \eqref{eq:3.14} we have
$$
\int_{B_\Omega(z,\sigma)}d(x)^{-1}\,\dee \mu(x)+W(r)
 \le Cr^{-\frac{1}{p}+\frac{N}{2}\left(1-\frac{1}{p}\right)} (\max\{r,\epsilon\}W'(r))^{\frac{1}{p}}
 \le Cr^{\frac{N}{2}\left(1-\frac{1}{p}\right)} W'(r)^{\frac{1}{p}}
$$
for $r\in(2\sigma^2,T)$. 
Then, by Lemma~\ref{Lemma:3.1} we obtain
$$
\int_{B_\Omega(z,\sigma)}d(x)^{-1}\,\dee \mu(x)
\le C\left(\int_{2\sigma^2}^T r^{-\frac{N}{2}(p-1)}\,\dee r\right)^{-\frac{1}{p-1}}. 
$$
This implies \eqref{eq:3.7}, and the proof is complete.
$\Box$
\vspace{5pt}

By Lemma~\ref{Lemma:3.2} and condition~({\bf D}) 
we obtain interior and boundary estimates of solvable initial data of problem~\eqref{eq:PP}.
\begin{lemma}
\label{Lemma:3.3}
Let $R_*$ satisfy condition~{\rm ({\bf D})}. 
Then there exists $C=C(\Omega,p)>0$  such that, 
if problem~\eqref{eq:PP} possesses a solution in~$Q_T$ for some $T\in(0,R_*^2)$,
then 
\begin{equation}
\label{eq:3.15}
\mu(B_\Omega(z,\sigma))
\le C(d(z)+\sigma)\sigma^{N-\frac{2}{p-1}}
\end{equation}
for $z\in\overline{\Omega}$ and $\sigma\in(0,\sqrt{T})$. 
\end{lemma}
{\bf Proof.}
Assume that problem~\eqref{eq:PP} possesses a solution~$u$ in~$Q_T$ for some $T\in(0,R_*^2)$. 
Let $z\in\overline{\Omega}$ and $\sigma\in(0,\sqrt{T})$. 
By Lemma~\ref{Lemma:2.2}, 
for a.a.~$\tau\in(\sigma^2/4,\sigma^2/2)$, we have:
\begin{itemize}
  \item[{\rm (i)}] 
  $\displaystyle{u(x,\tau)\ge \int_{\overline{\Omega}}K(x,y,\tau)\,\dee\mu(y)}$ for a.a.~$x\in\Omega$;
  \item[{\rm (ii)}] 
  $u_\tau$ is a solution to problem~\eqref{eq:PP} in $Q_{T-\tau}$ 
  with $\dee\mu=d(\cdot)u_\tau(\cdot,0)\,\dee x$, where $u_\tau$ is as in Lemma~\ref{Lemma:2.2}-(2)-(b). 
\end{itemize}
%
\underline{Step 1.} 
Let $d(z)>\sigma/2$. 
Then
\begin{equation*}
\begin{array}{ll} 
d(x)\ge d(z)-\sigma/4\ge\sigma/4\asymp \sqrt{\tau} & \quad\mbox{for  $x\in B(z,\sigma/4)$},\vspace{3pt}\\
 d(y)+\sqrt{\tau}\preceq d(z)+\sigma\preceq d(z) & \quad\mbox{for $y\in B(z,\sigma)$}.
\end{array}
\end{equation*}
These together property~(i) and condition~({\bf D})-(2) imply that
\begin{align*}
d(z)u(x,\tau) 
 & \ge C\tau^{-\frac{N}{2}}\frac{d(z)d(x)}{d(x)+\sqrt{\tau}}
\int_{B_\Omega(z,\sigma)}\exp\left(-\frac{|x-y|^2}{C\tau}\right)\frac{\dee \mu(y)}{d(y)+\sqrt{\tau}}\\
 & \ge C\sigma^{-N}\mu(B_\Omega(z,\sigma))
\end{align*}
for a.a.~$x\in B(z,\sigma/4)$. 
Then 
\begin{equation}
\label{eq:3.16}
\mu(B_\Omega(z,\sigma))
\le Cd(z)\int_{B(z,\sigma/4)}u(x,\tau)\,\dee x
\le C(d(z)+\sigma)\int_{B(z,\sigma/4)}u_\tau(x,0)\,\dee x.
\end{equation}
On the other hand, 
since $T-\tau>\sigma^2-\tau>\sigma^2/4$ and $d(z)>\sigma/2$, 
applying Lemma~\ref{Lemma:3.2} with $T=\sigma^2/4$ to $u_\tau$, 
we obtain 
$$
\int_{B(z,\sigma/4)}u_\tau(x,0)\,\dee x\le C\left(\int_{\sigma^2/8}^{\sigma^2/4} r^{-\frac{N}{2}(p-1)}\,\dee r\right)^{-\frac{1}{p-1}}
\le C\sigma^{N-\frac{2}{p-1}}.
$$
This together with \eqref{eq:3.16} implies \eqref{eq:3.15} in the case of $d(z)\ge \sigma/2$. 
\vspace{3pt}
\newline
\underline{Step 2.} 
Let $0\le d(z)\le\sigma/2$.
Set $z_*=z$ if $z\in\partial\Omega$. 
If $z\in\Omega$, by condition~({\bf D})-(1) and  $\sigma<\sqrt{T}<R_*$ we find $z_*\in\partial\Omega$ be such that 
\begin{equation}
\label{eq:3.17}
z=z_*+d(z)n_{z_*},\quad
\partial B(z,d(z))\cap \partial\Omega=\{z_*\}.
\end{equation}
Furthermore, setting $z_0:=z_*+(3\sigma/4)\,n_{z_*}$, 
we have
\begin{align*}
 & \,\, B(z_0,3\sigma/4)\subset\Omega,\quad \partial B(z_0,3\sigma/4)\cap\partial\Omega=\{z_*\},\\
 &
\begin{array}{ll}
d(x)\ge d(z_0)-\sigma/4=\sigma/2\succeq\sqrt{\tau} & \quad\mbox{for $x\in B(z_0,\sigma/4)$},\vspace{3pt}\\
d(y)+\sqrt{\tau}\preceq d(z)+\sigma\preceq \sigma & \quad\mbox{for $y\in B(z,\sigma)$}.
\end{array}
\end{align*}
These together property~(i) and condition~({\bf D})-(2) imply that
\begin{equation*}
\begin{split}
u(x,\tau) & 
\ge C\tau^{-\frac{N}{2}}\frac{d(x)}{d(x)+\sqrt{\tau}}\int_{B_\Omega(z,\sigma)}
\exp\left(-\frac{|x-y|^2}{C\tau}\right)\frac{\dee \mu(y)}{d(y)+\sqrt{\tau}}\\
 & \ge C\sigma^{-N-1}\mu(B_\Omega(z,\sigma))
\end{split}
\end{equation*} 
for a.a.~$x\in B(z_0,\sigma/4)$. 
Then 
\begin{equation}
\label{eq:3.18}
\mu(B_\Omega(z,\sigma))
\le C\sigma\int_{B(z_0,\sigma/4)}u(y,\tau)\,\dee y
\le C(d(z)+\sigma)\int_{B(z_0,\sigma/4)}u_\tau(y,0)\,\dee y.
\end{equation}
On the other hand, 
similarly to Step~1, 
applying Lemma~\ref{Lemma:3.2} with $T=\sigma^2/4$ to $u_\tau$, we obtain 
$$
 \int_{B(z_0,\sigma/4)}u_\tau(y,0)\,\dee y
 \le C\left(\int_{\sigma^2/8}^{\sigma^2/4} r^{-\frac{N}{2}(p-1)}\,\dee r\right)^{-\frac{1}{p-1}}
 \le C\sigma^{N-\frac{2}{p-1}}.
$$
This together with \eqref{eq:3.18} implies  \eqref{eq:3.15} in the case of $0\le d(z)\le\sigma/2$.
Thus \eqref{eq:3.18} holds for $z\in\overline{\Omega}$ and $\sigma\in(0,\sqrt{T})$, 
and the proof is complete. 
$\Box$
\vspace{5pt}

Furthermore, in the case of $p=p_N$, we have: 
\begin{lemma}
\label{Lemma:3.4}
Let $p=p_N$. Let $R_*$ satisfy condition~{\rm ({\bf D})}. 
Then there exists $C=C(\Omega)>0$  such that, 
if problem~\eqref{eq:PP} possesses a solution in $Q_T$ for some $T\in(0,R_*^2)$,
then 
\begin{equation}
\label{eq:3.19}
\mu(B_\Omega(z,\sigma))
\le C(d(z)+\sigma)\left[\log\left(e+\frac{\min\{d(z),\sqrt{T}\}}{\sigma}\right)\right]^{-\frac{N}{2}}
\end{equation}
for $z\in\overline{\Omega}$ and $\sigma\in(0,\sqrt{T})$. 
\end{lemma}
{\bf Proof.}
Let $z\in\Omega^c(\sqrt{T})$. 
If $\sigma\in(0,\sqrt{T/2})$, then Lemma~\ref{Lemma:3.2} implies \eqref{eq:3.19}. 
If not, then
$$
\frac{\min\{d(z),\sqrt{T}\}}{\sigma}\le \frac{\sqrt{T}}{\sigma}\le \sqrt{2},
$$
which together with Lemma~\ref{Lemma:3.3} and $p=p_N$ implies that 
$$
\mu(B_\Omega(z,\sigma))\le C(d(z)+\sigma)\le C(d(z)+\sigma)\left[\log\left(e+\frac{\min\{d(z),\sqrt{T}\}}{\sigma}\right)\right]^{-\frac{N}{2}},
$$
that is, \eqref{eq:3.19} holds for $\sigma\in[\sqrt{T/2},\sqrt{T})$. Thus Lemma~\ref{Lemma:3.4} follows 
in the case of $z\in\Omega^c(\sqrt{T})$. 
Similarly, Lemma~\ref{Lemma:3.4} also follows from Lemma~\ref{Lemma:3.3} in the case of $z\in\partial\Omega$. 

Let $z\in\Omega(\sqrt{T})\setminus\partial\Omega$. 
Since $d(x)\le d(z)+\sigma$ for $x\in B_\Omega(z,\sigma)$ and $0<d(z)\le\sqrt{T}$, 
by Lemma~\ref{Lemma:3.2} with $T=d(z)^2/2$ we have 
\begin{align*}
\frac{1}{d(z)+\sigma}\mu(B_\Omega(z,\sigma))
 & \le\int_{B_\Omega(z,\sigma)}d(x)^{-1}\,\dee \mu(x)
\le C\left(\int_{2\sigma^2}^{d(z)^2/2} r^{-\frac{N}{2}(p-1)}\,\dee r\right)^{-\frac{1}{p-1}}\\
 & \le C\left[\log\frac{d(z)^2}{4\sigma^2}\right]^{-\frac{N}{2}}
\le C\left[\log\left(e+\frac{d(z)}{\sigma}\right)\right]^{-\frac{N}{2}}\\
 & =C\left[\log\left(e+\frac{\min\{d(z),\sqrt{T}\}}{\sigma}\right)\right]^{-\frac{N}{2}}
\end{align*}
for $\sigma\in(0,d(z)/4)$. 
This implies \eqref{eq:3.19} in the case of $\sigma\in(0,d(z)/4)$.
Furthermore, since
$$
\frac{\min\{\sqrt{T},d(z)\}}{\sigma}\le\frac{d(z)}{\sigma}\le 4\quad\mbox{for}\quad
\sigma\in\left[\frac{d(z)}{4},\sqrt{T}\right),
$$
it follows from Lemma~\ref{Lemma:3.3} and $p=p_N$ that
$$
\mu(B_\Omega(z,\sigma))\le C(d(z)+\sigma)\le C(d(z)+\sigma)\left[\log\left(e+\frac{\min\{d(z),\sqrt{T}\}}{\sigma}\right)\right]^{-\frac{N}{2}}
$$
for $\sigma\in[d(z)/4,\sqrt{T})$. 
Thus \eqref{eq:3.19} holds in the case of $z\in\Omega(\sqrt{T})\setminus\partial\Omega$. 
The proof of Lemma~\ref{Lemma:3.4} is complete.
$\Box$\vspace{5pt}

We are ready to prove Theorem~\ref{Theorem:1.2}-(1), (2). 
\vspace{5pt}
\newline
{\bf Proofs of Theorem~\ref{Theorem:1.2}-(1), (2).}
Assume that problem~\eqref{eq:E} possesses a solution $u$ in~$Q_T$ for some $T\in(0,R_*^2)$. 
Let $\epsilon\in(0,T)$. 
By Lemma~\ref{Lemma:2.2}-(2)-(b) we see that, 
for a.a.~$\tau\in(0,\epsilon)$, $u_\tau$ is a solution to problem~\eqref{eq:PP} 
in $Q_{T-\epsilon}$ with the initial data $d u(\cdot,\tau)$. 
Applying Lemmas~\ref{Lemma:3.3} and \ref{Lemma:3.4} to the solution~$u_\tau$, we obtain 
\begin{equation}
\label{eq:3.20}
\begin{split}
 & \int_{B_\Omega(z,\sigma)}d(y)u(y,\tau)\,\dee y\\
 & \le
\left\{
\begin{array}{ll}
\displaystyle{C(d(z)+\sigma)\sigma^{N-\frac{2}{p-1}}} & \mbox{if}\quad p\not=p_N,\vspace{7pt}\\
\displaystyle{C(d(z)+\sigma)\left[\log\left(e+\frac{\min\{d(z),\sqrt{T-\epsilon}\}}{\sigma}\right)\right]^{-\frac{N}{2}}} & \mbox{if}\quad p=p_N,
\end{array}
\right.
\end{split}
\end{equation}
for $z\in\overline{\Omega}$, $\sigma\in(0,\sqrt{T-\epsilon})$, and a.a.~$\tau\in(0,\epsilon)$.
This together with Theorem~\ref{Theorem:1.1}-(1) implies that
\begin{equation}
\label{eq:3.21}
\nu(B_\Omega(z,\sigma))
\le
\left\{
\begin{array}{ll}
\displaystyle{C(d(z)+\sigma)\sigma^{N-\frac{2}{p-1}}} & \mbox{if}\quad p\not=p_N,\vspace{7pt}\\
\displaystyle{C(d(z)+\sigma)\left[\log\left(e+\frac{\min\{d(z),\sqrt{T-\epsilon}\}}{\sigma}\right)\right]^{-\frac{N}{2}}} & \mbox{if}\quad p=p_N,
\end{array}
\right.
\end{equation}
for $z\in\overline{\Omega}$, $\sigma\in(0,\sqrt{T-\epsilon})$. 
Indeed, let $\sigma'\in(0,\sigma)$ and 
$\psi\in C_0({\mathbb R}^N)$ be such that $\mbox{supp}\,\psi\subset B(z,\sigma)$, $\psi=1$ in $B(z,\sigma')$, and $0\le\psi\le 1$ in $B(z,\sigma)$.
Then it follows from Theorem~\ref{Theorem:1.1}-(1) and \eqref{eq:3.20} that 
\begin{align*}
\nu(B_\Omega(z,\sigma'))
 & \le\int_{\overline{\Omega}}\psi(y)\,\dee\nu(y)=\underset{t\to +0}{\mbox{{\rm ess lim}}}
  \int_\Omega \psi(y)d(y)u(y,t)\,\dee y\\
 & \le
\left\{
\begin{array}{ll}
\displaystyle{C(d(z)+\sigma)\sigma^{N-\frac{2}{p-1}}} & \mbox{if}\quad p\not=p_N,\vspace{7pt}\\
\displaystyle{C(d(z)+\sigma)\left[\log\left(e+\frac{\min\{d(z),\sqrt{T-\epsilon}\}}{\sigma}\right)\right]^{-\frac{N}{2}}} & \mbox{if}\quad p=p_N.
\end{array}
\right.
\end{align*}
Since $\sigma'\in(0,\sigma)$ is arbitrary, we obtain \eqref{eq:3.21}. 
Then, letting $\epsilon\to +0$, we obtain 
$$
\nu(B_\Omega(z,\sigma))
\le
\left\{
\begin{array}{ll}
\displaystyle{C(d(z)+\sigma)\sigma^{N-\frac{2}{p-1}}} & \mbox{if}\quad p\not=p_N,\vspace{7pt}\\
\displaystyle{C(d(z)+\sigma)\left[\log\left(e+\frac{\min\{d(z),\sqrt{T}\}}{\sigma}\right)\right]^{-\frac{N}{2}}} & \mbox{if}\quad p=p_N,
\end{array}
\right.
$$
for $z\in\overline{\Omega}$ and $\sigma\in(0,\sqrt{T})$. 
Thus Theorem~\ref{Theorem:1.2}-(1), (2) follow. 
$\Box$
\subsection{Proofs of Theorem~\ref{Theorem:1.2}-(4), Theorem~\ref{Theorem:1.3}, and Corollaries~\ref{Corollary:1.1} and \ref{Corollary:1.2}}
We modify the arguments in the proof of Lemma~\ref{Lemma:3.2} to prove the following lemma, 
which is crucial for the proof of Theorem~\ref{Theorem:1.2}-(4), Theorem~\ref{Theorem:1.3}, and Corollaries~\ref{Corollary:1.1} and \ref{Corollary:1.2}.
\begin{lemma}
\label{Lemma:3.5}
Assume the same conditions as in Theorem~{\rm\ref{Theorem:1.3}}.
Then there exists $C>0$ such that, 
if problem~\eqref{eq:PP} possesses a solution in~$Q_T$ for some $T\in(0,R^2]$, 
then 
\begin{equation}
\label{eq:3.22}
\int_{\Omega(\sigma)}\frac{\phi(y)}{d(y)}\,\dee \mu(y)\le C\left(\int_{2\sigma^2}^{T}\left(\int_{\Omega(\sqrt{r})}\phi(y)\,\dee y\right)^{-(p-1)}\,\dee r\right)^{-\frac{1}{p-1}}
\end{equation}
for $\sigma\in(0,\sqrt{T/2})$. Here $\phi$ and $R$ are as in Theorem~{\rm\ref{Theorem:1.3}}. 
\end{lemma}
{\bf Proof.} 
Let $\phi$ and $R$ be as in Theorem~\ref{Theorem:1.3}. 
Assume that problem~\eqref{eq:PP} possesses a solution~$u$ in~$Q_T$ for some $T\in(0,R^2]$. 
Let $\eta$ and $\eta^*$ be as in Section~3.1. 
For any $r\in(0,T)$, we set 
$$
\tilde{\psi}_r(x,t):=\eta\left(\frac{2d(x)^{2}+2t}{r}\right),
\quad	
\tilde{\psi}_r^*(x,t):=\eta^*\left(\frac{2d(x)^{2}+2t}{r}\right),
\quad (x,t)\in Q_T.
$$
Since $\partial\Omega$ is compact, $d\in C^2(\Omega(R))$, $|\nabla d|\le 1$ in $\Omega(R)$,
$r<T$, $\tilde{\psi}_r=0$ if $x\in\Omega^c(\sqrt{r})$, and $\tilde{\psi}_r=0$ if $t\ge r$, 
we see that $\tilde{\psi}_r\in C^2_0(\overline{\Omega}\times[0,T))$. 
Furthermore, by \eqref{eq:1.14} and \eqref{eq:3.3} we have 
\begin{equation}
\label{eq:3.23}
\begin{split}
|(\partial_t\tilde{\psi}_r)(x,t)| & \le\frac{C}{r}\left|\eta'\left(\frac{2d(x)^{2}+2t}{r}\right)\right|\le\frac{C}{r}\tilde{\psi}_r^*(x,t)^{\frac{1}{p}},\\
|(\nabla\tilde{\psi}_r)(x,t)| & \le\frac{Cd(x)|(\nabla d)(x)|}{r}\left|\eta'\left(\frac{2d(x)^{2}+2t}{r}\right)\right|
 \le\frac{Cd(x)}{r}\tilde{\psi}_r^*(x,t)^{\frac{1}{p}},\\
|(\Delta\tilde{\psi}_r)(x,t)| & \le \frac{Cd(x)^2|(\nabla d)(x)|^2}{r^2}\left|\eta''\left(\frac{2d(x)^2+2t}{r}\right)\right|\\
 & \quad
 +C\frac{d(x)|(\Delta d)(x)|+|(\nabla d)(x)|^2}{r}\left|\eta'\left(\frac{2d(x)^{2}+2t}{r}\right)\right|
 \le\frac{C}{r}\tilde{\psi}_r^*(x,t)^{\frac{1}{p}},
\end{split}
\end{equation}
for $(x,t)\in Q_T$. 
By \eqref{eq:1.13} we set 
$$
\lambda:=\sup_{z\in\Omega}\frac{|(\Delta \phi)(z)|}{\phi(z)}<\infty.
$$
Set $\varphi(x,t):=e^{\lambda t}\phi(x)$ for $(x,t)\in\overline{\Omega}\times[0,\infty)$. 
Then 
\begin{equation}
\label{eq:3.24}
-\partial_t\varphi-\Delta\varphi\le e^{\lambda t}(-\lambda\phi+|\Delta\phi|)\le 0\quad\mbox{in}\quad\Omega.
\end{equation}
It follows from \eqref{eq:1.13}, \eqref{eq:3.17}, \eqref{eq:3.23}, and \eqref{eq:3.24} that
\begin{equation}
\label{eq:3.25}
\begin{split}
 & (-\partial_t-\Delta)(\tilde{\psi}_r\varphi)\\
 & \le|\partial_t\tilde{\psi}_r|\varphi+|\Delta\tilde{\psi}_r|\varphi
 +2|\nabla\tilde{\psi}_r||\nabla\varphi|+\tilde{\psi}_r(-\partial_t\varphi-\Delta\varphi)\\
 & \le\frac{C}{r}(\tilde{\psi}_r^*)^{\frac{1}{p}}\varphi+\frac{C}{r}de^{\lambda t}|\nabla\phi(x)|(\tilde{\psi}_r^*)^{\frac{1}{p}}
 \le \frac{C}{r}(\tilde{\psi}_r^*)^{\frac{1}{p}}\varphi\quad\mbox{in}\quad Q_T.
\end{split}
\end{equation}
Since $\partial\Omega$ is compact, by \eqref{eq:1.13} and \eqref{eq:3.25}
we apply Lemma~\ref{Lemma:2.3} to obtain 
\begin{equation}
\label{eq:3.26}
\begin{split}
 & \int_{\overline{\Omega}} \tilde{\psi}_r(x,0)\frac{\phi(x)}{d(x)}\,\dee \mu+\int_0^r\int_\Omega u(x,t)^p\tilde{\psi}_r(x,t)\varphi(x,t)\,\dee x\,\dee t\\
 & =\int_0^r\int_\Omega u(x,t)(-\partial_t-\Delta)(\tilde{\psi}_r \varphi)(x,t)\,\dee x\,\dee t
 \le \frac{C}{r}\int_0^r \int_\Omega u(x,t)\tilde{\psi}_r^*(x,t)^{\frac{1}{p}}\varphi(x,t)\,\dee x\,\dee t\\
 & \le \frac{C}{r}\left(\int_0^r\int_\Omega \chi_{\{\tilde{\psi}_r^*(x,t)>0\}}\varphi(x,t)\,\dee x\,\dee t\right)^{1-\frac{1}{p}}
\left(\int_0^r\int_\Omega u(x,t)^p\tilde{\psi}_r^*(x,t)\varphi(x,t)\,\dee x\,\dee t\right)^{\frac{1}{p}}.
\end{split}
\end{equation}

Let $\sigma\in(0,\sqrt{T/2})$ and $r\in(2\sigma^2,T)$. Since $\sigma<\sqrt{r/2}$ and $\tilde{\psi}_r(x,0)=1$ for $x\in\Omega(\sqrt{r/2})$, 
it follows that 
\begin{equation}
\label{eq:3.27}
\int_{\overline{\Omega}} \tilde{\psi}_r(x,0)\frac{\phi(x)}{d(x)}\,\dee \mu
 \ge \int_{\Omega(\sqrt{r/2})}\frac{\phi(x)}{d(x)}\,\dee\mu(x)
 \ge \int_{\Omega(\sigma)}\frac{\phi(x)}{d(x)}\,\dee\mu(x).
\end{equation}
In addition, since $d(x)\le\sqrt{r}$ if $\tilde{\psi}_r^*(x,t)>0$ and $\lambda=0$ if $R=\infty$ (see \eqref{eq:1.14}),
we have 
$$
\int_\Omega \chi_{\{\tilde{\psi}_r^*(x,t)>0\}}\varphi(x,t)\,\dee x
\le C\int_{\Omega(\sqrt{r})}e^{\lambda t}\phi(x)\,\dee x
\le C\int_{\Omega(\sqrt{r})}\phi(x)\,\dee x
$$ 
for $t\in(0,r)\subset(0,T)$. 
This implies that 
\begin{equation}
\label{eq:3.28}
\frac{1}{r}\left(\int_0^r\int_\Omega \chi_{\{\tilde{\psi}_r^*(x,t)>0\}}\varphi(x,t)\,\dee x\,\dee t\right)^{1-\frac{1}{p}}
\le Cr^{-\frac{1}{p}}\left(\int_{\Omega(\sqrt{r})}\phi(x)\,\dee x\right)^{1-\frac{1}{p}}.
\end{equation}
Combining \eqref{eq:3.26}, \eqref{eq:3.27}, and \eqref{eq:3.28}, we obtain
\begin{equation}
\label{eq:3.29}
\begin{split}
 & \int_{\Omega(\sigma)}\frac{\phi(x)}{d(x)}\,\dee\mu(x)+\int_0^r\int_\Omega u(x,t)^p\tilde{\psi}_r(x,t)\varphi(x,t)\,\dee x\,\dee t\\
 & \le Cr^{-\frac{1}{p}}\left(\int_{\Omega(\sqrt{r})}\phi(x)\,\dee x\right)^{1-\frac{1}{p}}
 \left(\int_0^r\int_\Omega u(x,t)^p\tilde{\psi}_r^*(x,t)\varphi(x,t)\,\dee x\,\dee t\right)^{\frac{1}{p}}.
\end{split}
\end{equation}

Let $\epsilon\in(0,2\sigma^2)$, and set
$$
\tilde{W}(\xi):=\int_0^\xi\int_0^r\int_\Omega \min\{s^{-1},\epsilon^{-1}\} u(x,t)^p\tilde{\psi}_s^*(x,t)\,\dee x\,\dee t\,\dee s,
\quad\xi\in(0,T).
$$
Applying the same arguments as in \eqref{eq:3.13} and \eqref{eq:3.14},
by \eqref{eq:3.29} we obtain 
\begin{align*}
\int_{\Omega(\sigma)}\frac{\phi(x)}{d(x)}\,\dee\mu(x)+\tilde{W}(r)
 & \le Cr^{-\frac{1}{p}}\left(\int_{\Omega(\sqrt{r})}\phi(x)\,\dee x\right)^{1-\frac{1}{p}}
(\max\{r,\epsilon\}\tilde{W}'(r))^{\frac{1}{p}}\\
 & \le C\left(\int_{\Omega(\sqrt{r})}\phi(x)\,\dee x\right)^{1-\frac{1}{p}}\tilde{W}'(r)^{\frac{1}{p}}
\end{align*}
for $r\in[2\sigma^2,T)$. 
Then, by Lemma~\ref{Lemma:3.1} we obtain 
$$
\int_{\Omega(\sigma)}\frac{\phi(x)}{d(x)}\,\dee\mu(x)
\le C\left(\int_{2\sigma^2}^{T}\left(\int_{\Omega(\sqrt{r})}\phi(x)\,\dee x\right)^{-(p-1)}\,\dee r\right)^{-\frac{1}{p-1}}.
$$
This implies \eqref{eq:3.22}. Thus Lemma~\ref{Lemma:3.5} follows. 
$\Box$\vspace{5pt}

Now we are rea\dee y to complete the proof of Theorem~\ref{Theorem:1.3}, 
\vspace{5pt}
\newline
{\bf Proof of Theorem~\ref{Theorem:1.3}.}
Assume that problem~\eqref{eq:PP} possesses a solution $u$ in~$Q_T$ for some $T\in(0,R^2]$. 
Let $\epsilon\in(0,T)$. 
By Lemma~\ref{Lemma:2.2}-(2)-(b), 
for a.a.~$\tau\in(0,\epsilon)$, 
$u_\tau$ is a solution to problem~\eqref{eq:PP} in $Q_{T-\epsilon}$ with the initial data $du(\cdot,\tau)$. 
Then, by Lemma~\ref{Lemma:3.5} we have 
$$
\int_{\Omega(\sigma)}\phi(y)u(y,\tau)\,\dee y\le C\left(\int_{2\sigma^2}^{T-\epsilon}\left(\int_{\Omega(\sqrt{r})}\phi(y)\,\dee y\right)^{-(p-1)}\,\dee r\right)^{-\frac{1}{p-1}}
$$
for $\sigma\in(0,\sqrt{(T-\epsilon)/2})$ and a.a.~$\tau\in(0,\epsilon)$.
Then, similarly to \eqref{eq:3.21}, we see that 
$$
\int_{\Omega(\sigma)}\frac{\phi(y)}{d(y)}\,\dee \nu(y)\le C\left(\int_{2\sigma^2}^{T-\epsilon}\left(\int_{\Omega(\sqrt{r})}\phi(y)\,\dee y\right)^{-(p-1)}\,\dee r\right)^{-\frac{1}{p-1}}
$$
for $\sigma\in(0,\sqrt{(T-\epsilon)/2})$. 
Letting $\epsilon\to +0$, we obtain \eqref{eq:1.15}. 
Thus Theorem~\ref{Theorem:1.3} follows.
$\Box$
\vspace{5pt}

\noindent
{\bf Proof of Corollary~\ref{Corollary:1.1}}.
Let $z\in\partial\Omega$. 
By $\partial\Omega\in C^2$ 
we find $R_z>0$ such that 
there exists a bounded $C^2$-smooth domain $\tilde{\Omega}$ in ${\mathbb R}^N$ such that 
 \begin{equation}
  \label{eq:3.30}
  \Omega\cap B(z,2R_z)\subset\tilde{\Omega}\subset\Omega,
  \quad 
  d_{\partial\tilde{\Omega}}\in C^2(\tilde{\Omega}(R_z)).
  \end{equation}
  Here
  $$
  d_{\partial\tilde{\Omega}}(x):=\inf\{|x-y|\,:\,y\in\partial\tilde{\Omega}\},\quad
  \tilde{\Omega}(r):=\{x\in\overline{\tilde{\Omega}}\,:\,0\le d_{\partial\tilde{\Omega}}(x)<r\},
  $$
  where $x\in{\mathbb R}^N$ and $r>0$. 
Then 
  \begin{equation}
  \label{eq:R1}
  d(x)\asymp  d_{\partial\tilde{\Omega}}(x),\quad z\in B(z,R_z).
  \end{equation}
Assume that problem~\eqref{eq:E} with $p\ge 2$ possesses a solution~$u$ in~$Q_T$ for some $T\in(0,R_z^2)$. 
Then, for a.a.~$\tau\in(0,T)$, 
$u_\tau$ defined by $u_\tau(x,t):=u(x,t+\tau)$ is a solution to problem~\eqref{eq:PP} in $Q_{T-\tau}$ 
with $\mu=d(\cdot)u(\cdot,\tau)$.
Fix $\tau\in(0,T)$ and let $z\in\partial\Omega$. 
Since $G_{\tilde{\Omega}}\le G_\Omega$ in $\tilde{\Omega}\times\tilde{\Omega}\times(0,\infty)$, 
$u_\tau$ is a supersolution to problem~\eqref{eq:PP} in~$\tilde{\Omega}\times(0,T-\tau)$ with $\mu=d(\cdot)u(\cdot,\tau)$. 
Then Lemma~\ref{Lemma:2.1} implies that 
there exists a solution to problem~\eqref{eq:PP} in~$\tilde{\Omega}\times(0,T-\tau)$ with $\mu=d(\cdot)u(\cdot,\tau)$. 

Let $\phi$ be the first Dirichlet eigenfunction in $\tilde{\Omega}$ for $-\Delta$. 
Then it follows from Hopf's lemma that 
\begin{equation}
\label{eq:3.31}
\phi(x)\ge Cd_{\partial\tilde{\Omega}}(x),\quad x\in\tilde{\Omega}.
\end{equation} 
This together with $\phi\in C^1(\overline{\tilde{\Omega}})$ and the boundedness of $\tilde{\Omega}$ implies that
$$
d_{\partial\tilde{\Omega}}(x)|\nabla\phi(x)|\le Cd_{\partial\tilde{\Omega}}(x)\le C\phi(x)\mbox{ for $x\in\tilde{\Omega}$}.
$$
Furthermore, since $\phi=0$ on $\partial\tilde{\Omega}$, 
it follows from $\|\nabla\phi\|_{L^\infty(\tilde{\Omega})}<\infty$ that 
$$
\int_{\tilde{\Omega}(\sigma)}\phi(y)\,\dee y\le C\sigma^2,\quad \sigma\in(0,\sqrt{T}).
$$
Then, applying Theorem~\ref{Theorem:1.3} with Theorem~\ref{Theorem:1.1}-(2), by \eqref{eq:R1} and\eqref{eq:3.31} we see that 
\begin{equation*}
\begin{split}
 & \int_{B(z,\sqrt{T-\tau})\cap \Omega(\sigma)}d(y)u(y,\tau)\,\dee y
 \le C\int_{B(z,\sqrt{T-\tau})\cap \Omega(\sigma)}d_{\partial\tilde{\Omega}}(y)u(y,\tau)\,\dee y\\
  & \qquad\quad
  \le C\int_{\tilde{\Omega}(\sigma)} \phi(y)u(y,\tau)\,\dee y\\
 & \qquad\quad
 \le C\left(\int_{2\sigma^2}^{T-\tau}
\left(\int_{\tilde{\Omega}(\sqrt{r})}\phi(y)\,\dee y\right)^{-(p-1)}\,\dee r\right)^{-\frac{1}{p-1}}
\le C\left(\int_{2\sigma^2}^{T-\tau} r^{-(p-1)}\,\dee r\right)^{-\frac{1}{p-1}}\\
 & \qquad\quad
 \le 
 \left\{
 \begin{array}{ll}
 C\sigma^{2\frac{p-2}{p-1}} & \mbox{if}\quad p>2,\vspace{5pt}\\
 \displaystyle{C\left[\log\left(e+\frac{\sqrt{T-\tau}}{\sigma}\right)\right]^{-1}} & \mbox{if}\quad p=2,
 \end{array}
 \right.
\end{split}
\end{equation*}
for $\sigma\in(0,\sqrt{(T-\tau)/4})$ and a.a.~$\tau\in(0,T)$. 
Therefore, 
similarly to \eqref{eq:3.21}, we see that 
for any $z\in\partial\Omega$, 
we find $R_z>0$ and $C_z>0$ such that, 
if problem~\eqref{eq:PP} with $p\ge 2$ possesses a solution in~$Q_T$ for some $T\in(0,R_z^2)$, 
then 
\begin{equation}
\label{eq:3.32}
\nu(B(z,\sqrt{T})\cap \Omega(\sigma))
 \le 
 \left\{
 \begin{array}{ll}
 C_z\sigma^{2\frac{p-2}{p-1}} & \mbox{if}\quad p>2,\vspace{5pt}\\
 \displaystyle{C_z\left[\log\left(e+\frac{\sqrt{T}}{\sigma}\right)\right]^{-1}} & \mbox{if}\quad p=2,
 \end{array}
 \right.
\end{equation}
for $\sigma\in(0,\sqrt{T/4})$. 
Since $\Omega$ is uniformly regular of class $C^{2,\theta}$ for some $\theta\in(0,1)$, 
we see that $\inf_{z\in\partial\Omega}R_z>0$ and $\sup_{z\in\partial\Omega}C_z<\infty$. 
Then Theorem~\ref{Theorem:1.2}-(4) follows.
$\Box$
\vspace{5pt}

Similarly, we prove Theorem~\ref{Theorem:1.2}-(4).
\vspace{3pt}
\newline
{\bf Proof of Theorem~\ref{Theorem:1.2}-(4).}
Let $p\ge 2$. 
Assume that problem~\eqref{eq:E} possess a solution $u$ in $Q_T$ for some $T>0$. 
By \eqref{eq:3.32}, for any $z\in\partial\Omega$, 
the initial trace $\nu$ of $du$ satisfies 
$\nu(B(z,\sigma)\cap\partial\Omega)=0$ for some $\sigma>0$. 
This implies that $\nu(\partial\Omega)=0$, and Theorem~\ref{Theorem:1.2}-(4) follows. 
$\Box$
\vspace{5pt}

At the end of this subsection we prove Corollary~\ref{Corollary:1.2}.
\vspace{5pt}
\newline
{\bf Proof of Corollary~\ref{Corollary:1.2}.}
Let $\Omega$ be the exterior of $C^2$-smooth compact set in ${\mathbb R}^N$, where $N\ge 2$. 
The proof is by contradiction. 
Assume that problem~\eqref{eq:PP} possesses a non-trivial global-in-time solution~$u$. 
Then, by Definition~\ref{Definition:1.1}-(2) and Lemma~\ref{Lemma:2.2}-(2)-(b)
we find $\tau>0$ such that 
\begin{equation}
\label{eq:3.33}
u(x,\tau)>0,\quad\mbox{a.a.~$x\in\Omega$}
\end{equation}
and 
$u_\tau$ is a solution to problem~\eqref{eq:PP} in $Q_\infty$ with the initial data $d(\cdot)u(\cdot,\tau)$. 
Let $R_0>0$ be such that $\Omega_0:={\mathbb R}^N\setminus\overline{B(0,R_0)}\subset\Omega$. 
Since $G_{\Omega_0}\le G_\Omega$ in $\Omega_0\times\Omega_0\times(0,\infty)$,  
$u_\tau$ is a supersolution to problem~\eqref{eq:PP} in $\Omega_0\times(0,\infty)$. 
Then Lemma~\ref{Lemma:2.1} implies that 
there exists a global-in-time solution to problem~\eqref{eq:PP} 
with $\Omega$ and $\mu$ replaced by $\Omega_0$ and  $d(\cdot)u(\cdot,\tau)$, respectively.

Let $d_{\Omega_0}(x):=|x|-R_0$ for $x\in{\mathbb R}^N$. Let
$$
\phi(x):=1-\left(\frac{|x|}{R_0}\right)^{-(N-2)}\quad\mbox{if}\quad N\ge 3,
\qquad
\phi(x)=\log\frac{|x|}{R_0}\quad\mbox{if}\quad N=2,
$$
for $x\in\Omega_0$. 
Then $\Delta\phi=0$ in $\Omega_0$, $\phi>0$ in $\Omega_0$, $\phi=0$ on $\partial\Omega_0$, 
$\phi\asymp d_{\Omega_0}$ near $\partial\Omega$, and $\sup_{x\in\Omega_0}(1+d_{\Omega_0})|\nabla\phi(x)|<\infty$. 
These imply that \eqref{eq:1.4} and \eqref{eq:1.5} hold with $R=\infty$ and $d=d_{\Omega_0}$.
Furthermore, setting
$$
\Omega_0(r):=\{x\in{\mathbb R}^N\,:\,R_0\le|x|\le R_0+r\},\quad r>0,
$$
we have 
\begin{equation}
\label{eq:3.34}
\int_{\Omega_0(\sqrt{r})}\phi(x)\,\dee x
\le
\left\{
\begin{array}{ll}
Cr^{\frac{N}{2}}\quad & \mbox{if}\quad N\ge 3,\vspace{5pt}\\
Cr\log r\quad & \mbox{if}\quad N=2,
\end{array}
\right.
\end{equation}
for $r$ large enough. 
We apply Theorem~\ref{Theorem:1.1}-(2) and Theorem~\ref{Theorem:1.3} with $R=T=\infty$ to obtain 
$$
\int_{\Omega_0(\sigma)} \phi(y)u(y,\tau)\,\dee y
\le C\left(\int_{\sigma^2}^\infty 
\left(\int_{\Omega_0(\sqrt{r})}\phi(x)\,\dee x\right)^{-(p-1)}\,\dee r\right)^{-\frac{1}{p-1}}=0
$$
for $\sigma\in(0,\infty)$. 
This contradicts \eqref{eq:3.33}. 
Thus problem~\eqref{eq:PP} possesses no non-trivial global-in-time solutions, 
and Corollary~\ref{Corollary:1.2} follows.
$\Box$
\section{Improved boundary estimates of solvable initial data}
We improve boundary estimates of solvable initial data of problem~\eqref{eq:PP} in the case of $p=p_{N+1}$, 
and complete the proof of Theorem~\ref{Theorem:1.2}. 
\begin{lemma}
\label{Lemma:4.1}
Let $p=p_{N+1}$. Let $R_*$ satisfy condition~{\rm ({\bf D})}. 
Then there exists $C=C(\Omega)>0$ 
such that, if problem~\eqref{eq:PP} possesses a solution  in~$Q_T$ for some $T\in(0,R_*^2)$, 
then
\begin{equation}
\label{eq:4.1}
\mu(B_\Omega(z,\sigma))
\le C\biggr[\log\biggr(e+\frac{\sqrt{T}}{\sigma}\biggr)\biggr]^{-\frac{N+1}{2}}
\end{equation} 
for $z\in\partial\Omega$ and $\sigma\in(0,\sqrt{T})$.
\end{lemma}
{\bf Proof.}
Let $p=p_{N+1}$. 
Assume that problem~\eqref{eq:PP} possesses a solution~$u$ in~$Q_T$ for some $T\in(0,R_*^2)$. 
By Lemma~\ref{Lemma:3.3} we have 
$$
\mu(B_\Omega(z,\sigma))\le C\sigma^{N+1-\frac{2}{p-1}}=C\le C\biggr[\log\biggr(e+\frac{\sqrt{T}}{\sigma}\biggr)\biggr]^{-\frac{N+1}{2}}
$$
for $z\in\partial\Omega$ and $\sigma\in[\sqrt{T/2},\sqrt{T})$. 
Thus it suffices to prove \eqref{eq:4.1} in the case of $\sigma\in(0,\sqrt{T/2})$.

Let $z\in\partial\Omega$ and $\sigma\in(0,\sqrt{T/2})$, and fix them. 
For any $r\in(0,T/4)$, set
$$
z_r:=z+2\sqrt{r}n_z. 
$$
Then it follows from ({\bf D})-(2) and $2\sqrt{r}<\sqrt{T}< R_*$ that 
\begin{equation}
\label{eq:4.2}
B(z_r,\sqrt{r})\subset B_r:=B(z_r,2\sqrt{r})\subset\Omega.
\end{equation}
Let $\psi_r$ and $\psi_r^*$ be as in \eqref{eq:3.2}. 
Set 
$$
\hat{\psi}_{r,z}(x,t):=\psi_r(x-z_r,t),\quad \hat{\psi}_{r,z}^*(x,t):=\psi_r^*(x-z_r,t),
\quad (x,t)\in Q_T.
$$
Then, by \eqref{eq:3.4} and \eqref{eq:4.2} 
we see that $\hat{\psi}_{r,z}=0$ on $[\Omega\setminus B(z_r,\sqrt{r})]\times[0,\infty)$ and on $\Omega\times[r,T)$. 
Let $\phi_*$ be the first Dirichlet eigenfunction in $B(0,1)$ for $-\Delta$ and $\lambda_*$ the corresponding first Dirichlet eigenvalue. 
Then $\|\nabla\phi_*\|_{L^\infty(B(0,1))}<\infty$ and 
$C^{-1}\mbox{dist}\,(x,\partial B(0,1))\le \phi_*(x)\le C\mbox{dist}\,(x,\partial B(0,1))$ for $x\in B(0,1)$. 
Set 
$$
\phi_r(x):=2\sqrt{r}\phi_*\left((2\sqrt{r})^{-1}\left(x-z_r\right)\right)\mbox{ for $x\in B_r$},
\quad
\lambda:=\lambda_*/(4r). 
$$
Then $\phi_r$ satisfies 
\begin{align}
\label{eq:4.3}
 & -\Delta \phi_r=\lambda\phi_r\quad\mbox{in}\quad B_r,\quad
\phi_r>0\quad\mbox{in}\quad B_r,\quad
\phi_r=0\quad\mbox{on}\quad\partial B_r,\\
\label{eq:4.4}
 & \|\nabla\phi_r\|_{L^\infty(B_r)}\le C,
 \quad
 C^{-1}\mbox{dist}\,(x,\partial B_r)\le\phi_r(x)\le C\mbox{dist}\,(x,\partial B_r)\mbox{ for $x\in B_r$}. 
\end{align}
{Here the constants $C$ in \eqref{eq:4.4} are independent of $r\in(0,T/4)$. 
Since
$|x-z_r|\le\sqrt{r}\le \mbox{dist}\,(x,\partial B_r)$ for $x\in B(z_r,\sqrt{r})$ (see \eqref{eq:4.2}), 
by \eqref{eq:3.3}, \eqref{eq:4.3}, and \eqref{eq:4.4} we have 
\begin{equation}
\label{eq:4.5}
\begin{split}
 & -\partial_t(e^{\lambda t}\phi_r\hat{\psi}_{r,z})-\Delta(e^{\lambda t}\phi_r\hat{\psi}_{r,z})\\
 & =-e^{\lambda t}\phi_r\partial_t\hat{\psi}_{r,z}-\lambda e^{\lambda t}\phi_r\hat{\psi}_{r,z}-e^{\lambda t}\hat{\psi}_{r,z}\Delta\phi_r
 -2e^{\lambda t}\langle\nabla\phi_r,\nabla\hat{\psi}_{r,z}\rangle- e^{\lambda t}\phi_r\Delta \hat{\psi}_{r,z}\\
 &  =-e^{\lambda_*\frac{t}{4r}}\phi_r\partial_t\hat{\psi}_{r,z}
 -2e^{\lambda_* \frac{t}{4r}}\langle\nabla\phi_r,\nabla\hat{\psi}_{r,z}\rangle-e^{\lambda_*\frac{t}{4r}}\phi_r\Delta \hat{\psi}_{r,z}\\
 & \le \frac{C}{r}\phi_r(\hat{\psi}_{r,z}^*)^{\frac{1}{p}}+C\frac{|x-z_r|}{r}(\hat{\psi}_{r,z}^*)^{\frac{1}{p}}
\le \frac{C}{r}\phi_r(\hat{\psi}_{r,z}^*)^{\frac{1}{p}}
\end{split}
\end{equation}
for $x\in B(z_r,\sqrt{r})$ and $t\in(0,r)$. 
Then, 
setting $u_r(x,t):=u(x,t+r)$ for $(x,t)\in Q_T$, 
by Lemma~\ref{Lemma:2.2}-(2)-(b), Lemma~\ref{Lemma:2.3}, and \eqref{eq:4.5}
we have 
\begin{equation}
\label{eq:4.6}
\begin{split}
 & \int_{\Omega} u_r(x,0)\phi_r(x)\hat{\psi}_{r,z}(x,0) \,\dee x+\int_0^r\int_\Omega u_r(x,t)^pe^{\lambda t}\phi_r(x)\hat{\psi}_{r,z}(x,t)\,\dee x\,\dee t\\
 & =\int_0^r\int_\Omega u_r(x,t)(-\partial_t-\Delta)(e^{\lambda t}\phi_r\hat{\psi}_{r,z})(x,t)\,\dee x\,\dee t\\
 & \le \frac{C}{r}\int_0^r\int_\Omega u_r(x,t)\phi_r(x)\hat{\psi}_{r,z}^*(x,t)^{\frac{1}{p}}\,\dee x\,\dee t\\
 & \le \frac{C}{r}\left(\int_0^r\int_\Omega \chi_{\{\hat{\psi}_{r,z}^*(x,t)>0\}}\phi_r(x)\,\dee x\,\dee t\right)^{1-\frac{1}{p}}
 \left(\int_0^r\int_\Omega u_r(x,t)^p\hat{\psi}_{r,z}^*(x,t)\phi_r(x)\,\dee x\,\dee t\right)^{\frac{1}{p}}
\end{split}
\end{equation}
for a.a.~$r\in(0,T/4)$. 

Let $r\in (\sigma^2/4,T/4)$. Since 
$\hat{\psi}_{r,z}(x,0)=1$ for $x\in B(z_r,\sqrt{r/2})\subset\Omega$, 
we observe from \eqref{eq:4.2} and \eqref{eq:4.4} that 
\begin{equation}
\label{eq:4.7}
\int_{\Omega} u_r(x,0)\phi_r(x)\hat{\psi}_{r,z}(x,0) \,\dee x
 \ge C\sqrt{r}\int_{B(z_r,\sqrt{r/2})}u(x,r)\,\dee x.
\end{equation}
Furthermore, it follows from condition~({\bf D})-(2) that
\begin{align*}
u(x,r) & \ge\int_{\overline{\Omega}}K(x,y,r)\,\dee\mu(y)\\
 & \ge Cr^{-\frac{N}{2}}\frac{d(x)}{d(x)+\sqrt{r}}\int_{B_\Omega(z,\sigma)}\exp\left(-C\frac{|x-y|^2}{r}\right)\frac{\dee\mu(y)}{d(y)+\sqrt{r}}
\ge Cr^{-\frac{N+1}{2}}\mu(B(z,\sigma))
\end{align*}
for $x\in B(z_r,\sqrt{r/2})$. 
This together with \eqref{eq:4.7} implies that
\begin{equation}
\label{eq:4.8}
\int_{\Omega} u_r(x,0)\phi_r(x)\hat{\psi}_{r,z}(x,0) \,\dee x
\ge C\mu(B(z,\sigma)).
\end{equation}
On the other hand, 
since $x\in B(z_r,\sqrt{r})$ if $\hat{\psi}_{r,z}^*(x,t)>0$, 
by \eqref{eq:4.4} we have
\begin{equation}
\label{eq:4.9}
\int_0^r \int_\Omega \chi_{\{\hat{\psi}_{r,z}^*(x,t)>0\}}\phi_r(x)\,\dee x\,\dee t
\le Cr\sqrt{r}|B(z_r,\sqrt{r})|_N\le Cr^{1+\frac{N+1}{2}}. 
\end{equation}
Combining \eqref{eq:4.6}, \eqref{eq:4.8}, and \eqref{eq:4.9}, 
we obtain 
\begin{equation}
\label{eq:4.10}
\begin{split}
 & \mu(B(z,\sigma))+\int_0^r\int_\Omega u(x,t+r)^p\phi_r(x)\hat{\psi}_{r,z}(x,t)\,\dee x\,\dee t\\
 & \le Cr^{-\frac{1}{p}+\frac{N+1}{2}\left(1-\frac{1}{p}\right)}
 \left(\int_0^r\int_\Omega u(x,t+r)^p\phi_r(x)\hat{\psi}_{r,z}^*(x,t)\,\dee x\,\dee t\right)^{\frac{1}{p}}
\end{split}
\end{equation}
for a.a.~$r\in(\sigma^2/4,T/4)$. Then \eqref{eq:4.10} holds for all $r\in[\sigma^2/4,T/4)$.

Let $\epsilon\in(0,\sigma^2/4)$.  
Set 
$$
W(\xi):=\int_0^\xi\int_0^s\int_\Omega \min\{s^{-1},\epsilon^{-1}\}u(x,t+s)^p\phi_{s}(x)\hat{\psi}_{s,z}^*(x,t) \,\dee x\,\dee t\,\dee s,
\quad \xi\in\left(0,\frac{T}{4}\right). 
$$
Then, applying the same arguments as in \eqref{eq:3.13} and \eqref{eq:3.14},
by \eqref{eq:4.10} we obtain
$$
\mu(B(z,\sigma))+W(r)
\le Cr^{-\frac{1}{p}+\frac{N+1}{2}\left(1-\frac{1}{p}\right)} (\max\{r,\epsilon\}W'(r))^{\frac{1}{p}}\\
\le Cr^{\frac{N+1}{2}\left(1-\frac{1}{p}\right)} W'(r)^{\frac{1}{p}}
$$
for $r\in[\sigma^2/4,T/4)$. 
Since $p=p_{N+1}$, by Lemma~\ref{Lemma:3.1} we obtain 
\begin{align*}
\mu(B(z,\sigma)) & \le C\left(\int_{\sigma^2/4}^{T/4} r^{-\frac{N+1}{2}(p-1)}\,\dee r\right)^{-\frac{1}{p-1}}\\
 & \le C\biggr[\log\frac{T}{\sigma^2}\biggr]^{-\frac{N+1}{2}}
\le C\biggr[\log\biggr(e+\frac{\sqrt{T}}{\sigma}\biggr)\biggr]^{-\frac{N+1}{2}}.
\end{align*}
This implies \eqref{eq:4.1} for $\sigma\in(0,\sqrt{T/2})$, and Lemma~\ref{Lemma:4.1} follows.
$\Box$\vspace{5pt}

Now we are ready to complete the proof of Theorem~\ref{Theorem:1.2}.
\vspace{5pt}
\newline
{\bf Proof of Theorem~\ref{Theorem:1.2}.}
By the results in Section~3 
it suffices to prove Theorem~\ref{Theorem:1.2}-(3). 
Assume that problem~\eqref{eq:PP} with $p=p_{N+1}$ possesses a solution in~$Q_T$ for some $T\in(0,R_*^2)$. 
Let $\epsilon\in(0,T)$.  
By Lemma~\ref{Lemma:2.2}-(2)-(b) and Lemma~\ref{Lemma:4.1}
we have 
$$
\int_{B_\Omega(z,\sigma)}u(y,\tau)\,\dee y\le C\left[\log\left(e+\frac{\sqrt{T-\epsilon}}{\sigma}\right)\right]^{-\frac{N+1}{2}}
$$
for $z\in\partial\Omega$, $\sigma\in(0,\sqrt{T-\epsilon})$, and a.a.~$\tau\in(0,\epsilon)$.
Then, similarly to \eqref{eq:3.21}, 
by Theorem~\ref{Theorem:1.1}-(1) we see that 
$$
\nu(B_\Omega(z,\sigma))\le C\left[\log\left(e+\frac{\sqrt{T-\epsilon}}{\sigma}\right)\right]^{-\frac{N+1}{2}}
$$
for $z\in\partial\Omega$ and $\sigma\in(0,\sqrt{T-\epsilon})$. 
Letting $\epsilon\to +0$, we obtain the desired estimate of $\nu(B_\Omega(z,\sigma))$ 
for $z\in\partial\Omega$ and $\sigma\in(0,\sqrt{T})$, 
and Theorem~\ref{Theorem:1.2}-(3) holds. 
Thus Theorem~\ref{Theorem:1.2} follows.
$\Box$
\section{Solvability of solutions}
In this section we develop the arguments in \cite{HIT02} to obtain sufficient conditions 
for the existence of solutions to problem~\eqref{eq:PP}. 
The arguments in this section show the optimality of 
our interior and boundary estimates of initial traces of solutions to problem~\eqref{eq:E}. 
By Remark~\ref{Remark:1.1}-(2)  
we find $R_{**}\in(0,R_*]$ with the following properties:
\begin{itemize}
  \item[({\bf D}')] 
  \begin{itemize}
  \item[(1)] ({\bf D})-(1), (2) hold;
  \item[(2)] 
there exists $C_\Omega>0$ such that 
  $$
  \sup_{z\in\overline{\Omega}}\int_{\partial\Omega\cap B(z,\sigma)}\,\dee S(y)\le C_\Omega\sigma^{N-1}
  $$
  for $\sigma\in(0,R_{**})$;
  \item[(3)] there exists $C_G'>0$ such that 
  $$
  G(x,y,t)\le C_G't^{-\frac{N}{2}}\frac{d(x)}{d(x)+\sqrt{t}}\frac{d(y)}{d(y)+\sqrt{t}}\exp\left(-\frac{|x-y|^2}{C_G't}\right)
  $$
  for $(x,y,t)\in \Omega\times \Omega\times(0,R_{**}^2)$.
  \end{itemize} 
\end{itemize}
Similarly to Remark~\ref{Remark:1.1}-(1), 
if $\Omega$ is a uniformly regular domain of class $C^2$, 
then condition~({\bf D}') holds for some $R_{**}>0$. 
Notice that $R_{**}=\infty$ if $\Omega={\mathbb R}^N_+$. 
(See Remark~\ref{Remark:1.1}-(2) and \cite{HIT02}*{Lemma~2.2}.)
\subsection{Preliminary lemmas}
For any $\mu\in{\mathcal M}(\overline{\Omega})$, 
set
$$
[\K(t)\mu](x) := \int_{\overline{\Omega}} K(x,y,t)\,\dee\mu(y),
\quad
[\G(t)\mu](x) := \int_{\Omega} G(x,y,t)\,\dee\mu(y),
\quad (x,t)\in Q_\infty.
$$ 
Then the following two lemmas hold.
\begin{lemma}
\label{Lemma:5.1}
Let $R_{**}$ satisfy condition~{\rm ({\bf D}')}. 
Then there exists $C>0$ such that 
$$
\|\K(t)\mu\|_{L^\infty(\Omega)}\le Ct^{-\frac{N}{2}}\sup_{z\in\overline{\Omega}}
\int_{B_\Omega(z,\sqrt{t})}\frac{\dee \mu(y)}{d(y)+\sqrt{t}},\quad t\in(0,R_{**}^2),
$$
for $\mu\in{\mathcal M}(\overline{\Omega})$.
\end{lemma}
{\bf Proof.}
Let $\mu\in{\mathcal M}(\overline{\Omega})$. 
It follows from \cite{HI01}*{Lemma~2.1} and condition~({\bf D}')-(3) that
\begin{equation}
\label{eq:5.1}
\|\K(t)\mu\|_{L^\infty(\Omega)}\le Ct^{-\frac{N}{2}}\sup_{z\in{\mathbb R}^N}\int_{B_\Omega(z,C\sqrt{t})}
\frac{\dee \mu(y)}{d(y)+\sqrt{t}}
\end{equation}
for $t\in(0,R_{**}^2)$. 
On the other hand, we have:
\begin{itemize}
  \item 
  for any $z\in{\mathbb R}^N$, $L\ge 1$, and $r>0$, 
  there exists a set $\{z_k\}_{k=1}^m\subset{\mathbb R}^N$ such that
  \begin{equation}
  \label{eq:5.2}
  B(z,Lr)\subset\bigcup_{k=1}^m B(z_k,r),
  \end{equation}
  where $m$ depends only on $N$ and $L$ (see e.g., \cite{ISato}*{Lemma 2.1});  
  \item 
  if $B_\Omega(z,r/2)\not=\emptyset$, we find $z'\in\overline{\Omega}$ such that $B_\Omega(z,r/2)\subset B_\Omega(z',r)$.
\end{itemize}
Then, by \eqref{eq:5.1} we see that 
\begin{equation}
\label{eq:5.3}
\|\K(t)\mu\|_{L^\infty(\Omega)}\le Ct^{-\frac{N}{2}}\sup_{z\in{\mathbb R}^N}\int_{B_\Omega(z,\sqrt{t}/2)}\frac{\dee \mu(y)}{d(y)+\sqrt{t}}
\le Ct^{-\frac{N}{2}}\sup_{z\in\overline{\Omega}}\int_{B(z,\sqrt{t})}
\frac{\dee \mu(y)}{d(y)+\sqrt{t}}
\end{equation}
for $t\in(0,R_{**}^2)$. Thus Lemma~\ref{Lemma:5.1} follows.
$\Box$
\begin{lemma}
\label{Lemma:5.2}
Let $\Psi$ be a strictly increasing, nonnegative, and convex function in $[0,\infty)$ such that $\Psi(0)=0$. 
\begin{itemize}
  \item[{\rm (1)}]
  Let $f$ be a nonnegative measurable function in $\Omega$. 
  Then 
  $$
  [\G(t)f](x)\le\Psi^{-1}([\G(t)\Psi(f)](x)),\quad (x,t)\in Q_\infty.
  $$
  \item[{\rm (2)}] 
  Let $R_{**}$ satisfy condition~{\rm ({\bf D}')}. 
  Set 
  $$
  k(x,t):=\int_{\partial\Omega}K(x,y,t)\,\dee S(y),\quad (x,t)\in Q_\infty.
  $$
  Then, for any $T\in(0,R_{**}^2]$, there exists $C_*>0$ such that 
  $$
  \|k(t)\|_{L^\infty(\Omega)}\le C_*t^{-1},\quad t\in(0,T). 
  $$
  Furthermore, for any nonnegative measurable function $h$ in $\partial\Omega$, 
  $$
  [\tilde{\K}(t)h](x)\le\Psi^{-1}([\tilde{\K}(t)\Psi(h)](x)),\quad (x,t)\in Q_T,
  $$
  where
  \begin{equation*}
  [\tilde{\K}(t)h](x):=C_*^{-1}t\int_{\partial\Omega}K(x,y,t)h(y)\,\dee S(y).
  \end{equation*}
\end{itemize}
\end{lemma}
{\bf Proof.}
We prove assertion~(1). 
It follows from Jensen's inequality that 
\begin{equation}
\label{eq:5.4}
[\G(t)f](x)\le [\G(t)1](x)\Psi^{-1}\left(\frac{1}{[\G(t)1](x)}\int_\Omega G(x,y,t)\Psi(f(y))\,\dee y\right),
\quad (x,t)\in Q_\infty.
\end{equation}
On the other hand, since $\Psi^{-1}$ is concave and $\Psi^{-1}(0)=0$, we have
\begin{equation}
\label{eq:5.5}
\lambda\Psi^{-1}(r)=\lambda\Psi^{-1}(\lambda^{-1}\lambda r)
\ge \lambda\left[(1-\lambda^{-1})\Psi^{-1}(0)+\lambda^{-1}\Psi^{-1}(\lambda r)\right]\ge\Psi^{-1}(\lambda r)
\end{equation}
for $r\in[0,\infty)$ and $\lambda\ge 1$.
This together with the fact that $[\G(t) 1](x)\le 1$ for $(x,t)\in Q_\infty$ implies that 
\begin{equation}
\label{eq:5.6}
\Psi^{-1}\left(\frac{1}{[\G(t)1](x)}\int_\Omega G(x,y,t)\Psi(f(y))\,\dee y\right)\le\frac{1}{[\G(t)1](x)}\Psi^{-1}\left([\G(t)\Psi(f)](x)\right)
\end{equation}
for $(x,t)\in Q_\infty$. 
By \eqref{eq:5.4} and \eqref{eq:5.6} we obtain assertion~(1). 

We prove assertion~(2). 
For any $T\in(0,R_{**}^2]$, by Lemma~\ref{Lemma:5.1} and condition~({\bf D}')-(2)
we find $C_*>0$ such that
$$
\|k(t)\|_{L^\infty(\Omega)}\le Ct^{-\frac{N+1}{2}}\sup_{z\in\overline{\Omega}}\int_{\partial\Omega\cap B(z,\sqrt{t})}\,\dee S(y)\le C_*t^{-1},
\quad t\in(0,T). 
$$
Then it follows from Jensen's inequality and \eqref{eq:5.5} that
\begin{align*}
[\tilde{\K}(t)h](x) & \le C_*^{-1}tk(x,t)\Psi^{-1}\left(\frac{C_*^{-1}t}{C_*^{-1}tk(x,t)}\int_{\partial\Omega} K(x,y,t)\Psi(h(y))\,\dee S(y)\right)\\
 & \le \Psi^{-1}\left(C_*^{-1}t\int_{\partial\Omega} K(x,y,t)\Psi(h(y))\,\dee S(y)\right)
 =\Psi^{-1}\left([\tilde{\K}(t)\Psi(h)](x)\right)
\end{align*}
for $(x,t)\in Q_T$. Thus assertion~(2) follows, and the proof of Lemma~\ref{Lemma:5.2} is complete.
$\Box$
\subsection{Proof of Theorem~\ref{Theorem:1.4}}
For the proof of Theorem~\ref{Theorem:1.4}, we prepare the following lemma. 
\begin{lemma}
\label{Lemma:5.3}
Let $R_{**}$ satisfy condition~{\rm ({\bf D}')}. 
Then there exists $\gamma=\gamma(\Omega,p)>0$ such that,
if $\mu\in\mathcal{M}(\overline{\Omega})$ satisfies
\begin{equation}
\label{eq:5.7}
\int_0^Ts^{-\frac{N}{2}(p-1)} \left(\sup_{z\in \overline{\Omega}} \int_{B_\Omega(z,\sqrt{s})}\frac{\dee\mu(y)}{d(y)+\sqrt{s}}\right)^{p-1}\,\dee s \le \gamma
\end{equation}
for some $T\in (0,R_{**}^2]$, then problem \eqref{eq:PP} possesses a solution in~$Q_T$.
\end{lemma}
{\bf Proof.}
Assume \eqref{eq:5.7} for some $T\in (0,R_{**}^2]$. Set
$w(x,t) := 2[\K(t)\mu](x)$ for $(x,t)\in Q_T$. 
It follows from Lemma~\ref{Lemma:5.1} that 
$$
\|w(t)\|_{L^\infty(\Omega)}\le Ct^{-\frac{N}{2}}\sup_{z\in\overline{\Omega}}\,\int_{B_\Omega(z,\sqrt{t})}\frac{\dee\mu(y)}{d(y)+\sqrt{t}},
\quad t\in(0,R_{**}^2).
$$
Then, by Fubini's theorem and \eqref{eq:1.8} we have
\begin{equation*}
\begin{split}
&\int_{\overline{\Omega}} K(x,y,t) \, \dee\mu(y) + \int_0^t \int_{\Omega} G(x,y,t-s) w(y,s)^p\,\dee y\,\dee s\\
&\le \frac{1}{2}w(x,t) + \int_0^t \|w(s)\|_{L^\infty(\Omega)}^{p-1} \int_{\Omega} G(x,y,t-s) w(y,s)\,\dee y\,\dee s\\
&\le \frac{1}{2}w(x,t) + 2\int_0^t \|w(s)\|_{L^\infty(\Omega)}^{p-1} \int_{\overline{\Omega}}\int_{\Omega} G(x,y,t-s) K(y,z,s)\,\dee y\,\dee \mu(z)\,\dee s\\
&= \frac{1}{2}w(x,t) + 2\int_0^t \|w(s)\|_{L^\infty(\Omega)}^{p-1} \,\dee s\int_{\overline{\Omega}} K(x,z,t)\,\dee\mu(z)\\
&\le 
\left(\frac{1}{2}+C\int_0^T s^{-\frac{N}{2}(p-1)} \left(\sup_{z\in\overline{\Omega}} 
\int_{B_\Omega(z,\sqrt{s})}\frac{\dee \mu(y)}{d(y)+\sqrt{s}}\right)^{p-1} \,\dee s\right)w(x,t)\\
&\le \left(\frac{1}{2}+C\gamma\right)w(x,t)
\end{split}
\end{equation*}
for $(x,t) \in Q_T$. Taking $\gamma>0$ small enough if necessary, we see that $w$ is a supersolution to problem~\eqref{eq:PP} in~$Q_T$. 
Then Lemma~\ref{Lemma:5.3} follows from Lemma~\ref{Lemma:2.1}. 
$\Box$
\vspace{5pt}
\newline
{\bf Proof of Theorem~\ref{Theorem:1.4}.}
We prove assertion~(1). 
Assume that problem~\eqref{eq:PP} possesses a solution~$u$ in~$Q_T$ for some $T>0$. 
We can assume, without loss of generality, that $T<\min\{R_*^2,16\}$, where $R_*$ is as in condition~({\bf D}). 
For any $z\in\Omega(\sqrt{T}/4)$, 
we find $z_*\in\partial\Omega$ such that $B_\Omega(z,\sqrt{T}/4)\subset B_\Omega(z_*,\sqrt{T}/2)$.
Then it follows from Theorem~\ref{Theorem:1.1}-(2) and Theorem~\ref{Theorem:1.2}-(1) that 
\begin{align*}
 & \sup_{z\in\Omega(\sqrt{T}/4)}\frac{\mu(B_\Omega(z,\sqrt{T}/4))}{1+d(z)}\le
\sup_{z\in\partial\Omega}\mu(B_\Omega(z,\sqrt{T}/2))\le CT^{\frac{N+1}{2}-\frac{1}{p-1}},\\
 & \sup_{z\in\Omega^c(\sqrt{T}/4)}\frac{\mu(B_\Omega(z,\sqrt{T}/4))}{1+d(z)}
\le C\sup_{z\in\Omega^c(\sqrt{T}/4)}\frac{d(z)+\sqrt{T}}{1+d(z)}T^{\frac{N}{2}-\frac{1}{p-1}}
\le CT^{\frac{N}{2}-\frac{1}{p-1}}.
\end{align*}
These imply that
\begin{equation}
\label{eq:5.8}
\sup_{z\in\overline{\Omega}}\frac{\mu(B_\Omega(z,\sqrt{T}/4))}{1+d(z)}<\infty.
\end{equation}

Let $z\in\overline{\Omega}$. By \eqref{eq:5.2} we find a set $\{z_k\}_{k=1}^m\in{\mathbb R}^N$ such that 
\begin{equation}
\label{eq:5.9}
B_\Omega(z,1)\subset\bigcup_{k=1}^m B_\Omega(z_k,\sqrt{T}/8),
\quad
B_\Omega(z,1)\cap B_\Omega(z_k,\sqrt{T}/8)\not=\emptyset.
\end{equation}
Here there exists $m_*=m_*(N,T)>0$ such that $m\le m_*$. 
For any $k=1,\dots,m$, let $z_k^*\in B_\Omega(z_k,\sqrt{T}/8)$. 
Then
$$
B(z_k,\sqrt{T}/8)\subset B(z_k^*,\sqrt{T}/4),
\quad
d(z_k^*)\le d(z_k)+\frac{\sqrt{T}}{8}\le d(z)+1+\frac{\sqrt{T}}{4}\le d(z)+2.
$$
These together with \eqref{eq:5.8} and \eqref{eq:5.9} imply that
$$
\mu(B_\Omega(z,1))\le\sum_{k=1}^m\mu(B(z_k^*,\sqrt{T}/4))\le C\sum_{k=1}^m (1+d(z_k^*))
\le C(1+d(z))
$$
for $z\in\overline{\Omega}$. Thus \eqref{eq:1.16} holds. 

Next, we assume \eqref{eq:1.16}. 
It follows from \eqref{eq:1.16} that 
\begin{equation*}
\begin{array}{ll}
\displaystyle{\int_{B_\Omega(z,\sqrt{s})}\frac{\dee\mu(y)}{d(y)+\sqrt{s}}\le s^{-\frac{1}{2}}\mu(B_\Omega(z,\sqrt{s}))
\le Cs^{-\frac{1}{2}}(1+d(z))\le Cs^{-\frac{1}{2}}} & \quad\mbox{if}\quad d(z)\le 1,\vspace{7pt}\\
 \displaystyle{\int_{B_\Omega(z,\sqrt{s})}\frac{\dee\mu(y)}{d(y)+\sqrt{s}}\le \frac{\mu(B_\Omega(z,\sqrt{s}))}{d(z)}
\le 2\frac{\mu(B_\Omega(z,\sqrt{s}))}{1+d(z)}\le C} & \quad\mbox{if}\quad d(z)>1,
\end{array}
\end{equation*}
for $z\in\overline{\Omega}$ and $s\in(0,1)$. 
Then, since $1<p<p_{N+1}$, 
we obtain
$$
\int_0^Ts^{-\frac{N}{2}(p-1)} \left(\sup_{z\in\overline{\Omega}} \int_{B_\Omega(z,\sqrt{s})}\frac{\dee\mu(y)}{d(y)+\sqrt{s}}\right)^{p-1}\,\dee s
\le C\int_0^Ts^{-\frac{N+1}{2}(p-1)}\,\dee s
 \le CT^{-\frac{N+1}{2}(p-1)+1}
$$
for $T\in(0,1)$. 
By Lemma~\ref{Lemma:5.3}, taking $T\in(0,1)$ small enough if necessary, 
we see that problem~\eqref{eq:PP} possesses a solution in~$Q_T$. 
Thus assertion~(1) follows. 
Similarly, we prove assertion~(2). 
Assertion~(3) follows from Theorem~\ref{Theorem:1.1}-(2) and Theorem~\ref{Theorem:1.2}-(4). 
The proof of Theorem~\ref{Theorem:1.4} is complete. 
$\Box$
\subsection{Sufficient conditions for the existence of solutions}
In this subsection we collect propositions for the proof of Theorem~\ref{Theorem:1.5}. 
We use the same notation as in Section~5.1.
\begin{proposition}
\label{Proposition:5.1}
Let $f\in L^1_{{\rm loc}}(\Omega,d\,\dee x)$ and $h\in L^1_{{\rm loc}}(\partial\Omega)$ be nonnegative in $\Omega$ and on $\partial\Omega$, respectively.
Let $\Psi$ be a strictly increasing, nonnegative, and convex function in $[0,\infty)$ such that $\Psi(0)=0$.
Set
\begin{align*}
& v(x,t) := 2\Psi^{-1} ([\G(t) \Psi(f)](x)),\\
& w(x,t):= 2C_*t^{-1}\Psi^{-1} ([\tilde{\K}(t)\Psi(h)](x)),\quad (x,t)\in Q_\infty,
\end{align*}
where $C_*$ is as in Lemma~{\rm\ref{Lemma:5.2}}.
Define
\[
A(r) := \frac{\Psi^{-1}(r)^p}{r},\quad B(r) := \frac{r}{\Psi^{-1}(r)}, \quad \mbox{for} \quad r>0.
\]
Let $R_{**}$ satisfy condition~{\rm ({\bf D}')}. 
Then there exists $\gamma>0$ such that, if
\begin{equation}
\label{eq:5.10}
\begin{split}
&\sup_{t\in(0,T)} \left\{\|B(\G(t)\Psi(f))\|_{L^\infty(\Omega)} \int_0^t \|A(\G(s)\Psi(f))\|_{L^\infty(\Omega)}\,\dee s\right\} \le \gamma,\\	
&\sup_{t\in(0,T)} \left\{\|B(\tilde{\K}(t)\Psi(h))\|_{L^\infty(\Omega)}\int_0^t s^{-(p-1)} \|A(\tilde{\K}(s)\Psi(h))\|_{L^\infty(\Omega)}\,\dee s\right\} 
\le \gamma,
\end{split}
\end{equation}
for some $T\in (0,R_{**}^2]$, then problem~\eqref{eq:PP} with $\dee \mu=fd\,\dee x+h\,\dee S\in{\mathcal M}(\overline{\Omega})$ 
possesses a solution~$u$ in~$Q_T$, with $u$ satisfying
\[
0\le u(x,t) \le v(x,t) + w(x,t), \quad (x,t)\in Q_T.
\]
\end{proposition}
{\bf Proof.}
Let $\gamma>0$ be small enough, and assume \eqref{eq:5.10}. 
We show that $v+w$ is a supersolution to problem \eqref{eq:PP} in~$Q_T$.
By Lemma~\ref{Lemma:5.2}
we have 
\begin{equation}
\label{eq:5.11}
\begin{split}
 & \int_{\overline{\Omega}} K(x,y,t)\,\dee\mu(y)
 =[\G(t)f](x)+C_*t^{-1}[\tilde{\K}(t)h](x)\\
&\le \Psi^{-1}([\G(t)\Psi(f)](x)) + C_*t^{-1}\Psi^{-1}([\tilde{\K}(t)\Psi(h)](x))
= \frac{v(x,t)+w(x,t)}{2}
\end{split}
\end{equation}
for $(x,t)\in Q_\infty$.
Furthermore, 
\begin{equation}
\label{eq:5.12}
\begin{split}
 & \int_0^t\G(t-s)v(s)^p\,\dee s
 \le 2^p\int_0^t \G(t-s)\left[\Psi^{-1}(\G(s) \Psi(f))\right]^p\,\dee s\\
& \le 2^p \int_0^t \G(t-s) \left\|\frac{[\Psi^{-1}(\G(s)\Psi(f))]^p}{\G(s)\Psi(f)}\right\|_{L^\infty(\Omega)} \G(s)\Psi(f)\,\dee s \\
&= 2^p \G(t)\Psi(f)\int_0^t \|A(\G(s)\Psi(f))\|_{L^\infty(\Omega)} \,\dee s\\
&\le 2^{p-1}v(x,t) \|B(\G(t)\Psi(f))\|_{L^\infty(\Omega)}\int_0^t \|A(\G(s)\Psi(f))\|_{L^\infty(\Omega)} \,\dee s
\le C\gamma v(x,t)
\end{split}
\end{equation}
for $(x,t)\in Q_T$. 
Thanks to \eqref{eq:1.8}, 
similarly, we have 
\begin{equation}
\label{eq:5.13}
\begin{split}
&\int_0^t \G(t-s)w(s)^p\,\dee s
\le (2C_*)^p \int_0^t s^{-p}\G(t-s)  \left[\Psi^{-1}(\tilde{\K}(s)\Psi(h))\right]^p\,\dee s\\
& \le C \int_0^t s^{-p} \left\| \frac{[\Psi^{-1}(\tilde{\K}(s)\Psi(h))]^p}{\tilde{\K}(s)\Psi(h)}\right\|_{L^\infty(\Omega)}\G(t-s)\tilde{\K}(s)\Psi(h) \,\dee s\\
& =C \int_0^t s^{-(p-1)} \|A(\tilde{\K}(s)\Psi(h))\|_{L^\infty(\Omega)}\G(t-s)\K(s)\Psi(h) \,\dee s\\
& =C\K(t)\Psi(h)\int_0^t s^{-(p-1)}\|A(\tilde{\K}(s)\Psi(h))\|_{L^\infty(\Omega)}\,\dee s\\
& \le Ct^{-1}\left\|\frac{\tilde{\K}(t)\Psi(h)}{\Psi^{-1}(\tilde{\K}(t)\Psi(h))}\right\|_{L^\infty(\Omega)}
\Psi^{-1}(\tilde{\K}(t)\Psi(h))
\int_0^t s^{-(p-1)} \|A(\tilde{\K}(s)\Psi(h))\|_{L^\infty(\Omega)}\,\dee s\\
& \le C w(x,t)\|B(\tilde{\K}(t)\Psi(h))\|_{L^\infty(\Omega)}
\int_0^t s^{-(p-1)} \|A(\tilde{\K}(s)\Psi(h))\|_{L^\infty(\Omega)}\,\dee s
\le C\gamma w(x,t)
\end{split}
\end{equation}
for $(x,t)\in Q_T$. 
Since $(a+b)^p\le 2^{p-1}(a^p+b^p)$ for $a$, $b\in[0,\infty)$
and $\gamma>0$ is small enough, 
combining \eqref{eq:5.11}, \eqref{eq:5.12}, and \eqref{eq:5.13}, 
we obtain 
\begin{equation*}
\begin{split}
&\int_{\overline{\Omega}} K(x,y,t) \,\dee\mu(y) +\int_0^t\G(t-s)(v(s)+w(s))^p\,\dee s\\
&\le  \frac{v(x,t)+w(x,t)}{2} +2^{p-1}C\gamma[v(x,t)+w(x,t)]
\le v(x,t)+w(x,t)
\end{split}
\end{equation*}
for $(x,t)\in Q_T$. Thus $v+w$ is a super solution to problem~\eqref{eq:PP} in~$Q_T$. 
Then, by Lemma~\ref{Lemma:2.1} we see that 
problem~\eqref{eq:PP} possesses our desired solution in~$Q_T$, 
and the proof is complete.
$\Box$\vspace{5pt}

We apply Proposition~\ref{Proposition:5.1} to obtain the following three propositions 
on the existence of solutions to problem~\eqref{eq:PP}.
\begin{proposition}
\label{Proposition:5.2}
Let $f\in L^1_{{\rm loc}}(\overline{\Omega})$ be nonnegative in $\Omega$. 
Let $h\in L^1_{{\rm loc}}(\partial\Omega)$ be nonnegative if $1<p<2$ and $h=0$ on $\partial\Omega$ if $p\ge 2$. 
Let $R_{**}$ satisfy condition~{\rm ({\bf D}')}. 
Then, for any $\alpha>1$, there exists $\gamma = \gamma(\Omega, p,\alpha)>0$ with the following property:
if there exists $T\in (0,R_{**}^2]$ such that
\begin{equation}
\label{eq:5.14}
\begin{split}
&\sup_{z\in \overline{\Omega}} \int_{B_\Omega(z,\sigma)} \frac{d(y)}{d(y)+\sigma}f(y)^\alpha \,\dee y \le \gamma \sigma^{N-\frac{2\alpha}{p-1}},\\
&\sup_{z\in \partial \Omega} \int_{B_\Omega (z,\sigma) \cap \partial\Omega} h(y)^\alpha \,\dee S(y) \le \gamma \sigma^{N-1-2\alpha\frac{2-p}{p-1}},
\end{split}
\end{equation}
for $\sigma \in (0,\sqrt{T})$, then problem \eqref{eq:PP} 
with $\dee\mu=fd\,\dee x+h\,\dee S$
possesses a solution~$u$ in~$Q_T$,
with~$u$ satisfying
\begin{equation}
\label{eq:5.15}
0\le u(x,t) \le 2[\G(t) f^\alpha](x)^\frac{1}{\alpha}+2C_*t^{-1} [\tilde{\K}(t)h^\alpha](x)^\frac{1}{\alpha},
\quad (x,t)\in Q_T.
\end{equation}
\end{proposition}
{\bf Proof.}
Assume \eqref{eq:5.14}. We can assume, without loss of generality, that $\alpha\in (1,p)$.
Indeed, if $\alpha\ge p$, let $1<\beta<p$. 
Then it follows from H\"{o}lder's inequality and condition~({\bf D'})-(2) that
\begin{equation*}
\begin{split}
 & \sup_{z\in \overline{\Omega}} \int_{B_\Omega(z,\sigma)} \frac{d(y)}{d(y)+\sigma}f(y)^\beta \, \dee y\\
&\le \sup_{z\in\overline{\Omega} } \left[\int_{B_\Omega(z,\sigma)} \frac{d(y)}{d(y)+\sigma} \,\dee y\right]^{1-\frac{\beta}{\alpha}} 
\left[\int_{B_\Omega(z,\sigma)} \frac{d(y)}{d(y)+\sigma} f(y)^\alpha\, \dee y\right]^\frac{\beta}{\alpha}\\
& \le C\sigma^{N\left(1-\frac{\beta}{\alpha}\right)}\left(\gamma\sigma^{N-\frac{2\alpha}{p-1}}\right)^{\frac{\beta}{\alpha}}
\le C\gamma^{\frac{\beta}{\alpha}} \sigma^{N-\frac{2\beta}{p-1}},\\
& \sup_{z\in \partial \Omega} \int_{B_\Omega(z,\sigma)\cap \partial \Omega} h(y)^\beta\, \dee S(y)\\
&\le \sup_{z\in \partial \Omega}\left\{\left[\int_{B_\Omega(z,\sigma)\cap \partial \Omega}\,\dee S(y)\right]^{1-\frac{\beta}{\alpha}} 
 \left[\int_{B_\Omega(z,\sigma)\cap \partial \Omega} h(y)^\alpha\,\dee S(y)\right]^\frac{\beta}{\alpha}\right\}\\
& \le C\sigma^{(N-1)\left(1-\frac{\beta}{\alpha}\right)}\left(\gamma \sigma^{N-1-2\alpha\frac{2-p}{p-1}}\right)^{\frac{\beta}{\alpha}}
\le C\gamma^\frac{\beta}{\alpha} \sigma^{N-1-2\beta\frac{2-p}{p-1}}\quad\mbox{if}\quad 1<p<2,
\end{split} 
\end{equation*}
for $\sigma\in (0,\sqrt{T})$.
Then \eqref{eq:5.14} holds with $\alpha$ replaced by $\beta$.
Furthermore, if \eqref{eq:5.15} holds with $\alpha$ replaced by $\beta$, then
it follows from Lemma~\ref{Lemma:5.2} with $\Psi(r)=r^{\alpha/\beta}$ that
\[
[\G(t)f^\beta](x)^{\frac{1}{\beta}} \le [\G(t)f^\alpha](x)^\frac{1}{\alpha}, \quad [\tilde{\K}(t)h^{\beta}](x)^\frac{1}{\beta}\le  [\tilde{\K}(t)h^{\alpha}](x)^\frac{1}{\alpha},
\quad (x,t)\in Q_T. 
\]
This implies that \eqref{eq:5.15} holds. 
Thus it suffices to consider the case of $\alpha\in(1,p)$. 

We apply Proposition~\ref{Proposition:5.1} with $\Psi(r)=r^\alpha$ to prove Proposition~\ref{Proposition:5.2}.
Let $A$ and $B$ be as in Proposition~\ref{Proposition:5.1}.
Then $A(r) =r^{(p/\alpha)-1}$ and $B(r) =r^{1-(1/\alpha)}$.
It follows from Lemma~\ref{Lemma:5.1} and \eqref{eq:5.14} that
\begin{equation*}
\begin{split}
[\G(t)f^\alpha](x) 
 &\le C t^{-\frac{N}{2}} \sup_{z\in \overline{\Omega}} \int_{B_\Omega(z,\sqrt{t})}\frac{d(y)}{d(y)+\sqrt{t}}f(y)^\alpha\,\dee y \le C\gamma t^{-\frac{\alpha}{p-1}},\\
[\tilde{\K}(t)h^\alpha](x) 
&=C_*^{-1} t\int_{\partial\Omega} K(x,y,t) h(y)^\alpha\, \dee S(y)\\
&\le Ct^{-\frac{N}{2}+1} \sup_{z\in \partial \Omega} \int_{B_\Omega(z,\sqrt{t})\cap \partial \Omega} \frac{h(y)^\alpha}{d(y)+\sqrt{t}}\,\dee S(y)
\le C\gamma t^{-\alpha\frac{2-p}{p-1}}\quad\mbox{if}\quad 1<p<2,
 \end{split}
\end{equation*}
for $(x,t)\in Q_T$. 
Then 
\begin{align*}
 & \sup_{t\in(0,T)} \left\{\|B(\G(t)\Psi(f))\|_{L^\infty(\Omega)} \int_0^t \|A(\G(s)\Psi(f))\|_{L^\infty(\Omega)}\,\dee s\right\}\\
 & \le \sup_{t\in(0,T)}\left\{\left[C\gamma t^{-\frac{\alpha}{p-1}}\right]^{1-\frac{1}{\alpha}}
\int_0^t  \left[C\gamma s^{-\frac{\alpha}{p-1}}\right]^{\frac{p}{\alpha}-1}\,\dee s\right\}\le C\gamma^{\frac{p-1}{\alpha}}
\end{align*}
and 
\begin{align*}
 & \sup_{t\in(0,T)} \left\{\|B(\tilde{\K}(t)\Psi(h))\|_{L^\infty(\Omega)}\int_0^t s^{-(p-1)} \|A(\tilde{\K}(s)\Psi(h))\|_{L^\infty(\Omega)}\,\dee s\right\}\\
 & \le\sup_{t\in(0,T)}\left\{\left[C\gamma t^{-\alpha\frac{2-p}{p-1}}\right]^{1-\frac{1}{\alpha}}
 \int_0^t s^{-(p-1)}\left[C\gamma s^{-\alpha\frac{2-p}{p-1}}\right]^{\frac{p}{\alpha}-1}\,\dee s\right\}
 \le C\gamma^{\frac{p-1}{\alpha}}\quad\mbox{if}\quad 1<p<2.
\end{align*}
Then we apply Proposition~\ref{Proposition:5.1} to obtain the desired conclusion. Thus Proposition~\ref{Proposition:5.2} follows.~$\Box$
\begin{proposition}
\label{Proposition:5.3}
Let $p=p_{N+\ell}$ with $\ell\in\{0,1\}$ and $\beta>0$. Set $\Psi(r) :=r[\log(e+r)]^\beta$ for $r\in[0,\infty)$. 
Let $R_{**}$ satisfy condition~{\rm ({\bf D}')}. 
Then, for any $T\in(0,R_{**}^2)$, there exists $\gamma=\gamma(\Omega,p,\beta,T)>0$ such that, 
if a nonnegative measurable function $f$ in $\Omega$ satisfies
\begin{equation*}
\sup_{z\in \overline{\Omega}} \int_{B_\Omega(z,\sigma)} d(y)^\ell \Psi(T^\frac{1}{p-1} f(y)) \,\dee y \le 
\gamma T^{\frac{N+\ell}{2}} \left[\log\left(e+\frac{\sqrt{T}}{\sigma}\right)\right]^{\beta-\frac{N+\ell}{2}}
\end{equation*}
for $\sigma\in (0,\sqrt{T})$,
then problem \eqref{eq:PP} with $\dee\mu =fd\,\dee x$ possesses a solution $u$ in~$Q_T$, with $u$ satisfying
\[
0\le u(x,t) \le C\Psi^{-1}\left([\G(t)\Psi(T^\frac{1}{p-1}f)](x)\right),\quad (x,t)\in Q_T,
\]
for some $C>0$.
\end{proposition} 
{\bf Proof.}
Let $0<\epsilon <p-1$. We  find $L\in[e,\infty)$ with the following properties:
\begin{itemize}
\item[(a)] the function $\Psi_L:[0,\infty)\ni r\mapsto  r[\log(L + r)]^\beta$ is convex;
\item[(b)] the function $[0,\infty)\ni r\mapsto r^p/\Psi_L(r)$ is increasing;
\item[(c)] the function $[0,\infty)\ni r\mapsto r^\epsilon [\log(L+r)]^{-p\beta}$ is increasing.
\end{itemize} 
Since $\Psi(s)\asymp \Psi_L(s)$ for $s\in (0,\infty)$, 
it suffices to prove Proposition~\ref{Proposition:5.3} with $\Psi$ replaced by~$\Psi_L$. 

Assume that 
\begin{equation}
\label{eq:5.16}
\sup_{z\in \overline{\Omega}} \int_{B_\Omega(z,\sigma)} d(y)^\ell \Psi_L(T^\frac{1}{p-1} f(y)) \,\dee y 
\le \gamma T^{\frac{N+\ell}{2}} \left[\log\left(e+\frac{\sqrt{T}}{\sigma}\right)\right]^{\beta-\frac{N+\ell}{2}}
\end{equation}
for $\sigma\in (0,\sqrt{T})$. Here we can assume, without loss of generality, that $\gamma\in(0,1)$.
Setting
$$
\tilde{\Psi}(r):=\Psi_L(T^{\frac{1}{p-1}}r)=\Psi_L(T^{\frac{N+\ell}{2}}r),\quad r\in[0,\infty), 
$$
we apply Proposition~\ref{Proposition:5.1} with $\Psi=\tilde{\Psi}$ to prove Proposition~\ref{Proposition:5.3}. 
It follows from Lemma~\ref{Lemma:5.1} and \eqref{eq:5.16} that 
\begin{equation}
\label{eq:5.17}
\begin{split}
\left\|\G(t)\Psi_L\left(T^{\frac{1}{p-1}}f\right)\right\|_{L^\infty(\Omega)}
 & \le Ct^{-\frac{N}{2}}\sup_{z\in\overline{\Omega}}\int_{B_\Omega(z,\sqrt{t})}\frac{d(y)}{d(y)+\sqrt{t}}\Psi_L(T^{\frac{1}{p-1}}f(y))\,\dee y\\
 & \le Ct^{-\frac{N+\ell}{2}}\sup_{z\in\overline{\Omega}}\int_{B_\Omega(z,\sqrt{t})} d(y)^\ell \Psi_L(T^{\frac{1}{p-1}}f(y))\,\dee y\\
 & \le C\gamma {t_T}^{-\frac{N+\ell}{2} } |\log t_T|^{\beta-\frac{N+\ell}{2}}
\le C{t_T}^{-\frac{N+\ell}{2} } |\log t_T|^{\beta-\frac{N+\ell}{2}}
\end{split}
\end{equation}
for $t\in (0,T)$, where $t_T := t/(2T) \in (0,1/2)$.

Set 
\begin{align*}
 & A(r):= \frac{\tilde{\Psi}^{-1}(r)^p}{r}=\frac{(T^{-\frac{1}{p-1}}\Psi_L^{-1}(r))^p}{r}=T^{-\frac{p}{p-1}}\frac{\Psi_L^{-1}(r)^p}{r}
\asymp T^{-\frac{p}{p-1}}r^{p-1}[\log(L+r)]^{-p\beta},\\
 & B(r):=\frac{r}{\tilde{\Psi}^{-1}(r)}=T^{\frac{1}{p-1}}\frac{r}{\Psi_L^{-1}(r)}\asymp T^{\frac{1}{p-1}}[\log(L+r)]^\beta,
\end{align*}
for $r\in[0,\infty)$. 
Here we used the relation that $\Psi_L^{-1}(r)\asymp  r[\log(L+r)]^{-\beta}$ for $r\in[0,\infty)$.
Then, thanks to properties (b) and (c), by \eqref{eq:5.17} we have
\begin{align*}
 & \|A(\G(s)\tilde{\Psi}(f))\|_{L^\infty(\Omega)}\le A\left(\|\G(s)\tilde{\Psi}(f)\|_{L^\infty(\Omega)}\right)\\
 & \le CT^{-\frac{p}{p-1}}\|\G(s)\tilde{\Psi}(f)\|_{L^\infty(\Omega)}^{p-1-\epsilon}\|\G(s)\tilde{\Psi}(f)\|_{L^\infty(\Omega)}^\epsilon\left[\log\left(L+\|\G(s)\tilde{\Psi}(f)\|_{L^\infty(\Omega)}\right)\right]^{-p\beta}\\
 & \le CT^{-\frac{p}{p-1}}
 \left(C\gamma {s_T}^{-\frac{N+\ell}{2}} |\log s_T|^{\beta-\frac{N+\ell}{2}}\right)^{p-1-\epsilon}\\
 & \qquad
 \times\left(C{s_T}^{-\frac{N+\ell}{2} } |\log s_T|^{\beta-\frac{N+\ell}{2}}\right)^{\epsilon}
 \left[\log\left(L+C{s_T}^{-\frac{N+\ell}{2} } |\log s_T|^{\beta-\frac{N+\ell}{2}}\right)\right]^{-p\beta}\\
  & \le C\gamma^{p-1-\epsilon}T^{-\frac{p}{p-1}}(s_T)^{-\frac{N+\ell}{2}(p-1)}|\log s_T|^{\beta(p-1)-\frac{N+\ell}{2}(p-1)-p\beta}\\
  & =C\gamma^{p-1-\epsilon}T^{-\frac{p}{p-1}}(s_T)^{-1}|\log s_T|^{-1-\beta}
\end{align*}
for $s\in(0,T)$, where $s_T=s/2T$.
Similarly, we have 
\begin{align*}
\|B(\G(t)\tilde{\Psi}(f))\|_{L^\infty(\Omega)} & \le CT^{\frac{1}{p-1}}\left\|\left[\log\left(L+\G(t)\tilde{\Psi}(f)\right)\right]^\beta\right\|_{L^\infty(\Omega)}\\
 & \le CT^{\frac{1}{p-1}}\left[\log\left(L+\|\G(t)\tilde{\Psi}(f)\|_{L^\infty(\Omega)}\right)\right]^\beta\\
 & \le CT^{\frac{1}{p-1}}\left[\log\left(L+C{t_T}^{-\frac{N+\ell}{2} } |\log t_T|^{\beta-\frac{N+\ell}{2}}\right)\right]^\beta
 \le CT^{\frac{1}{p-1}}\left|\log t_T\right|^\beta
\end{align*}
for $t\in(0,T)$.
These imply that 
\begin{align*}
 & \sup_{t\in(0,T)} \left\{\|B(\G(t)\Psi(f))\|_{L^\infty(\Omega)} \int_0^t \|A(\G(s)\Psi(f))\|_{L^\infty(\Omega)}\,\dee s\right\}\\
 & \le C\gamma^{p-1-\epsilon}T^{-1}
\sup_{t\in(0,T)}\left\{\left|\log t_T\right|^\beta
\int_0^t(s_T)^{-1}|\log s_T|^{-1-\beta}\,\dee s\right\}\\
 & =C\gamma^{p-1-\epsilon}
 \sup_{t\in(0,T)}\left\{\left|\log t_T\right|^\beta
 \int_0^{t_T} s|\log s|^{-1-\beta}\,\dee s\right\}
 \le C\gamma^{p-1-\epsilon}.
\end{align*}
Therefore, if $\gamma>0$ is small enough, then, 
by Proposition~\ref{Proposition:5.1} we find our desired solution to problem~\eqref{eq:PP} in~$Q_T$. 
Thus Proposition~\ref{Proposition:5.3} follows.
$\Box$
\begin{proposition}
\label{Proposition:5.4}
Let $p=p_{N+1}<2$ and $\beta>0$. 
Set $\Psi(r):=r[\log (e+r)]^\beta$ for $\tau\ge0$.
Let $R_{**}$ satisfy condition~{\rm ({\bf D}')}. 
For any $T\in(0,R_{**}^2)$, there exists $\gamma= \gamma(\Omega, \beta,T)>0$ such that, 
if a locally integrable nonnegative function $h$ on $\partial\Omega$ satisfies
\begin{equation}
\label{eq:5.18}
\sup_{z\in \partial \Omega} \int_{B_\Omega (z,\sigma)\cap \partial\Omega} 
\Psi \left(T^\frac{1}{p-1} h(y)\right) \, \dee S(y) \le \gamma T^{\frac{N-1}{2}} \left[\log\left(e+\frac{\sqrt{T}}{\sigma}\right)\right]^{\beta-\frac{N+1}{2}}
\end{equation}
for $\sigma\in (0,\sqrt{T})$, then problem \eqref{eq:PP} with $\dee\mu = h\,\dee S$ possesses a solution $u$ in~$Q_T$, with~$u$ satisfying
\[
0\le u(x,t) \le Ct^{-1} \Psi^{-1} ([\tilde{\K}(t)\Psi(h)](x)),\quad (x,t)\in Q_T,
\]
for some $C>0$.
\end{proposition}
{\bf Proof.}
Assume \eqref{eq:5.18}. We can assume, without loss of generality, that $\gamma\in(0,1)$.
Let $\Psi_L$ and $\tilde{\Psi}$ be as in the proof of Proposition~\ref{Proposition:5.3}. Let $\epsilon\in(0,p-1)$. 
It follows from Lemma~\ref{Lemma:5.1} and \eqref{eq:5.18} that 
\begin{align*}
\left\|\tilde{\K}(t)\Psi_L\left(T^{\frac{1}{p-1}}h\right)\right\|_{L^\infty(\Omega)}
 & \le Ct^{-\frac{N}{2}+1}\sup_{z\in\partial\Omega}\int_{\partial\Omega\cap B(z,\sqrt{t})}\frac{\Psi_L\left(T^{\frac{1}{p-1}}h(y)\right)}{d(y)+\sqrt{t}}\,\dee S(y)\\
 & \le C\gamma t_T^{-\frac{N-1}{2}}|\log t_T|^{\beta-\frac{N+1}{2}}
 \le Ct_T^{-\frac{N-1}{2}}|\log t_T|^{\beta-\frac{N+1}{2}}
\end{align*}
for $t\in(0,T)$, where $t_T=t/2T$. 
Then
\begin{align*}
 & \|A(\tilde{\K}(s)\Psi(h))\|_{L^\infty(\Omega)}\le A(\|\tilde{\K}(s)\Psi(h)\|_{L^\infty(\Omega)})\\
 & \le CT^{-\frac{p}{p-1}}\|\tilde{\K}(s)\tilde{\Psi}(h)\|_{L^\infty(\Omega)}^{p-1-\epsilon}
 \|\tilde{\K}(s)\tilde{\Psi}(h)\|_{L^\infty(\Omega)}^\epsilon\left[\log\left(L+\|\tilde{\K}(s)\tilde{\Psi}(h)\|_{L^\infty(\Omega)}\right)\right]^{-p\beta}\\
 & \le CT^{-\frac{p}{p-1}}\left(C\gamma s_T^{-\frac{N-1}{2}}|\log s_T|^{\beta-\frac{N+1}{2}}\right)^{p-1-\epsilon}\\
 & \qquad
 \times\left(Cs_T^{-\frac{N-1}{2}}|\log t_T|^{\beta-\frac{N+1}{2}}\right)^\epsilon
 \left[\log\left(L+Cs_T^{-\frac{N-1}{2}}|\log s_T|^{\beta-\frac{N+1}{2}}\right)\right]^{-p\beta}\\
 & \le C\gamma^{p-1-\epsilon}T^{-\frac{p}{p-1}}(s_T)^{-\frac{N-1}{2}(p-1)}|\log s_T|^{\left(\beta-\frac{N+1}{2}\right)(p-1)-p\beta}\\
 & =C\gamma^{p-1-\epsilon}T^{-\frac{p}{p-1}}s_T^{-\frac{N-1}{2}(p-1)}|\log s_T|^{-1-\beta}
\end{align*}
for $s\in(0,T)$, where $s_T=s/2T$. 
Similarly, we have 
\begin{align*}
\|B(\tilde{\K}(t)\tilde{\Psi}(f))\|_{L^\infty(\Omega)} & \le CT^{\frac{1}{p-1}}\left\|\left[\log\left(L+\tilde{\K}(t)\tilde{\Psi}(h)\right)\right]^\beta\right\|_{L^\infty(\Omega)}\\
 & \le CT^{\frac{1}{p-1}}\left[\log\left(L+\|\tilde{\K}(t)\tilde{\Psi}(h)\|_{L^\infty(\Omega)}\right)\right]^\beta\\
 & \le CT^{\frac{1}{p-1}}\left[\log\left(L+Ct_T^{-\frac{N-1}{2}}|\log t_T|^{\beta-\frac{N+1}{2}}\right)\right]^\beta
 \le CT^{\frac{1}{p-1}}\left|\log t_T\right|^\beta
\end{align*}
for $t\in(0,T)$.
These imply that
\begin{align*}
 & \sup_{t\in(0,T)} \left\{\|B(\tilde{\K}(t)\Psi(h))\|_{L^\infty(\Omega)}\int_0^t s^{-(p-1)} \|A(\tilde{\K}(s)\Psi(h))\|_{L^\infty(\Omega)}\,\dee s\right\}\\
 & \le C\gamma^{p-1-\epsilon}T^{-1}\sup_{t\in(0,T)}\left\{\left|\log t_T\right|^\beta\int_0^t s^{-(p-1)} s_T^{-\frac{N-1}{2}(p-1)}|\log s_T|^{-1-\beta}\,\dee s\right\}\\
 & \le C\gamma^{p-1-\epsilon}T^{-(p-1)}\sup_{t\in(0,T)}\left\{\left|\log t_T\right|^\beta\int_0^{t_T} s^{-1}|\log s|^{-1-\beta}\,\dee s\right\}
 \le C\gamma^{p-1-\epsilon}T^{-(p-1)}.
\end{align*}
Therefore, if $\gamma>0$ is small enough, then, 
by Proposition~\ref{Proposition:5.1} we find our desired solution to problem~\eqref{eq:PP} in~$Q_T$. 
Thus Proposition~\ref{Proposition:5.4} follows.
$\Box$
\subsection{Proof of Theorem~\ref{Theorem:1.5}}
We apply Theorem~\ref{Theorem:1.1}, Theorem~\ref{Theorem:1.2}, and propositions in Section~5.3 to prove Theorem~\ref{Theorem:1.5}. 
\vspace{5pt}
\newline
{\bf Proof of Theorem~\ref{Theorem:1.5}-(1).}
Let $z\in\Omega$ and $p\ge p_N$. 
Set 
$$
  f_z(x)=
  \left\{
  \begin{array}{ll}
  |x-z|^{-\frac{2}{p-1}}\chi_{B(z,1)} & \mbox{if}\quad p>p_N,\vspace{5pt}\\
  \displaystyle{|x-z|^{-N}\left[\log\left(e+\frac{1}{|x-z|}\right)\right]^{-\frac{N}{2}-1}\chi_{B_\Omega(z,1)}} & \mbox{if}\quad p=p_N.
  \end{array}
  \right.
$$
If $p>p_N$, then we find $\alpha>1$ such that
\begin{equation*}
\sup_{x\in \overline{\Omega}} \int_{B_\Omega(x,\sigma)} \frac{d(y)}{d(y)+\sigma} (\kappa f_z(y))^\alpha \,\dee y
\le \kappa^\alpha \sup_{x\in \overline{\Omega}}\int_{B(x,\sigma)} |y-z|^{-\frac{2\alpha}{p-1}} \, \dee y \le C\kappa^\alpha \sigma^{N-\frac{2\alpha}{p-1}}
\end{equation*}
for $\sigma\in(0,1)$. 
Indeed, the last relation in the above inequality follows from 
\begin{equation}
\label{eq:5.19}
\begin{split}
 & \sup_{x\in B(z,2\sigma)}\int_{B(x,\sigma)} |y-z|^{-\frac{2\alpha}{p-1}} \, \dee y
\le\int_{B(z,3\sigma)} |y-z|^{-\frac{2\alpha}{p-1}} \, \dee y\le C\sigma^{N-\frac{2\alpha}{p-1}},\\
 & \sup_{x\in {\mathbb R}^N\setminus B(z,2\sigma)}\int_{B(x,\sigma)} |y-z|^{-\frac{2\alpha}{p-1}} \, \dee y
\le C\sup_{x\in {\mathbb R}^N\setminus B(z,2\sigma)}\int_{B(z,\sigma)} \sigma^{-\frac{2\alpha}{p-1}} \, \dee y\le C\sigma^{N-\frac{2\alpha}{p-1}}.
\end{split}
\end{equation}
Similarly, 
if $p=p_N$, then, for any $\beta\in (0, N/2)$, we have
\begin{equation*}
\begin{split}
&\sup_{x\in \overline{\Omega}} \int_{B_\Omega(x,\sigma)} \kappa f_z(y)[\log(e+ \kappa f_z(y))]^\beta \, \dee y\\
&\le C\kappa \sup_{x\in \overline{\Omega}}\int_{B(z,\sigma)}|y-z|^{-N}|\log|y-z||^{-\frac{N}{2}-1+\beta}\, \dee y
\le C\kappa |\log\sigma|^{-\frac{N}{2}+\beta}\\
\end{split}
\end{equation*}
for $\sigma>0$ small enough and $\kappa\in(0,1)$.  
Then we apply Proposition~\ref{Proposition:5.2} in the case of $p>p_N$ 
and Proposition~\ref{Proposition:5.3} with $\ell=0$ in the case of $p=p_N$ to see that 
problem~\eqref{eq:PP} possesses a local-in-time solution for $\kappa>0$ small enough. 

Next, we assume that problem \eqref{eq:PP} possesses a solution in~$Q_T$. 
By Theorem~\ref{Theorem:1.1}-(2) and Theorem~\ref{Theorem:1.2} 
we have 
\begin{equation}
\label{eq:5.20}
\begin{split}
 & \kappa\int_{B_\Omega(z,\sigma)}d(y) f_z(y)\,\dee y\\
 & \le\left\{
\begin{array}{ll}
C(d(z)+\sigma)\sigma^{N-\frac{2}{p-1}} & \mbox{if}\quad p>p_N,\vspace{3pt}\\
\displaystyle{C(d(z)+\sigma)\left[\log\left(e+\frac{\min\{d(z),\sqrt{T}\}}{\sigma}\right)\right]^{-\frac{N}{2}}} & \mbox{if}\quad p=p_N,
\end{array}
\right.\\
 & \le\left\{
\begin{array}{ll}
Cd(z)\sigma^{N-\frac{2}{p-1}} & \mbox{if}\quad p>p_N,\vspace{3pt}\\
Cd(z)|\log \sigma|^{-\frac{N}{2}} & \mbox{if}\quad p=p_N,
\end{array}
\right.
\end{split}
\end{equation}
for $\sigma>0$ small enough. 
On the other hand, it follows that 
\begin{equation}
\label{eq:5.21}
\int_{B_\Omega(z,\sigma)} d(y)f_z(y) \,\dee y \ge 
\left\{
\begin{array}{ll}
	Cd(z)\sigma^{N-\frac{2}{p-1}} \quad &\mbox{if}\quad p>p_{N+1},\vspace{3pt}\\
	Cd(z) |\log \sigma|^{-\frac{N}{2}} \quad   &\mbox{if}\quad p=p_{N+1}, \vspace{3pt}
\end{array}
\right.
\end{equation}
for $\sigma>0$ small enough. 
By \eqref{eq:5.20} and \eqref{eq:5.21} we see that 
problem~\eqref{eq:PP} possesses no local-in-time solutions if $\kappa$ is large enough. 
Then, by Lemma~\ref{Lemma:2.1} we find a constant $\kappa_1(z)>0$ satisfying assertion~(1). 
Thus Theorem~\ref{Theorem:1.5}-(1) follows.
$\Box$\vspace{5pt}
\newline
{\bf Proof of Theorem~\ref{Theorem:1.5}-(2).}
Let $z\in\partial\Omega$ and $p\ge p_{N+1}$. 
Set  
  $$
  f_z(x)=
  \left\{
  \begin{array}{ll}
  |x-z|^{-\frac{2}{p-1}}\chi_{B(z,1)} & \mbox{if}\quad p>p_{N+1},\vspace{5pt}\\
  \displaystyle{|x-z|^{-N-1}\left[\log\left(e+\frac{1}{|x-z|}\right)\right]^{-\frac{N+1}{2}-1}\chi_{B_\Omega(z,1)}} & \mbox{if}\quad p=p_{N+1}.
  \end{array}
  \right.
  $$
If $p>p_{N+1}$, we find $\alpha>1$ such that
\begin{equation*}
\begin{split}
\sup_{x\in\overline{\Omega}}\int_{B_\Omega(x,\sigma)} \frac{d(y)}{d(y)+\sigma} (\kappa f_z(y))^\alpha \,\dee y
& \le C\kappa^\alpha\sigma^{-1} \sup_{x\in\overline{\Omega}}
\int_{B(x,\sigma)} |y-z|^{1-\frac{2 \alpha}{p-1}} \, \dee y
\le C\kappa^\alpha \sigma^{N-\frac{2\alpha}{p-1}}
\end{split}
\end{equation*}
for $\sigma\in(0,1)$. Here we used a similar argument to that of \eqref{eq:5.19}. 
Similarly, if $p=p_{N+1}$, for any $\beta\in (0, N/2)$, we have
\begin{equation*}
\begin{split}
&\sup_{x\in \overline{\Omega}} \int_{B_\Omega(x,\sigma)} \kappa d(y) f_z(y)[\log(e+ \kappa f_z(y))]^\beta \, \dee y\\
&\le C\kappa \sup_{x\in \overline{\Omega}}\int_{B(z,\sigma)}|y-z|^{-N}|\log|y-z||^{\beta-\frac{N+1}{2}-1} \, \dee y
\le C\kappa |\log\sigma|^{\beta-\frac{N+1}{2}}\\
\end{split}
\end{equation*}
for $\sigma>0$ small enough and $\kappa\in(0,1)$. 
Then we apply Proposition~\ref{Proposition:5.2} in the case of  $p>p_{N+1}$ 
and Proposition~\ref{Proposition:5.3} with $\ell=1$ and $T=1$ in the case of $p=p_{N+1}$ to see that 
problem~\eqref{eq:PP} possesses a local-in-time solution for $\kappa>0$ small enough. 

Next, we assume that problem \eqref{eq:PP} possesses a local-in-time solution. 
By Theorem~\ref{Theorem:1.1}-(2) and Theorem~\ref{Theorem:1.2}
we have 
\begin{equation}
\label{eq:5.22}
\kappa\int_{B_\Omega(z,\sigma)} f_z(y)d(y)\,\dee y
\le\left\{
\begin{array}{ll}
C\sigma^{N+1-\frac{2}{p-1}} & \mbox{if}\quad p>p_{N+1},\vspace{3pt}\\
\displaystyle{C\left[\log\left(e+\frac{\sqrt{T}}{\sigma}\right)\right]^{-\frac{N+1}{2}}} & \mbox{if}\quad p=p_{N+1},
\end{array}
\right.
\end{equation}
for $\sigma>0$ small enough. 
On the other hand,  
it follows that
\begin{equation}
\label{eq:5.23}
\int_{B_\Omega(z,\sigma)} f_z(y)d(y) \,\dee y \ge 
\left\{
\begin{array}{ll}
	C\sigma^{N+1-\frac{2}{p-1}} \quad &\mbox{if}\quad p>p_{N+1},\vspace{3pt}\\
	C|\log \sigma|^{-\frac{N+1}{2}} \quad   &\mbox{if}\quad p=p_{N+1}, \vspace{3pt}
\end{array}
\right.
\end{equation}
for $\sigma>0$ small enough. 
By \eqref{eq:5.22} and \eqref{eq:5.23} we see that 
problem~\eqref{eq:PP} possesses no local-in-time solutions if $\kappa$ is large enough. 
Then, by Lemma~\ref{Lemma:2.1} we find a constant $\kappa_2(z)>0$ satisfying assertion~(2). 
Thus Theorem~\ref{Theorem:1.5}-(2) follows.
$\Box$\vspace{5pt}

\noindent
{\bf Proof of Theorem~\ref{Theorem:1.5}-(3).}
Let $z\in\partial\Omega$  and $p_{N+1}\le p<2$. 
Set 
  $$
  h_z(x)=
  \left\{
  \begin{array}{ll}
  |x-z|^{-\frac{2(2-p)}{p-1}}\chi_{B(z,1)} & \mbox{if}\quad p>p_{N+1},\vspace{5pt}\\
  \displaystyle{|x-z|^{-N+1}\left[\log\left(e+\frac{1}{|x-z|}\right)\right]^{-\frac{N+1}{2}-1}\chi_{B_\Omega(z,1)}} & \mbox{if}\quad p=p_{N+1}.
  \end{array}
  \right.
  $$
Thanks to $C^1$-uniform regularity of $\Omega$, 
if $p>p_{N+1}$, we find $\alpha>1$ such that 
$$
\sup_{x\in\Omega}\int_{B_\Omega (x,\sigma) \cap \partial\Omega} (\kappa h_z(y))^\alpha \,\dee S(y)
\le \kappa^\alpha\sup_{x\in\Omega}\int_{B_\Omega (x,\sigma) \cap \partial\Omega} |x-z|^{-\alpha\frac{2(2-p)}{p-2}}\,\dee S(y)
\le C\kappa^\alpha\sigma^{N-1-2\alpha\frac{2-p}{p-1}}
$$
for $\sigma\in(0,1)$. 
Here we used a similar argument to that of \eqref{eq:5.19}. 
Similarly, if $p=p_{N+1}$, we find $\beta\in(0,N/2)$ such that 
\begin{equation*}
\begin{split}
 & \sup_{x\in\partial\Omega}\int_{\partial\Omega\cap B(x,\sigma)}\kappa h_z(y)[\log(e+\kappa h_z(y))]^\beta\,\dee S(y)\\
 & \le C\kappa\sup_{x\in\partial\Omega}\int_{\partial\Omega\cap B(x,\sigma)}|y-z|^{-N+1}|\log|y-z||^{\beta-\frac{N+1}{2}-1}\,\dee S(y)
 \le C|\log\sigma|^{\beta-\frac{N+1}{2}}
\end{split}
\end{equation*}
for $\sigma>0$ small enough and $\kappa\in(0,1)$. 
Then we apply Proposition~\ref{Proposition:5.2} in the case of $p>p_{N+1}$ 
and Proposition~\ref{Proposition:5.4} with $\ell=1$ in the case of $p=p_{N+1}$ to see that 
problem~\eqref{eq:PP} possesses a local-in-time solution for $\kappa>0$ small enough. 

Next, we assume that problem~\eqref{eq:PP} possesses a local-in-time solution. 
By Theorem~\ref{Theorem:1.1}-(2) and Theorem~\ref{Theorem:1.2} 
we have 
\begin{equation}
\label{eq:5.24}
\kappa\int_{\partial\Omega\cap B(z,\sigma)}h_z(y)\,\dee S(y)
\le\left\{
\begin{array}{ll}
C\sigma^{N+1-\frac{2}{p-1}} & \mbox{if}\quad p>p_{N+1},\vspace{3pt}\\
\displaystyle{C\left[\log\left(e+\frac{\sqrt{T}}{\sigma}\right)\right]^{-\frac{N+1}{2}}} & \mbox{if}\quad p=p_{N+1},
\end{array}
\right.
\end{equation}
for $\sigma>0$ small enough. 
On the other hand, 
we see that
\begin{equation}
\label{eq:5.25}
\int_{\partial\Omega\cap B(z,\sigma)}h_z(y)\,\dee S(y) \ge 
\left\{
\begin{array}{ll}
	C\sigma^{N+1-\frac{2}{p-1}} \quad &\mbox{if}\quad p>p_{N+1},\vspace{3pt}\\
	C|\log \sigma|^{-\frac{N+1}{2}} \quad   &\mbox{if}\quad p=p_{N+1}, \vspace{3pt}
\end{array}
\right.
\end{equation}
for $\sigma>0$ small enough. 
By \eqref{eq:5.24} and \eqref{eq:5.25} we see that 
problem~\eqref{eq:PP} possesses no local-in-time solutions if $\kappa$ is large enough. 
Then, by Lemma~\ref{Lemma:2.1} we find a constant $\kappa_3(z)>0$ satisfying assertion~(3). 
Thus Theorem~\ref{Theorem:1.5}-(3) follows.
$\Box$
\vspace{8pt}

\noindent
{\bf Acknowledgment.}
K. I. was supported in part by JSPS KAKENHI Grant Number 19H05599. 


\begin{bibdiv}
\begin{biblist}
\bib{AIS}{article}{
   author={Akagi, Goro},
   author={Ishige, Kazuhiro},
   author={Sato, Ryuichi},
   title={The Cauchy problem for the Finsler heat equation},
   journal={Adv. Calc. Var.},
   volume={13},
   date={2020},
   pages={257--278},
}
\bib{ADi}{article}{
   author={Andreucci, D.},
   author={DiBenedetto, E.},
   title={On the Cauchy problem and initial traces for a class of evolution
   equations with strongly nonlinear sources},
   journal={Ann. Scuola Norm. Sup. Pisa Cl. Sci. (4)},
   volume={18},
   date={1991},
   pages={363--441},
}
\bib{A}{article}{
   author={Aronson, D. G.},
   title={Non-negative solutions of linear parabolic equations},
   journal={Ann. Scuola Norm. Sup. Pisa Cl. Sci. (3)},
   volume={22},
   date={1968},
   pages={607--694},
}
\bib{AC}{article}{
   author={Aronson, D. G.},
   author={Caffarelli, L. A.},
   title={The initial trace of a solution of the porous medium equation},
   journal={Trans. Amer. Math. Soc.},
   volume={280},
   date={1983},
   pages={351--366},
}
\bib{BP}{article}{%
   author={Baras, Pierre},
   author={Pierre, Michel},
   title={Crit\`ere d'existence de solutions positives pour des \'{e}quations
   semi-lin\'{e}aires non monotones},
   journal={Ann. Inst. H. Poincar\'{e} Anal. Non Lin\'{e}aire},
   volume={2},
   date={1985},
   pages={185--212},
}
\bib{BCP}{article}{
   author={B\'{e}nilan, Philippe},
   author={Crandall, Michael G.},
   author={Pierre, Michel},
   title={Solutions of the porous medium equation in ${\mathbb R}^{N}$ under
   optimal conditions on initial values},
   journal={Indiana Univ. Math. J.},
   volume={33},
   date={1984},
   pages={51--87},
}
\bib{BCV}{article}{
   author={Bidaut-V\'{e}ron, Marie-Fran\c{c}oise},
   author={Chasseigne, Emmanuel},
   author={V\'{e}ron, Laurent},
   title={Initial trace of solutions of some quasilinear parabolic equations
   with absorption},
   journal={J. Funct. Anal.},
   volume={193},
   date={2002},
   pages={140--205},
}
\bib{BD}{article}{
   author={Bidaut-V\'{e}ron, Marie-Fran\c{c}oise},
   author={Dao, Nguyen Anh},
   title={Initial trace of solutions of Hamilton-Jacobi parabolic equation
   with absorption},
   journal={Adv. Nonlinear Stud.},
   volume={15},
   date={2015},
   pages={889--921},
}
\bib{BSV}{article}{
   author={Bonforte, Matteo},
   author={Sire, Yannick},
   author={V\'{a}zquez, Juan Luis},
   title={Optimal existence and uniqueness theory for the fractional heat
   equation},
   journal={Nonlinear Anal.},
   volume={153},
   date={2017},
   pages={142--168},
}
\bib{CV}{article}{
   author={Chen, Huyuan},
   author={V\'{e}ron, Laurent},
   title={Initial trace of positive solutions to fractional diffusion
   equations with absorption},
   journal={J. Funct. Anal.},
   volume={276},
   date={2019},
   pages={1145--1200},
}
\bib{Cho}{article}{
   author={Cho, Sungwon},
   title={Two-sided global estimates of the Green's function of parabolic
   equations},
   journal={Potential Anal.},
   volume={25},
   date={2006},
   pages={387--398},
   issn={0926-2601},
}
\bib{CKP}{article}{
   author={Cho, Sungwon},
   author={Kim, Panki},
   author={Park, Hyein},
   title={Two-sided estimates on Dirichlet heat kernels for time-dependent
   parabolic operators with singular drifts in $C^{1,\alpha}$-domains},
   journal={J. Differential Equations},
   volume={252},
   date={2012},
   pages={1101--1145},
}
\bib{DK}{article}{
   author={Dahlberg, Bj\"orn E. J.},
   author={Kenig, Carlos E.},
   title={Weak solutions of the porous medium equation in a cylinder},
   journal={Trans. Amer. Math. Soc.},
   volume={336},
   date={1993},
   pages={701--709},
}
\bib{DH}{article}{
   author={DiBenedetto, E.},
   author={Herrero, M. A.},
   title={On the Cauchy problem and initial traces for a degenerate
   parabolic equation},
   journal={Trans. Amer. Math. Soc.},
   volume={314},
   date={1989},
   pages={187--224},
}
\bib{DH02}{article}{
   author={DiBenedetto, E.},
   author={Herrero, M. A.},
   title={Nonnegative solutions of the evolution $p$-Laplacian equation.
   Initial traces and Cauchy problem when $1<p<2$},
   journal={Arch. Rational Mech. Anal.},
   volume={111},
   date={1990},
   pages={225--290},
}
\bib{EG}{book}{
   author={Evans, Lawrence C.},
   author={Gariepy, Ronald F.},
   title={Measure theory and fine properties of functions},
   series={Studies in Advanced Mathematics},
   publisher={CRC Press, Boca Raton, FL},
   date={1992},
   pages={viii+268},
}
\bib{Fr}{book}{
   author={Friedman, Avner},
   title={Partial differential equations of parabolic type},
   publisher={Prentice-Hall, Inc., Englewood Cliffs, NJ},
   date={1964},
   pages={xiv+347},
}
\bib{FHIL}{article}{
   author={Fujishima, Yohei},
   author={Hisa, Kotaro},
   author={Ishige, Kazuhiro},
   author={Laister, Robert},
   title={Solvability of superlinear fractional parabolic equations},
   journal={J. Evol. Equ.},
   volume={23},
   date={2023},
   pages={Paper No. 4, 38},
}
\bib{FHIL02}{article}{
   author={Fujishima, Yohei},
   author={Hisa, Kotaro},
   author={Ishige, Kazuhiro},
   author={Laister, Robert},
   title={Local solvability and dilation-critical singularities of
   supercritical fractional heat equations},
   journal={J. Math. Pures Appl.},
   volume={186},
   date={2024},
   pages={150--175},
}
\bib{FI01}{article}{
   author={Fujishima, Yohei},
   author={Ishige, Kazuhiro},
   title={Initial traces and solvability of Cauchy problem to a semilinear
   parabolic system},
   journal={J. Math. Soc. Japan},
   volume={73},
   date={2021},
   pages={1187--1219},
}
\bib{GT}{book}{
   author={Gilbarg, David},
   author={Trudinger, Neil S.},
   title={Elliptic partial differential equations of second order},
   series={Grundlehren der mathematischen Wissenschaften}, 
   volume={224},
   edition={2},
   publisher={Springer-Verlag, Berlin},
   date={1983},
   pages={xiii+513},
}
\bib{HW}{article}{
   author={Hartman, Philip},
   author={Wintner, Aurel},
   title={On the solutions of the equation of heat conduction},
   journal={Amer. J. Math.},
   volume={72},
   date={1950},
   pages={367--395},
}
\bib{HP}{article}{
   author={Herrero, Miguel A.},
   author={Pierre, Michel},
   title={The Cauchy problem for $u_t=\Delta u^m$ when $0<m<1$},
   journal={Trans. Amer. Math. Soc.},
   volume={291},
   date={1985},
   pages={145--158},
}
\bib{Hisa}{article}{%
   author={Hisa, Kotaro},
   title={Optimal singularities of initial data of a fractional semilinear heat equation in open sets},
   journal={preprint (https://doi.org/10.48550/arXiv.2312.10969)},
}
\bib{HI01}{article}{%
   author={Hisa, Kotaro},
   author={Ishige, Kazuhiro},
   title={Existence of solutions for a fractional semilinear parabolic
   equation with singular initial data},
   journal={Nonlinear Anal.},
   volume={175},
   date={2018},
   pages={108--132},
}
\bib{HI02}{article}{%
   author={Hisa, Kotaro},
   author={Ishige, Kazuhiro},
   title={Solvability of the heat equation with a nonlinear boundary
   condition},
   journal={SIAM J. Math. Anal.},
   volume={51},
   date={2019},
   pages={565--594},
}
\bib{HIT02}{article}{
   author={Hisa, Kotaro},
   author={Ishige, Kazuhiro},
   author={Takahashi, Jin},
   title={Initial traces and solvability for a semilinear heat equation on a
   half space of ${\mathbb R}^N$},
   journal={Trans. Amer. Math. Soc.},
   volume={376},
   date={2023},
   pages={5731--5773},
}
\bib{HC}{article}{
   author={Hui, Kin Ming},
   author={Chou, Kai-Seng},
   title={Nonnegative solutions of the heat equation in a cylindrical domain
   and Widder's theorem},
   journal={J. Math. Anal. Appl.},
   volume={532},
   date={2024},
   pages={Paper No. 127929, 18},
}
\bib{I}{article}{
   author={Ishige, Kazuhiro},
   title={On the existence of solutions of the Cauchy problem for a doubly
   nonlinear parabolic equation},
   journal={SIAM J. Math. Anal.},
   volume={27},
   date={1996},
   pages={1235--1260},
}
\bib{IS}{article}{
   author={Ikeda, Masahiro},
   author={Sobajima, Motohiro},
   title={Sharp upper bound for lifespan of solutions to some critical
   semilinear parabolic, dispersive and hyperbolic equations via a test
   function method},
   journal={Nonlinear Anal.},
   volume={182},
   date={2019},
   pages={57--74},
}
\bib{IS02}{article}{
   author={Ikeda, Masahiro},
   author={Sobajima, Motohiro},
   title={Remark on upper bound for lifespan of solutions to semilinear
   evolution equations in a two-dimensional exterior domain},
   journal={J. Math. Anal. Appl.},
   volume={470},
   date={2019},
   pages={318--326},
}
\bib{IKO}{article}{%
   author={Ishige, Kazuhiro},
   author={Kawakami, Tatsuki},
   author={Okabe, Shinya},
   title={Existence of solutions for a higher-order semilinear parabolic
   equation with singular initial data},
   journal={Ann. Inst. H. Poincar\'{e} Anal. Non Lin\'{e}aire},
   volume={37},
   date={2020},
   pages={1185--1209},
}
\bib{IJK}{article}{
   author={Ishige, Kazuhiro},
   author={Kinnunen, Juha},
   title={Initial trace for a doubly nonlinear parabolic equation},
   journal={J. Evol. Equ.},
   volume={11},
   date={2011},
   pages={943--957},
}
\bib{ISato}{article}{
   author={Ishige, Kazuhiro},
   author={Sato, Ryuichi},
   title={Heat equation with a nonlinear boundary condition and uniformly
   local $L^r$ spaces},
   journal={Discrete Contin. Dyn. Syst.},
   volume={36},
   date={2016},
   pages={2627--2652},
}
\bib{LSU}{book}{
   author={La\dee y\v zenskaja, O. A.},
   author={Solonnikov, V. A.},
   author={Ural\cprime ceva, N. N.},
   title={Linear and quasilinear equations of parabolic type},
   series={Translations of Mathematical Monographs},
   volume={Vol. 23},
   publisher={American Mathematical Society, Providence, RI},
   date={1968},
}
\bib{L}{article}{
   author={Lieberman, Gary M.},
   title={Boundary and initial regularity for solutions of degenerate
   parabolic equations},
   journal={Nonlinear Anal.},
   volume={20},
   date={1993},
   pages={551--569},
}
\bib{MV}{article}{
   author={Marcus, Moshe},
   author={V\'{e}ron, Laurent},
   title={Initial trace of positive solutions of some nonlinear parabolic
   equations},
   journal={Comm. Partial Differential Equations},
   volume={24},
   date={1999},
   pages={1445--1499},
}
\bib{MV02}{article}{
   author={Marcus, Moshe},
   author={V\'{e}ron, Laurent},
   title={Semilinear parabolic equations with measure boundary data and
   isolated singularities},
   journal={J. Anal. Math.},
   volume={85},
   date={2001},
   pages={245--290},
}
\bib{MV03}{article}{
   author={Marcus, Moshe},
   author={V\'{e}ron, Laurent},
   title={Initial trace of positive solutions to semilinear parabolic
   inequalities},
   journal={Adv. Nonlinear Stud.},
   volume={2},
   date={2002},
   pages={395--436},
}
\bib{MP}{article}{
   author={Mitidieri, \`E.},
   author={Pokhozhaev, S. I.},
   title={A priori estimates and the absence of solutions of nonlinear
   partial differential equations and inequalities},
   journal={Tr. Mat. Inst. Steklova},
   volume={234},
   date={2001},
   pages={1--384},
   issn={0371-9685},
}
\bib{QS}{book}{
   author={Quittner, Pavol},
   author={Souplet, Philippe},
   title={Superlinear parabolic problems},
   series={Birkh\"auser Advanced Texts: Basler Lehrb\"ucher. [Birkh\"auser
   Advanced Texts: Basel Textbooks]},
   edition={2},
   publisher={Birkh\"auser/Springer, Cham},
   date={2019},
   pages={xvi+725},
}
\bib{R}{article}{
   author={Riahi, Lotfi},
   title={Comparison of Green functions and harmonic measures for parabolic
   operators},
   journal={Potential Anal.},
   volume={23},
   date={2005},
   pages={381--402},
}
\bib{TY}{article}{
   author={Takahashi, Jin},
   author={Yamamoto, Hikaru},
   title={Solvability of a semilinear heat equation on Riemannian manifolds},
   journal={J. Evol. Equ.},
   volume={23},
   date={2023},
   pages={Paper No. 33, 55},
}
\bib{W1}{article}{
   author={Widder, D. V.},
   title={Positive temperatures on an infinite rod},
   journal={Trans. Amer. Math. Soc.},
   volume={55},
   date={1944},
   pages={85--95},
}
\bib{W2}{article}{
   author={Widder, D. V.},
   title={Positive temperatures on a semi-infinite rod},
   journal={Trans. Amer. Math. Soc.},
   volume={75},
   date={1953},
   pages={510--525},
}
\bib{Z}{article}{
   author={Zhang, Qi S.},
   title={The boundary behavior of heat kernels of Dirichlet Laplacians},
   journal={J. Differential Equations},
   volume={182},
   date={2002},
   pages={416--430},
   issn={0022-0396},
}
\bib{Zhao}{article}{
   author={Junning, Zhao},
   title={On the Cauchy problem and initial traces for the evolution
   $p$-Laplacian equations with strongly nonlinear sources},
   journal={J. Differential Equations},
   volume={121},
   date={1995},
   pages={329--383},
}
\bib{ZX}{article}{
   author={Zhao, Junning},
   author={Xu, Zhonghai},
   title={Cauchy problem and initial traces for a doubly nonlinear
   degenerate parabolic equation},
   journal={Sci. China Ser. A},
   volume={39},
   date={1996},
   pages={673--684},
}
\end{biblist}
\end{bibdiv}
\end{document}